\documentclass{article}

\usepackage{tikz}
\usetikzlibrary{positioning}
\usepackage{graphicx}
\usepackage{bm}
\usepackage{comment}
\usepackage{psfrag}
\usepackage{latexsym,amsmath,amsfonts,amscd, amsthm}
\usepackage{changebar}
\usepackage{color}
\usepackage{bm}
\usepackage{tikz}
\usepackage{multirow}
\usepackage{xcolor}
\usepackage{hyperref}
\usepackage{setspace}
\usepackage{amssymb}
\usepackage{mathrsfs}
\usepackage{caption}
\usepackage{subcaption}
\usetikzlibrary{arrows,backgrounds}

\usepackage{algcompatible}
\usepackage[ruled,vlined]{algorithm2e}

\usepackage{ulem}
\usepackage[title]{appendix}

\newtheoremstyle{plainNoItalics}{}{}{\normalfont}{}{\bfseries}{.}{ }{}

\theoremstyle{plain}
\newtheorem{thm}{Theorem}[section]

\theoremstyle{plainNoItalics}

\newtheorem{defn}[thm]{Definition}
\newtheorem{rem}[thm]{Remark}
\newtheorem{prop}[thm]{Proposition}

\newcommand{\beq}{\begin{equation}}
\newcommand{\eeq}{\end{equation}}
\newcommand{\beqa}{\begin{eqnarray}}
\newcommand{\eeqa}{\end{eqnarray}}
\newcommand{\bit}{\begin{itemize}}
\newcommand{\eit}{\end{itemize}}
\newcommand{\bedef}{\begin{defn}}
\newcommand{\edefn}{\end{defn}}
\newcommand{\bpro}{\begin{prop}}
\newcommand{\epro}{\end{prop}}

\newcommand{\Dx}{\Delta x}
\newcommand{\Dt}{\Delta t}

\newcommand{\MP}{\mathcal{P} }

\newcommand{\NN}{\mathcal{N}\mathcal{N}}
\newcommand{\phiNN}{{\phi}^{\mathcal{N}\mathcal{N}}}
\newcommand{\alphaNN}{{\alpha}^{\mathcal{N}\mathcal{N}}}
\newcommand{\NNphi}{{\mathcal{N}\mathcal{N}}^{\phi}}
\newcommand{\NNalpha}{{\mathcal{N}\mathcal{N}}^{\alpha}}

\newcommand{\half}{ \frac{1}{2} }

\newcommand{\mw}{\textcolor{black}}
\newcommand{\mwmw}{\textcolor{black}}
\newcommand{\lli}{\textcolor{black}}
\newcommand{\pzcrev}{\textcolor{black}}

\newcommand{\pzc}{\textcolor{black}}
\newcommand{\li}{\textcolor{black}}

\newcommand{\fli}{\textcolor{black}}

\newcommand{\pzcreview}{\textcolor{black}}
\newcommand{\minreview}{\textcolor{black}}

\definecolor{mygreen}{rgb}{0.0,0.65,0.5}
\definecolor{myorange}{rgb}{0.95,0.45,0.0}
\newcommand{\review}{\textcolor{black}}
\newcommand{\reviewone}{\textcolor{black}}
\newcommand{\reviewtwo}{\textcolor{black}}
\newcommand{\reviewthree}{\textcolor{black}}
\newcommand{\reviewmulti}{\textcolor{black}}
\newcommand{\pzcnew}{\textcolor{black}}
\providecommand{\keywords}[1]
{
  \small	
  \textbf{Keywords---} #1
}

\begin{document}
\title{A learning-based projection method for model order reduction of transport problems}

\author{
Zhichao Peng\thanks{Department of Mathematics, Michigan State University, East Lansing, MI 48824 U.S.A. Email: {\tt pengzhic@msu.edu}.}
\and
Min Wang\thanks{Department of Mathematics, Duke University, Durham, NC 27705 U.S.A. Email: {\tt wangmin@math.duke.edu}.}
\and
Fengyan Li\thanks{Department of Mathematical Sciences, Rensselaer Polytechnic Institute, Troy, NY 12180, U.S.A. Email: 
{\tt lif@rpi.edu}. Research is supported by NSF grants DMS-1719942 and DMS-1913072.}
}

\maketitle

\abstract{
The Kolmogorov $n$-width of the solution manifolds of transport-dominated problems can decay slowly. As a result, it can be challenging to design efficient and accurate reduced order models (ROMs) for such problems. To address this issue, we propose a new learning-based projection method to construct nonlinear adaptive ROMs for transport problems. The construction follows the offline-online decomposition. In the offline stage, we train a neural network to construct adaptive reduced basis dependent on time and model parameters.  In the online stage, we project the solution to the learned reduced manifold.
Inheriting the merits from both deep learning and the projection method, the proposed method is more 
efficient  than the conventional linear projection-based methods, and may reduce the generalization error of a solely learning-based ROM.
Unlike some learning-based projection methods, the proposed method does not need to take derivatives of the neural network in the online stage.
\\
}


\section{Introduction\label{sec:introduction}}

The transport \reviewone{phenomena arise} %arises 
in many important areas of applications,
such as fluid dynamics, plasma physics and electromagnetics. Efficient and accurate reduced order models (ROMs) of transport-dominated problems are highly desired. However, as pointed \mw{out by} \reviewone{\cite{benner2017model,ohlberger2016reduced,greif2019decay,cohen2021optimal}}, the Kolmogorov $n$-width \cite{pinkus2012n} of the solution manifold for a transport problem may decay slowly. If one applies  conventional ROM strategies such as proper orthogonal decomposition (POD) methods \cite{berkooz1993proper,gubisch2017proper} or reduced basis methods (RBM) \cite{hesthaven2016certified} straightforwardly,
 linear subspaces of rather high dimensions will be needed to obtain accurate reduced order approximations. The design of a model \mw{with genuine reduced order for} transport problems is \mw{hence} challenging.

Overcoming this difficulty 
\mw{has been an} active research topic in the ROM community, where two  \mw{general ideas} have been extensively considered:   1) to introduce online adaptivity \mw{to the reduced solution space}, and \pzc{2) to identify and conduct inherent nonlinear transformations 
which enable the construction of efficient ROMs.} 
\mw{Examples following the former include}
\cite{carlberg2015adaptive} and \cite{peherstorfer2020model}, where \cite{carlberg2015adaptive} adaptively enriches the reduced basis in the online \mw{stage} through splitting \mw{while} \cite{peherstorfer2020model} applies online adaptive sampling and basis reconstruction based on a local low-rank structure.
% The research of finding a suitable nonlinear transformation is very active. 
\mw{Examples following the latter \pzcreview{direction} include the}
shift POD method \cite{reiss2018shifted}, \mw{which} utilizes a shift operator determined by the underlying wave speed. 
Such information is also utilized under a Lagrangian framework \mw{in designing} the Lagrangian POD \cite{mojgani2017lagrangian}, the Lagrangian DMD \cite{lu2020lagrangian} and its recent extension to the case with shocks \cite{lu2021dynamic}. 
\mw{Besides, nonlinear} transformations can also be \mw{obtained} by the method of freezing \cite{rowley2003reduction,ohlberger2013nonlinear,beyn2004freezing}, the transformed snapshot interpolation \cite{welper2017h,welper2017interpolation,welper2020transformed}, the approximated Lax pairs \cite{gerbeau2014approximated,gerbeau2015reduced}, template fitting \cite{kirby1992reconstructing,rowley2000reconstruction}, the transport reversal method \cite{rim2018transport}, shock fitting \cite{taddei2015reduced}, the \reviewtwo{registration method \cite{taddei2020registration,torlo2020model,ferrero2021registration}}, and  the greedy Wasserstein barycenter algorithm \cite{ehrlacher2019nonlinear}.

Our proposed method is inspired by the series of work \cite{cagniart2017model,cagniart2018few,cagniart2019model,nonino2019overcoming}, and is also closely related to \reviewtwo{\cite{rim2019manifold,nair2019transported,torlo2020model,mojgani2020physics,taddei2020registration,ferrero2021registration}}. The key underlying idea  of \cite{cagniart2017model,cagniart2018few,cagniart2019model,nonino2019overcoming} is that, for a transport-dominated problem, there exists some suitable nonlinear transformation \pzc{$T_{t,\mu}$}, under which the transformed solution manifold displays  a much faster decay in its Kolmogorov $n$-width  and hence can be accurately approximated by a low-dimensional reduced space. Once the transformation is 
obtained by some means, a reduced order approximation is constructed in the form of 
\begin{align}
    u_r(x,t;\mu) =\sum_{i=1}^r \alpha_i(t,\mu)\psi_i( T_{t,\mu}(x;t,\mu) ) =  \sum_{i=1}^r \alpha_i(t,\mu)\phi_i(x; t, \mu).
    \label{eq:parametric_expansion}
\end{align}
Here, $x$ is the physical location, $t$ is time, $\mu$ is a scalar or vector-valued model parameter. Each $\psi_i$ is a reduced basis function in the transformed space and \pzc{$\phi_i(x;t,\mu)=\psi_i( T_{t,\mu}(x;t,\mu) )$} is the respective reduced basis function in the original solution manifold.  
The distinct feature of the 
adaptive reduced basis 
$\{\phi_i(x; t, \mu)\}_{i=1}^r$ is that each member depends on time $t$ and the model parameter $\mu$, in addition to $x$. The nonlinear transformation $T_{t,\mu}$ can be realized as coordinate transformations defined case by case,  such as for the 2D Euler equations in \cite{cagniart2017model},
the 1D Burgers equation in \cite{cagniart2019model} and the Navier-Stokes equations in  \cite{nonino2019overcoming}. 
The manifold approximations via transported subspaces (MATS) \cite{rim2019manifold}, \mw{on the other hand,} \mw{construct} the nonlinear transformation \mw{with an} optimal transport tool called displacement interpolation \cite{mccann1997convexity,villani2008optimal}, \mw{while}
the transported snapshot method \cite{nair2019transported} applies a polynomial approximation in \mw{such constructions}. \review{The Lagrangian registration based method \cite{taddei2020registration,ferrero2021registration}, the ALE registration method \cite{torlo2020model} and the ALE autoencoder \cite{mojgani2020physics} construct $T_{t,\mu}$ by solving nonlinear optimization problems.}

In this paper, we propose a new framework to design nonlinear reduced order models for transport problems by working with 
an adaptively constructed reduced basis $\{\phi_i(x; t,\mu)\}_{i=1}^r$. The framework can be decomposed into an offline stage and an online stage. In the offline stage, instead of \pzc{predetermining} \pzcreview{or  seeking} the transformation $T_{t,\mu}$, the adaptive basis will be learned directly from data
with neural networks, \review{since the ultimate goal is to obtain the reduced order space.} More specifically, $\{\alpha_i(t,\mu)\}_{i=1}^r$ and $\{\phi_i(x; t, \mu)\}_{i=1}^r$ in \eqref{eq:parametric_expansion} will be parameterized separately as feed-forward neural networks taking $t,\mu$ or $x,t,\mu$ as inputs,
and these neural networks \mw{will be} trained \mw{with} snapshot solutions (e.g. high-fidelity numerical solutions). In the online stage,  instead of making predictions entirely with \mw{the} trained neural networks as in a pure learning-based method, we only \mw{take} the learned basis $\{\phi_i(x; t,\mu)\}_{i=1}^r$  to construct a subspace \mw{that is local to $t$ and $\mu$}, and apply the standard projection method to find the reduced order approximation. 
\review{The performance of the proposed method is demonstrated numerically through a series of linear and nonlinear scalar hyperbolic conservation laws in 1D and 2D as well as the 1D  compressible Euler system. }
The methodology developed here can \review{also} be applied to obtain  reduced order approximations for other transport problems. \reviewtwo{Unlike the aforementioned methods that explicitly identify $T_{t,\mu}$ to build ROMs, our proposed method only learns the basis and the corresponding reduced order space, with the information of the transformation $T_{t,\mu}$ implicitly encoded in the learned space.  Identifying $T_{t,\mu}$ can sometimes be helpful, but the explicit form of it is not always necessary.} 

The proposed learning-based projection method is a combination of the deep learning \mw{method} and the standard  projection method. Compared with a conventional projection method  (e.g. POD), or  \pzc{a solely} learning-based method, the proposed method has several advantages \review{in constructing ROMs}. Our numerical experiments show that the proposed method is more accurate than the standard POD method for long time predictions with a similar \review{order reduction}. \mw{The reason behind is that the} inherent nonlinear information of the underlying transport problems can be implicitly embedded \review{in} the adaptive basis through training while a conventional linear projection method may suffer from the slow decay of the Kolmogorov $n$-width. Particularly, \mw{we observe that the} adaptive basis \mw{obtained in fact} travels \mw{at} the underlying wave speed for a simple 1D linear advection equation.
Another advantage over the conventional projection method is that the neural network takes $x$ as an input in a mesh-free manner, 
hence it is simple and flexible to \pzc{work with different meshes in the online and offline stages, or even learn from experimental \mw{mesh-free} data.}
%learn from data  associated with  different meshes or even from experimental data. 
\mw{Comparing} with a pure learning-based approach, our method can reduce the generalization error due to overfitting with an online projection procedure utilizing the underlying PDE information.

We also want to point out that we are not the first to combine projection methods \mw{with the} deep learning \mw{techniques} to build ROMs for transport problems. In \cite{lee2019deep,lee2020model,kim2020fast}, auto-encoders \mw{were used} to learn a low-dimensional nonlinear approximation of the solution manifold offline, and a reduced-order solution \mw{can then be obtained} through \review{an online} projection . 
The main differences between our method and \mw{aforementioned} methods are as follows. 
\mw{First,} \cite{lee2019deep,lee2020model,kim2020fast} are based on a different reduced-order representation.
Second, in the online computations, \mw{ different from our method,} these methods need to take the derivative of the decoder neural network, and this could be computationally expensive.
To improve \mw{their} online efficiency, \cite{kim2020fast} uses a decoder with a shallow architecture  and applies hyper-reduction techniques. 
\pzc{
Third, compared with these methods, our method does not require using the same spatial meshes at online and offline stages,  and \mw{is thus} advantageous for \reviewtwo{moving mesh methods or} $h$-adaptive simulations with locally refined meshes.}
Similar ideas as in \cite{lee2019deep,lee2020model,kim2020fast} are
also utilized to develop a learning based ROM \cite{fresca2021comprehensive} under a pure data-driven framework. There, a reduced order nonlinear trial manifold is constructed by the decoder part of an auto-encoder trained offline. Instead of using a projection-based method online, \cite{fresca2021comprehensive}  evolves the coordinate in the low-dimensional latent space through \pzcreview{an extra} neural network which is trained offline to learn the reduced dynamics. \reviewtwo{This data-driven online phase is more efficient than the  projection based online algorithm in \cite{lee2019deep,lee2020model,kim2020fast} and our method, but may suffer from larger generalization errors.} Another interesting work \mw{that uses} deep learning to construct ROMs for transport problems is \cite{rim2020depth}, which utilizes ideas \mw{similar} to the MATS \cite{rim2019manifold} to \mw{lighten} the architectures of neural networks.

\pzcrev{The main difference between our method and the optimal transport based method \cite{rim2019manifold}, optimization based registration/autoencoder \cite{torlo2020model,taddei2020registration,mojgani2020physics,ferrero2021registration}  and the transported snapshot method \cite{nair2019transported} is that we find the adaptive reduced order space without explicitly finding the transformation, while they explicitly find the transformation. We use a neural network ansatz for the adaptive reduced basis, while these methods use different ansatz. Our method can be seen as an alternative to these methods.}
%\li{(Maybe should be one paragraph?) }
%\pzcrev{
\li{We would also like to make more comments about the optimal transport based methods \cite{ehrlacher2019nonlinear,rim2019manifold}.} %are also promising.
\pzcrev{Compared with deep learning based methods, these optimal transport based methods are more interpretable, and theoretical understandings have been established. However, more developments are still needed for more complicated transport problems.
The greedy Wasserstein barycenter method \cite{ehrlacher2019nonlinear} and its recent generalization to the nonconservative flows in porous media \cite{battisti2022wasserstein} are limited to 1D for now, as they use an exponential map only valid in 1D to efficiently compute the Wasserstein distance. The MATS method \cite{rim2019manifold} transports reduced subspaces along characteristics and can not directly handle the case when characteristics intersect (e.g. problems with shocks). As demonstrated in our numerical examples, our proposed  method can
handle 2D problems and problems with shocks.
We also want to mention that \cite{rim2020depth} generalizes \cite{rim2019manifold} with deep learning techniques and could handle problems with shocks, but its connection with optimal transport  becomes less clear.}

The paper is organized as follows. In Section \ref{sec:background}, we briefly review the basic idea of  projection-based ROM methods and explain why designing an efficient ROM of transport problems could be challenging. Then, in Section \ref{sec:rom_algorithm}, we present the offline and online algorithms of the proposed learning-based projection method as well as the neural network architecture. In Section \ref{sec:numerical}, the performance of the proposed method is demonstrated through a series of numerical experiments. In Section \ref{sec:conclusions}, we draw our conclusions and discuss \mw{some} potential future research directions.

\section{Background\label{sec:background}}
We consider a parametric time-dependent transport problem:
\begin{align}
    \reviewone{\textrm{residual}(u;\mu) = 0,\;\mu\in X_\MP,}
    \label{eq:abstract_pde}
\end{align}
where $\mu$ is a scalar or vector-valued model parameter and \reviewone{$\textrm{residual}(\cdot;\mu)$}
%$r(\cdot;\mu)$ 
represents a time-dependent parametric hyperbolic partial differential  operator. The solution $u(x,t;\mu)$ depends on the location $x\in\Omega$, the time $t$ and the parameter $\mu$. To design a reduced order model,  we assume a full order model is available, which is 
a high-fidelity numerical scheme solving \eqref{eq:abstract_pde}  and can be written as 
\begin{align}\label{numerical_scheme}
    \reviewone{\textrm{residual}}_h(u_h;\mu)=0
\end{align} 
in its strong form.
We \minreview{view} the solution $u_h(x,t;\mu)$  as the ground truth, and denote the solution manifold as \pzcrev{$\{u_h(\cdot,t;\mu): 0\leq t\leq T_{\textrm{max}},\mu\in X_\MP\}$.} %\lli{/note: to revisit/}

A conventional projection-based ROM seeks a solution in the form of
\begin{align}
    u_r(x,t;\mu) = \sum_{i=1}^r \alpha_i(t,\mu)\varphi_i(x).
    \label{eq:linear_expansion_rom}
\end{align}
Here, $r$ is the order of the reduced order model, $\{\varphi_i(x)\}_{i=1}^r$ is a reduced basis, with the expansion coefficients
${\boldsymbol{\alpha}}(t,\mu)=(\alpha_1(t,\mu),\dots,\alpha_r(t,\mu))^T$.
The reduced basis depends on $x$, and is 
typically obtained through the offline training. \pzc{The coefficients are \mw{then} determined online through projection methods.} 
Such procedure provides efficient and accurate reduced order models for problems whose solution manifolds have \pzcreview{fast-decaying} Kolmogorov $n$-width \cite{melenk2000n, buffa2012priori}.

\begin{figure}[h!]
\centering
    \begin{tikzpicture}
         \node (img){\includegraphics[width=0.5\textwidth]{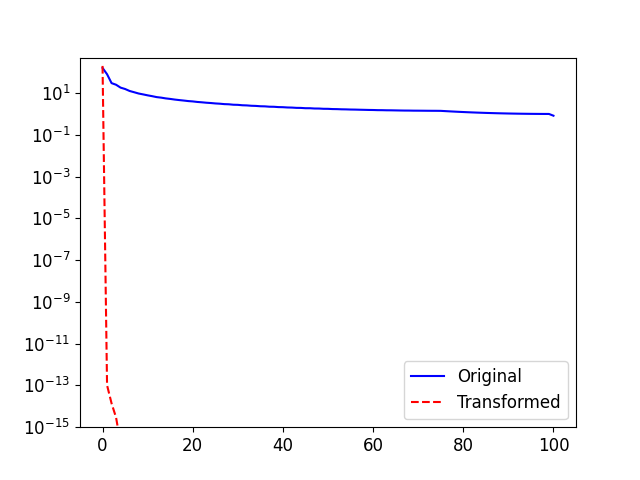}};
    \node[below =of img, node distance=0cm, yshift=1.5cm] {$n$};
  \node[left=of img, node distance=0cm, rotate=90, anchor=center,yshift=-1.2cm] {$\sigma_n$};
    \end{tikzpicture}
\caption{Singular values of the original and transformed snapshot matrices. \label{fig:svd_exact}}
\end{figure}

However, it is well-known that the Kolmogorov $n$-width 
for transport-dominated  problems
may decay slowly \reviewone{\cite{ohlberger2016reduced,benner2017model,greif2019decay,cohen2021optimal}}. On the other hand, the evidence in \cite{cagniart2017model} indicates that under some coordinate transformations, the transformed solution manifolds of transport-like problems can have much faster decay in their Kolmogorov $n$-width.  This can be illustrated by the following simple example. Consider the 1D advection equation
\begin{align}
    u_t+u_x = 0,\quad x\in[0,2],
  %  , \quad \fli{t\in(0,1]} 
\end{align}
with zero boundary conditions and the \review{following} initial condition
\begin{align}
u(x,0)=u_0(x)=\begin{cases}
           1, \quad 0.25\leq x \leq 0.5,\\
           0, \quad \text{otherwise.}
           \end{cases}
\end{align}
The exact solution  is  $u(x,t)=u_0(x-t)$ \fli{(over the time period $[0,1]$ of our consideration)}.
We introduce a uniform space-time mesh \mw{with $\Dx = \Dt=0.01$} and $x_i=i\Dx$,  $t_j=j\Dt$, and 
use the exact solution sampled from this mesh over $t\in[0,1]$ to define the solution manifold, whose  Kolmogorov n-width with respect to the \pzc{$l_2$} norm is measured as the $n$-th singular value of the snapshot matrix $S=(s_{ij})$, with  $s_{ij}=u(x_i, t_j)=u_0((i-j)\Dx)$. Meanwhile, we can follow \cite{cagniart2017model} and apply a $t$-dependent transformation $T(x;t)=x-t$, and obtain a transformed snapshot matrix $\tilde{S}=\{\tilde{s}_{ij}\}$ with $$\tilde{s}_{ij}=u(T^{-1}(x_i; t_j), t_j)=u(x_i+t_j, t_j)=u_0(x_i).$$ 
In Figure \ref{fig:svd_exact}, the singular values of the original and transformed snapshot matrices are presented. We observe that the singular values of the \pzc{original} snapshot matrix decay fairly slowly as expected, while those of the transformed snapshot matrix decay very fast.  As a result, standard ROMs such as POD will require a fairly large reduced space  for accurate approximation, yet with the  help of the \review{$t$-dependent} transformation $T(x;t)$, the dimension of the reduced space can be dramatically reduced.

The approximation power of the standard reduced order spaces (e.g. the POD space)  is limited \mw{as they are linearly expanded with basis functions that are static} with respect to $t$ and $\mu$, see \eqref{eq:linear_expansion_rom}. As shown above, one \mw{way} to conquer this limitation is to use a reduced order space with an adaptive time-and-parameter-dependent (i.e. ($t,\mu$)-dependent) basis $\{\phi_i(x;t,\mu)\}_{i=1}^r$, with the reduced order approximation given as:
\begin{align}
    u_r(x,t;\mu) = \sum_{i=1}^r \alpha_i(t,\mu)\phi_i(x;t,\mu).
    \label{eq:parametric_expansion_rom}
\end{align}
In the series of paper \cite{cagniart2017model,cagniart2019model,nonino2019overcoming}, the basis \mw{used takes the} form $\phi_i(x;t,\mu)=\psi_i(T_{t,\mu}(x))$, where $T_{t,\mu}$ is a suitable coordinate transformation defined problem by problem. Similar  ideas \mw{have} also been considered in \cite{rim2019manifold,nair2019transported,torlo2020model,taddei2020registration,mojgani2020physics,ferrero2021registration}.
\mw{In \cite{rim2019manifold},} $T_{t,\mu}$ is constructed \mw{with} the  displacement interpolation \cite{mccann1997convexity,villani2008optimal} \mw{while the} polynomial approximation \mw{is used} in \cite{nair2019transported}. \review{The Lagrangian and arbitrary Lagrangian Eulerian (ALE) registration method \cite{torlo2020model,taddei2020registration,ferrero2021registration} and the ALE autoencoder \cite{mojgani2020physics} {find}
$T_{t,\mu}$ by solving optimization problems.} 
 
Following \cite{cagniart2017model,cagniart2019model,nonino2019overcoming}, we propose a learning-based projection method. To avoid \mw{designing} suitable coordinate transformations $T_{t,\mu}$ \textit{a priori} for each problem, we  parameterize the basis functions $\{\phi_i(x;t,\mu)\}_{i=1}^r$ as 
a neural network $\NNphi(x,t,\mu;\theta^\phi)$ and \pzc{determine the trainable parameters of the neural network $\theta^\phi$} in the offline learning stage. In the online stage, we \mw{further utilize} \mw{a} projection method to predict \mw{solutions} either at a future time or corresponding to  a new value of the model parameter $\mu$. The ROM algorithms will be detailed in Section \ref{sec:rom_algorithm}.

\section{Reduced order model algorithms \label{sec:rom_algorithm}}
One natural extension to \cite{cagniart2017model,cagniart2019model,nonino2019overcoming} is \mw{utilizing} neural networks to learn \mw{a} coordinate transformation $T_{t,\mu}$. However, it is nontrivial to measure how a coordinate transformation reduces the Kolmogorov $n$-width of the transformed solution manifold for a general problem \pzcreview{with neural networks}, not to mention to systematically find a good one. 
Instead, we directly learn the $(t,\mu)$-dependent adaptive basis, \review{hence the adaptive reduced order space, using} neural networks, \review{ while the transformation $T_{t,\mu}$ is implicitly encoded.} 

We follow the standard offline-online decomposition of the projection-based ROM \pzc{method}.
 In Section \ref{sec:offline}, we present the offline stage algorithm which obtains the adaptive basis by training two neural networks \mw{with snapshot solutions}. In the online stage, with the learned $(t,\mu)$-dependent basis available, we apply a projection-based  method to construct a \mw{reduced representation of the solution}, and this will be \mw{detailed} in Section \ref{sec:online}. A discussion \mw{over} the architecture of the neural network is \mw{then} presented in Section \ref{sec:nn_architecture}.

%%%%%%%%%%%%%%%%%%%%%%%%%%%
\subsection{Offline stage\label{sec:offline}}

In the offline training stage, the snapshot solutions are high-fidelity numerical solutions from a full order model. For the time-dependent transport problems considered in this work, we assume the full order model is given by a one-step scheme 
 \begin{equation}
     \reviewone{\textrm{residual}_h(u_h^{n+1},u_h^n;\mu) = 0}
     \label{eq:residual}
 \end{equation}
 in its strong form, and it can be constructed through a finite difference/volume/element method, or a  discontinuous Galerkin method. Here $u_h^n=u_h(x,t_n;\mu)$ denotes
  the 
 numerical solution at time $t_n$. The snapshot solutions can also be obtained from measurements or practical experiments.
 
Suppose we have the snapshot solutions, $u_h(x_i,t_j;\mu_k)$, with
$$\{x_i\}_{i=0}^{N_x}\subset\Omega,\quad \{t_j\}_{j=1}^{N_t}\subset [0,T_{\textrm{train}}], \quad  \{\mu_k\}_{k=1}^{N_\mu}\subset \MP_{\textrm{train}}\subset X_\MP.$$
Note that the spatial degrees of freedom of $u_h$ used here are nodal values. Our proposed framework can be adapted to ROMs with other types of spatial degrees of freedom, such as moments. \lli{Though the  snapshot solutions are taken from $[0,T_{\textrm{train}}$], }
\pzcrev{
%Also,
the proposed method may be applied to predict solutions at a future time $t>T_\textrm{train}$.}

To obtain the $(t,\mu)$-dependent basis, we define the following \pzcreview{$r$-th order} neural network approximation:
\begin{align}
    u_{\NN}(x,t;\mu,\boldsymbol{\theta})
     &= \sum_{i=1}^{r}\alphaNN_i(t,\mu;\mw{\theta^{\alpha}})\ \phiNN_i(x,t,\mu;\mw{\theta^{\phi}})\notag\\
     &= \boldsymbol{\alphaNN}(t,\mu;\theta^\alpha)\cdot \boldsymbol{\phiNN}(x,t,\mu;\theta^\phi).
    \label{eq:NN_expansion}
\end{align}
Here, $\mw{r}\ll N_x$ is the order of the ROM. $\boldsymbol{\theta}=(\theta^\alpha,\theta^\phi)$ are trainable parameters \mw{of} the neural network. The $r$-dimensional vectors $\boldsymbol{\alphaNN}=(\alphaNN_1,\cdots,\alphaNN_r)^T$ and $\boldsymbol{\phiNN}=(\phiNN_1,\cdots,\phiNN_r)^T$ are the outputs of the neural networks $\NNalpha(t,\mu;\theta^\alpha)$ and $\NNphi(x,t,\mu;\theta^\phi)$, where $\NNalpha(t,\mu;\theta^\alpha)$ and $\NNphi(x,t,\mu;\theta^\phi)$ parameterize the expansion coefficients and the $(t,\mu)$-dependent basis, respectively.
The neural network is trained by minimizing the mean-square loss function:
\begin{align}\label{eq:loss_function}
\pzc{\boldsymbol{\theta^*} =(\theta^{\alpha,*},\theta^{\phi,*})= \arg\min_{\boldsymbol{\theta}}\frac{1}{(N_x+1)N_tN_\mu}\sum_{i=0}^{N_x}\sum_{j=1}^{N_t}\sum_{k=1}^{N_\mu}
\big(u_{\NN}(x_i,t_j;\mu_k,\boldsymbol{\theta})-u_h(x_i,t_j;\mu_k) \big)^2.}
\end{align}
The solution to the optimization problem in \eqref{eq:loss_function} is not unique in the sense that one can reorder and rescale basis functions. Keep in mind that what we truly want is not the basis itself but the associated adaptive reduced order space. {In practice, the ordering and the scaling of the learned basis will be influenced by the random initialization of the neural networks. Also, due to the non-convex nature of the underlying minimization problem, the optimizer may converge to a low dimensional manifold corresponding to a local minimum. But at least for numerical tests being considered, these issues do not seem to pose any problem in the present work, either to the training of the neural networks or to the overall performance of the proposed algorithm. Nevertheless, it will be left to our future exploration to learn more structured basis functions such as an orthonormal basis.}

Essentially, the neural network  $u_{\NN}(\cdot; \cdot, \mathbf{\theta^*})$ provides a reduced order model to approximate the true solution\mw{.} \mw{However,} due to \mw{the potential risk of over-fitting,}
\pzc{the purely learning-based reduced approximation $u_{\NN}$} is not always accurate and robust \mw{in} generalization.  We propose to only pass the learned adaptive basis $\{\phiNN_i\}_{i=1}^{r}$  onward to the online stage.

It is worth noticing that $\NNphi$ takes $x$ as one of its inputs in a mesh free manner, and therefore it is %\mw{thus able to be optimized}
simple and flexible to work
with \pzc{solution snapshots obtained} on different meshes tailored for different values of the model parameter, or even solution snapshots based on experimental data. 
However, a conventional projection-based method may need extra steps such as interpolations \cite{grassle2018pod} or full order offline re-computations \cite{yano2019discontinuous,yano2020goal}  in order to work with snapshots subject to different \pzc{meshes}. 

%\reviewone{
\begin{rem}
For now, the reduced order $r$ is pre-selected. One way to adaptively determining $r$ is to apply a block-by-block training strategy similar to \cite{dahmen2021nonlinear}. More specifically, start
with a pre-selected small $r_0$ and train the neural network $u_{\NN_0}=\sum_{i=1}^{r_0}\alphaNN_i(t,\mu;\theta^\alpha)\phiNN_i(x,t,\mu;\theta^\phi)$. If the prediction error of a validation set is too large with the pre-selected $r_0$, one can train a second neural network $u_{\NN_1}=\sum_{i=r_0+1}^{r_0+\Delta r}\alphaNN_i(t,\mu;\theta^\alpha)\phiNN_i(x,t,\mu;\theta^\phi)$ by minimizing the mismatch between $u_{\NN_0}+u_{\NN_1}$ and the full order solution $u_h$. Repeat this procedure until the prediction error of the validation set is small enough or reach the maximum number of basis functions allowed.
\end{rem}
%}

%\reviewmulti{
\begin{rem}
In practice, it is important to have an effective sampling strategy when  the dimension of the problem or the parameter $\mu$ is high. In order to reduce the computational cost of sampling, one may apply an adaptive sampling strategy based on the greedy algorithm \cite{maday2013locally,hesthaven2016certified} or the ``data exploration" in \cite{han2019uniformly}. Effective sampling strategy is not the focus of this paper, and for now, we do not apply these strategies and leave them for future investigations. \lli{Similarly,}
\pzcrev{the terminal training time $T_\textrm{train}$ can also be adaptively determined, e.g. based on the validation error or an error estimator.}
\end{rem}
%}
%%%%%%%%%%%%%%%%%%%%%%%%%%%
\subsection{Online stage\label{sec:online}}
In the online stage, we take the learned adaptive basis $\{\phiNN_i\}_{i=1}^{r}$ 
and seek a reduced-order numerical solution in their span
following a projection approach. For simplicity, we denote $\phi_i^{n+1}(x; \mu)=\phi_i^{\NN}(x, t_{n+1}, \mu; \theta^{\phi,*})$, $i=1, \cdots, r$. The reduced order space at $t_{n+1}$ associated with $\mu$ is then given as  
$$V^{n+1}_r(\mu) :=  \text{span}\{\phi_i^{n+1}(\cdot; \mu) \}_{i=1}^{r}.$$
We also introduce $\boldsymbol{\Phi}_\mu^{n+1}\in \mathbb{R}^{(N_x+1)\times r}$, with its $(i,j)$-th entry being
\begin{equation} 
\boldsymbol{\Phi}_\mu^{n+1}(i,j)=\phi_j^{n+1}(x_i;  \mu), \;\; i=0,\cdots, N_x,  \;\; j=1, \cdots,r.
\label{eq:ROM:mat}
\end{equation}

Suppose the full order model \eqref{eq:residual} has the following matrix-vector form:
\begin{equation}
\boldsymbol{F}_\mu({\bf u}_\mu^{n+1}) = {\bf{b}}_\mu^{n},
\label{eq:res:v}
\end{equation}
where ${\bf{u}}_\mu^{n+1}\in\mathbb{R}^{N_x+1}$ is the vector of nodal values (i.e.  the spatial degrees of freedom) of the  numerical solution $u_h(x,t_{n+1};\mu)$ at $\{x_j\}_{j=0}^{N_x}$ and time $t_{n+1}$. ${\bf{b}}_\mu^n\in \mathbb{R}^{N_x+1}$ is a known vector determined by ${\bf{u}}_\mu^n$, and the operator $\boldsymbol{F}_\mu: \mathbb{R}^{N_x+1}\mapsto \mathbb{R}^{N_x+1}$ is derived from the full order model   \eqref{eq:residual}. 
We look for the reduced order approximation of the solution $u_{r,h}^{n+1}\in V^{n+1}_r(\mu)$, whose nodal values at $\{x_j\}_{j=0}^{N_x}$ (i.e. spatial degrees of freedom) are $\boldsymbol{\Phi}_\mu^{n+1}\boldsymbol{\alpha}^{n+1}$, 
with some $\boldsymbol{\alpha}^{n+1}=(\alpha_1^{n+1},\cdots,\alpha_r^{n+1})^T\in\mathbb{R}^r$.
 Following a projection method, we require 
\begin{align}
    (\boldsymbol{\Phi}^{n+1}_{\mu})^T\boldsymbol{F}_\mu(\boldsymbol{\Phi}^{n+1}_{\mu}\boldsymbol{\alpha}^{n+1}) = (\boldsymbol{\Phi}^{n+1}_{\mu})^T{\bf{b}}^n_\mu.
\end{align}
Once ${\boldsymbol{\alpha}}^{n+1}$ is calculated, we reach a reduced-order approximation of the solution at $t_{n+1}$: $u_{r,h}^{n+1}=\sum_{i=1}^r\alpha_i^{n+1}\phi_i^{n+1}$. The online algorithm is summarized in Algorithm \ref{alg:onlinee}.

\begin{algorithm}[!h]
\caption{Online projection for one time step update \label{alg:onlinee} }
\begin{algorithmic}[1]
\STATE{{\bf{Input:}} given the parameter $\mu$; assume the reduced order  numerical solution $u_{r,h}^{n}$ at time $t_n$ is known, and it defines  ${\bf{b}}_\mu^n$ in  \eqref{eq:res:v} based on the full order model \eqref{eq:residual}.
}
\STATE{{\bf{Evaluate the $(t,\mu)$-dependent basis with the neural network $\NNphi$ trained offline:}}
\begin{equation}
\mw{\phi_i^{n+1}}(x;\mu) := \phiNN_i(x,t_{n+1},\mu;\theta^{\phi,*}), \quad i=1,\cdots,r.
\label{eq:online_basisGen}
\end{equation}
And form the matrix  $\boldsymbol{\Phi}_\mu^{n+1}\in \mathbb{R}^{{(N_x+1)}\times r}$ as in  \eqref{eq:ROM:mat}.
} 
\STATE{{\bf{Obtain $u_{r,h}^{n+1}$ through the projection:}} compute
$\boldsymbol{\alpha}^{n+1}$ by}
\begin{align}
    (\boldsymbol{\Phi}^{n+1}_{\mu})^T\boldsymbol{F}_\mu(\boldsymbol{\Phi}^{n+1}_{\mu}\boldsymbol{\alpha}^{n+1}) = (\boldsymbol{\Phi}^{n+1}_{\mu})^T{\bf{b}}_\mu^n.        \label{eq:online_projection}
\end{align}
\State{{\bf{Output}}: the reduced order numerical solution at $t_{n+1}$ for the parameter $\mu$ is  $$u_{r,h}^{n+1}=\sum_{i=1}^r\alpha_i^{n+1}\phi_i^{n+1}.$$} 
\end{algorithmic}
\end{algorithm}

Compared with \mw{a} standard projection-based method such as the POD method, in the online stage, we have \mw{the} extra step to \mw{evaluate} the $(t,\mu)$-dependent basis \mw{with} a trained neural network with known parameters \mw{$\theta^{\phi,*}$} at  each time step \mw{$t_{n+1}$}. In our numerical tests, we find that, even with 
\mw{this} additional cost, our method still leads to computational savings compared with a full order model \review{when the underlying mesh is well refined}.

As mentioned before, the neural network \mw{$\NNphi$} takes $x$ as one of its inputs in a mesh free manner. \mw{Although it is} not \review{deeply} investigated in this work, the proposed method has the capability to be integrated with $h$-adaptive numerical schemes for time dependent problems in the online stage. \reviewtwo{In Section \ref{sec:moving_mesh}, the proposed method is combined with a moving mesh method to show its capability.}

  We also want to compare our proposed method with the learning-based projection methods in \cite{lee2019deep,lee2020model,kim2020fast}. 
In the offline stage, taking the full order solution as the input, these methods find %s 
 a low-dimensional nonlinear manifold \mw{to represent the full order solution manifold} by training an auto-encoder. 
In the online stage, these methods seek a reduced order approximation in the form
\begin{align}
    u\approx u_{\textrm{ref}}+g(\hat{u}_r).
\end{align}
Here, $u_{\textrm{ref}}$ is a \mw{predetermined} reference frame, $g$ is the decoder part of the neural-network trained offline, which is a mapping from a low-dimensional latent space to the high-dimensional solution manifold, and $\mw{\hat{u}_r}$ can be seen as the coordinates in the  low-dimensional latent space. In the \mw{online stage}, these methods solve %s 
a minimal residual problem
\begin{align}
   \hat{u}_r^{n+1} = \arg\min ||\reviewone{\textrm{residual}}_h(g(\hat{u}_r^{n+1}),u_h^n;\mu) ||_2.
\end{align}
To obtain $\hat{u}^{n+1}_r$, \cite{lee2019deep,lee2020model} need to take the derivatives or even the second order derivatives of the neural-network $g$ online\mw{, which can be computationally expensive}.
To improve the online efficiency, \cite{kim2020fast} \mw{uses} a shallow architecture for the decoder neural network $g$ and \mw{applies the} hyper-reduction techniques. \mw{Unlike these methods}, 
%as can be seen in \eqref{eq:online_projection},
we only need to compute $\boldsymbol{\alpha}^{n+1}$ as in \eqref{eq:online_projection} and
there is no need to take derivatives of neural networks at all \mw{in our method}. 
Furthermore, since the input of our method is not full/reduced-order solutions as of the auto-encoders,  
our method is relatively more flexible to work with  different meshes in the online and offline stages. \reviewtwo{We also want to point out that \cite{fresca2021comprehensive} proposes an alternative data-driven online phase for this auto-encoder and deep learning based framework. The data-driven online phase of \cite{fresca2021comprehensive} is more efficient than the projection-based online phase in \cite{lee2019deep,lee2020model,kim2020fast} and our method, while the projection-based online algorithm may lead to better accuracy.}

% \li{/note: shall we delete this, as this entire paragraph is also in introduction?/} The main difference between our method and the optimal transport based method \cite{rim2019manifold}, optimization based registration/autoencoder \cite{torlo2020model,taddei2020registration,mojgani2020physics,ferrero2021registration}  and the transported snapshot method \cite{nair2019transported} is as follows. We find the adaptive reduced order space without explicitly finding the transformation, while they explicitly find the transformation. Our method can be seen as an alternative to these methods. 

\begin{rem}
In our method, 
we require the projection of the residual of the full order model \eqref{eq:residual} onto the local reduced order space to vanish.
Alternatively, in our method, one can also use a
minimal residual formulation  similar to \cite{lee2020model,kim2020fast,cagniart2017model}:
\begin{align}
    \boldsymbol{\alpha}^{n+1} = \arg\min_{ \boldsymbol{\alpha\in \mathbb{R}^r}} ||\reviewone{\textrm{residual}}_h(\sum_{i=1}^r\alpha_i\phi_i^{n+1},u_h^n;\mu) ||_2.
\end{align}
In the matrix-vector formulation, this alternative online step becomes 
\begin{align}
    (\boldsymbol{\Phi}^{n+1}_{\mu})^T(\boldsymbol{J}_\mu^{n+1})^T\boldsymbol{F}_\mu(\boldsymbol{\Phi}^{n+1}_{\mu}\boldsymbol{\alpha}^{n+1}) = (\boldsymbol{\Phi}^{n+1}_{\mu})^T\reviewtwo{(\boldsymbol{J}_\mu^{n+1})^T}{\bf{b}}^n_\mu.
\end{align}
Here, $\boldsymbol{J}_\mu^{n+1}$ is the Jacobian matrix of $\boldsymbol{F}_\mu$ with respect to $\boldsymbol{u}_{\mu}^{n+1}$, and this can be seen as a weighted projection.  
\end{rem}

\begin{rem}
\reviewtwo{To work with the operator $(\boldsymbol{\Phi}^{n+1}_{\mu})^T\boldsymbol{F}_\mu(\boldsymbol{\Phi}^{n+1}_{\mu})$  in \eqref{eq:online_projection}, an $O(N_x)$ computational cost will be needed per time step. Unlike ROM methods (e.g. POD)  with static bases, the reduced basis in our method is  updated dynamically and can not be precomputed. As a result, the proposed method will not save computational time without hyper-reduction when explicit time steppers are applied. Our method  is indeed more suitable when implicit time integrators are preferred, for example, in the cases of very non-uniform spatial meshes or for stiff problems.}

\pzcrev{Hyper-reduction strategies have been considered for the ROMs of transport dominant problems. 
\lli{To predict solutions for unseen values of physical parameter $\mu$ (yet over $[0,  T_\textrm{train}]$), the autoencoder neural network based projection methods
\cite{kim2020fast} and \cite{romor2022non} apply  the DEIM-SNS/GNAT-SNS  \cite{choi2020sns} and the reduced over-collocation method  \cite{chen2021eim} for the hyper-reduction, respectively.}
For steady state problems, the $l_0$ minimization based strategy in \cite{ferrero2021registration} has been designed. When predicting the solutions for unseen \lli{values of} %physical 
%parameters
$\mu$ at time \lli{$t\in [0,  T_\textrm{train}]$,}
%$0\leq t\leq T_\textrm{train}$,
a hyper-reduction strategy similar to the aforementioned methods may be combined with our method. However, to accurately predict solutions at future time $t>T_\textrm{train}$, an adaptive hyper-reduction strategy may be necessary, and the design of such adaptive strategy could be challenging.} 

\reviewtwo{ Alternative to the hyper-reduction techniques, one can also enhance the computational efficiency of explicit time steppers by freezing the basis over a few time steps as in  \cite{peherstorfer2020model}.  This will lead to an $O(N_x)$ computational cost for one time step followed by an $O(r)$ cost for a few time steps. The key is to design an adaptive strategy determining when to update the basis for a good balance of efficiency and accuracy.}

\end{rem}

\subsection{Neural network architecture \label{sec:nn_architecture}}

As indicated by \eqref{eq:NN_expansion}, we  parameterize the reduced-order basis and the corresponding coefficient as independent neural networks, namely $\NNalpha(\mu,t;\mw{\theta^{\alpha}})$ and $\NNphi(x,t,\mu;\mw{\theta^{\phi}})$, respectively.
More specifically, we require each neural network to output a vector of dimension $r$, where each component stands for one dimension in the reduced latent space. \pzc{The inner product between the output vectors of the two neural networks} will provide a reduced order approximation for the solution. 

There is much flexibility in choosing  the neural network architecture. In this paper, we take each neural network to be simply a feed-forward fully-connected neural network. A standard $n$-layer architecture of such neural network has a repeated compositional structure, where each layer is composed of a linear function and a nonlinear activation functions $\sigma$:
$$ NN(\mathbf{x};\boldsymbol{\theta}) =W_n\sigma\left(\cdots\sigma \left( W_2\sigma \left(W_1\mathbf{x}+b_1\right)+b_2 \right)\cdots\right)+b_n.$$
Here, $\boldsymbol{\theta}$ stands for the collection of all trainable parameters, that is, $\boldsymbol{\theta} =\{W_1,b_1, \cdots, W_n,b_n\}$.
Other architectures can also be considered within our proposed reduced-order \pzc{model} framework. An interesting choice will be the partition of unity (POU) neural network \cite{lee2021partition}, which is related to $h,p$-adaptive methods and can \review{potentially} be beneficial for the construction of ROMs for transport problems.  \mw{Another} architecture that can be potentially used is the DeepONet \cite{lu2019deeponet}. The DeepONet can take  $\mu$ and $x,t$ in different ways, and this is natural in parameterizing a parameter-dependent function.

\begin{figure}[h!]
    \centering
    \includegraphics[width = 0.6\textwidth]{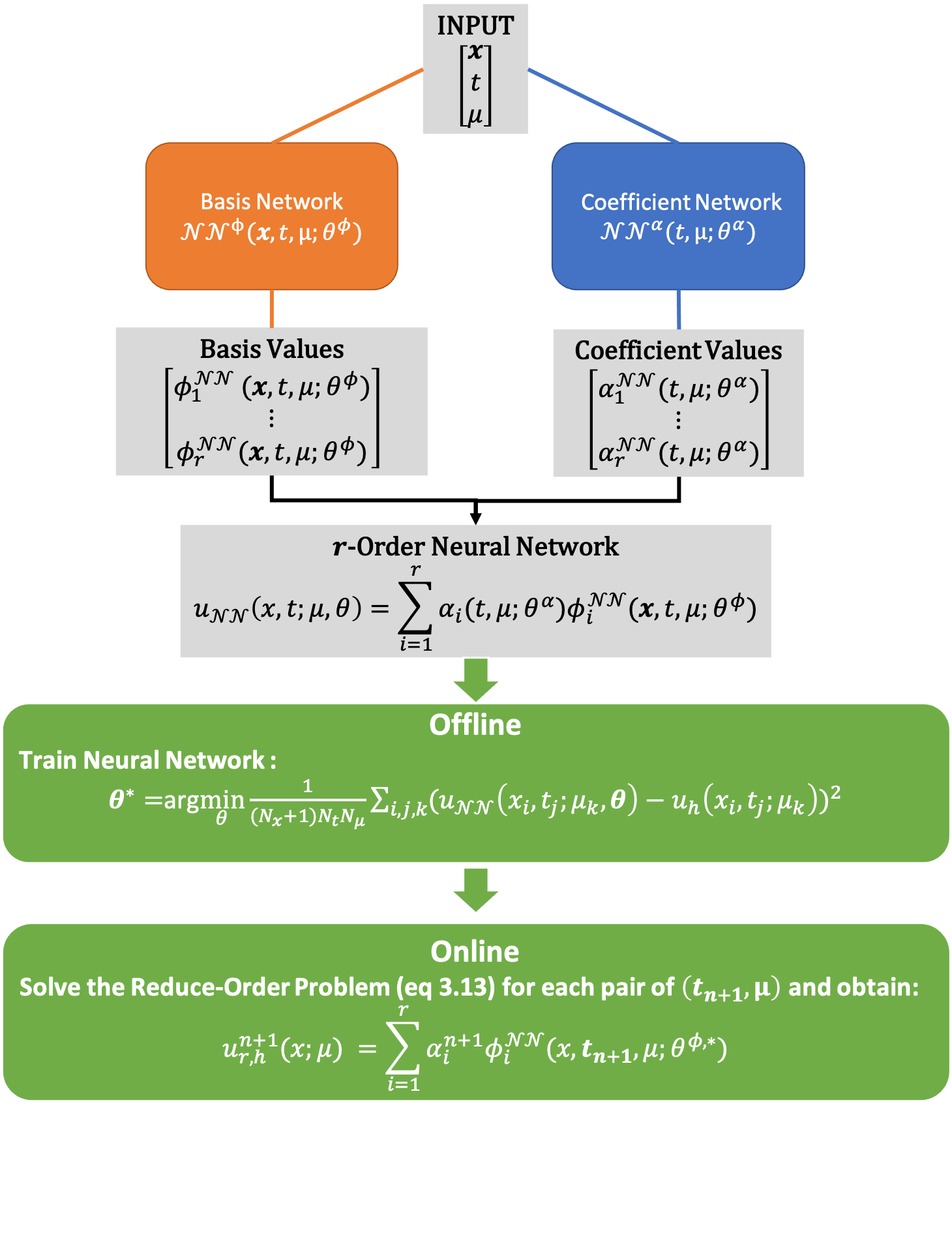}
    \caption{The flow chart of the proposed algorithm}
    \label{fig:nn_architecture}
\end{figure}

The entire algorithm can be summarized as the flow chart in Figure \ref{fig:nn_architecture}. 

\section{Numerical experiments\label{sec:numerical}}

In this section,
we demonstrate the performance of the proposed learning-based projection (LP) method by applying it to
\review{linear and nonlinear 
%scalar or system of 1D \pzcreview{or 2D} 
hyperbolic conservation laws in 1D and 2D.} Particularly, the \review{1D and 2D} linear advection equation, the \review{1D and 2D} Burgers equation, and the \review{1D} compressible Euler system of gas dynamics are tested. 
 We also compare our method with the standard POD method and  the prediction directly provided by the neural network $u_{\NN}$ in \eqref{eq:NN_expansion}.
Throughout this section, the reduced order solution obtained with the learning-based projection method  is referred to as ``LP"  and the one  directly predicted by the neural network $u_{\NN}$
is referred to as ``learning".

In our numerical experiments, we adopt \pzc{the feed forward \mw{neural network} architecture} for both $\NNphi$ and $\NNalpha$.  
\mw{T}he activation function is taken as $\tanh$, \mw{e}xcept for the output layers \mw{which are defined as} 
\begin{align}
    {\bf{y}}_{\textrm{output}} = \textrm{softmax}(W{\mathbf{x}})+b,
\end{align}
\mw{for both  the basis and coefficients.}
\mw{Here,} $W$ is the weight and $b$ is the bias and $\mathbf{x}$ stands for the output of the previous layer.
Among different activation functions of the output layer we have numerically experimented with, $\textrm{softmax}$ outperforms other choices.
\mw{As} pointed out in \cite{lee2021partition},  the $\textrm{softmax}$ activation function introduces $h$-adaptivity to some extent, which may have made the neural network more suitable for the transport problems with local structures. \reviewthree{As presented in Section \ref{sec:2d_advection}, the proposed method is not sensitive to the choice of the width and the depth of the hidden layers of the neural network.} \review{Unless specified}, the neural network $\NNphi$ has $4$ hidden layers \mw{of} $25$  neurons \review{per layer} and an output layer with $r$ neurons. The neural network $\NNalpha$ has $3$ hidden layers \mw{of} $25$ neurons \review{per layer} and an output layer with $r$ neurons. 

Our code is implemented under the framework of Tensorflow 2.0 \cite{abadi2016tensorflow}.  \reviewtwo{When an implicit time integrator such as backward Euler or Crank-Nicolson method is applied, {the associated linear or linearized  systems are solved} by the GMRES method. If the equation is nonlinear, the Newton's method will be applied first as an outer loop iterative  solver.}

\subsection{Preliminaries 
\label{sec:numerical_full}}

\review{To set the stage, we start with a \mwmw{1D} scalar} 
hyperbolic conservation law 
\begin{align}
    u_t + \left( f\left(u\right) \right)_x = 0, \;x\in [x_L,x_R]
\end{align}
with an initial condition and some suitable boundary conditions.

Let $x_L = x_0< x_1 \dots < x_{N_x}=x_R$ be a uniform mesh in space, with $x_i = x_L+i\Dx$ and $\Dx=\frac{x_R-x_L}{N_x}$, and let $\Dt$ be the time step size. One of the full order models used   in our experiments is the following first-order  method, that involves the   backward Euler method in time and  a first-order conservative finite difference  method in space:  
\begin{equation}
\frac{u_j^{n+1}-u^n_j}{\Dt} +\frac{f^{n+1}_{j+\half}-f^{n+1}_{j-\half}}{\Dx} =0.
\label{eq:rev:med1}
\end{equation}
Here, $u_j^n\approx u(x_j,t_n)$, $f^n_{j+\frac{1}{2}}$ is a numerical flux at $x_{j+\frac{1}{2}}=\frac{x_j+x_{j+1}}{2}$, which will be specified \review{later.}
{Boundary conditions are imposed through numerical fluxes.} \review{The method in \eqref{eq:rev:med1} can be easily extended to non-uniform meshes, and it can also be  interpreted as a finite volume method, with $u_j^n$ as the element average.}

Throughout this section, the $L_2$ error and the relative $L_2$ error for \mw{the solution associated with} a test parameter $\mu$ at $t_n$ are defined as
\begin{subequations}
\begin{align}
    &\mathcal{E}^{n}_{L_2}(\mu) = \sqrt{\sum_{j=0}^{N_x} (u_{\textrm{reduced},j}^n-u_{\textrm{full},j}^n)^2\Dx },\\
    &\mathcal{E}^{n}_{\textrm{relative}}(\mu) = \frac{\mathcal{E}^{n}_{L_2}(\mu)}{\sqrt{\sum_{j=0}^{N_x} (u_{\textrm{full},j}^n)^2\Dx }},
\end{align}
\end{subequations}
respectively. The average relative $L_2$ error for a test set $\{\mu_i\}_{i=1}^{N_\mu}$ from $t_{n_0}$ to $t_{n_m}$ is defined as
\begin{align}
    \mathcal{E}_{\textrm{average}} = \frac{1}{N_\mu}\sum_{i=1}^{N_\mu}\max_{n_0\leq n\leq n_m}\left\{\mathcal{E}^{n}_{\textrm{relative}}(\mw{\mu_i})\right\}.
\end{align}
\review{The 2D analog of these quantities can be defined similarly.}
%%%%%%%%%%%%%%%%%%%%%%%%%%%%%%%%%%%%%%%%%%%%%%%%%%%%%%%%%%%%%%%%%%%%%%%%%%%%%%%%%
% Basis shifting
%%%%%%%%%%%%%%%%%%%%%%%%%%%%%%%%%%%%%%%%%%%%%%%%%%%%%%%%%%%%%%%%%%%%%%%%%%%%%%%%%
\subsection{\review{1D linear advection equation: learned basis functions}}
\label{sec:linear_advection}
We first consider the \review{1D} linear advection equation of the constant velocity 1: 
\begin{align}\label{eq:advection}
    u_t+u_x = 0,
\end{align}
with the initial condition
\begin{align}
    u(x,0)=\begin{cases}
           2, \quad 0.5\leq x \leq 1.5,\\
           0, \quad \text{otherwise.}
           \end{cases}
\end{align}
The computational domain is $[0,5]$, and it is partitioned \mw{with} a uniform mesh with $\Dx=0.01$. The time step size is taken as \review{$\Dt=\Dx$}. \review{We use the exact solution from $t=[0,0.75]$ to train the neural network with {the reduced order} $r=3$. The purpose of this example is to test whether the proposed method is able to capture the underlying {$t$-dependent transformation}.} 

\begin{figure}[h!]
\centering
\begin{subfigure}{.45\textwidth}
  \centering
  \includegraphics[width=\textwidth]{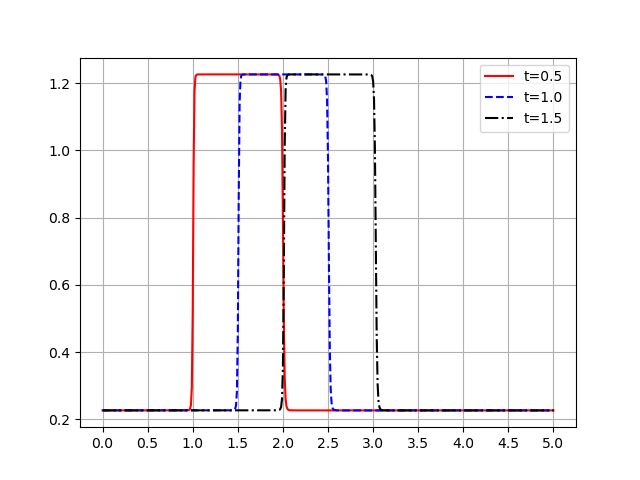}
  \caption{Basis function 1}
 \end{subfigure}
 \begin{subfigure}{.45\textwidth}
  \centering
  \includegraphics[width=\textwidth]{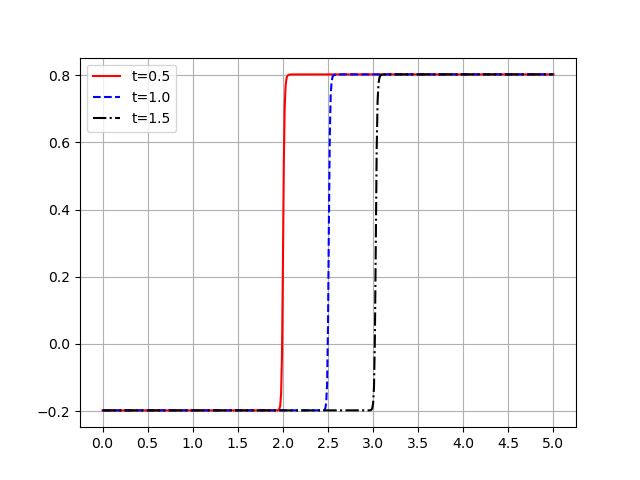}
  \caption{Basis function 2}
 \end{subfigure}
 \begin{subfigure}{.45\textwidth}
  \centering
  \includegraphics[width=\textwidth]{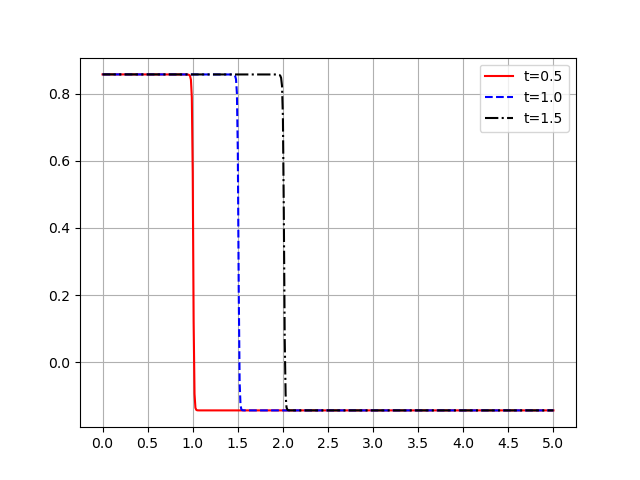}
  \caption{Basis function  3}
 \end{subfigure}
\caption{The learned basis functions \mw{$\phiNN_i$ for $i =1,2,3$} at different times for the linear advection equation %\eqref{eq:advection}
in Section \ref{sec:linear_advection}.
\label{fig:linear_advection_basis}} 
\end{figure}

\begin{figure}[h!]
\centering
  \includegraphics[width=0.5\textwidth]{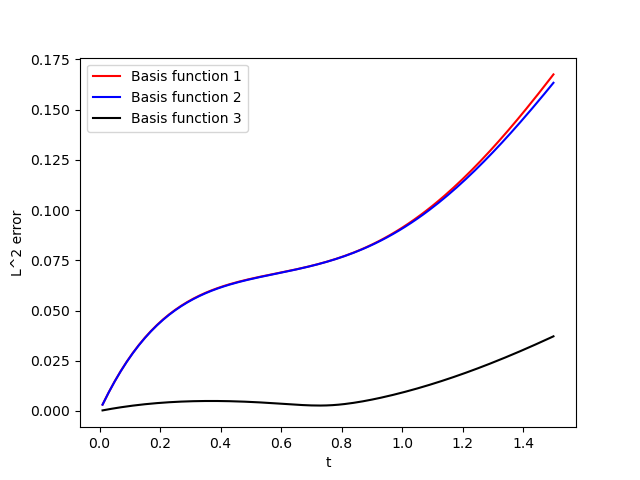}
\caption{\reviewone{The difference between the learned basis functions and the shifted initial basis functions for the linear advection equation in Section \ref{sec:linear_advection}.}
\label{fig:linear_advection_basis_error}}
\end{figure}

In Figure \ref{fig:linear_advection_basis}, we present  the learned basis functions \review{at three different times} obtained offline. We observe that, as time evolves, all  three basis functions propagate from the left to the right with \review{a wave speed approximately $1$}. This qualitatively verifies that the proposed LP method is able to capture the underlying physics. \reviewone{In Figure \ref{fig:linear_advection_basis_error}, we further qualitatively measure the difference between the learned basis functions and shifted initial basis functions with wave speed $1$. As time evolves, the difference will be larger.}

%%%%%%%%%%%%%%%%%%%%%%%%%%%%%%
%%%%%%%%%%%%%%%%%%%%%%%%%%%%%%
%%%%%%%%%%%%%%%%%%%%%%%%%%%%%%

%%%%%%%%%%%%%%%%%%%%%%%%%%%%%%%%%%%%%%%%%%%%%%%%%%%%%%%%%%%%%%%%%%%%%%%%%%%%%%%%%
% Moving mesh
%%%%%%%%%%%%%%%%%%%%%%%%%%%%%%%%%%%%%%%%%%%%%%%%%%%%%%%%%%%%%%%%%%%%%%%%%%%%%%%%%
\subsection{\review{1D linear advection equation with a moving mesh method}\label{sec:moving_mesh}}
\reviewtwo{The aim of this example is to show the flexibility and the effectiveness of our method to couple with a  full order moving mesh method. }

\reviewtwo{We consider the same linear advection equation as in Section \ref{sec:linear_advection}, but with a slightly different initial condition
\begin{align}
    u(x,0)=\begin{cases}
           2, \quad 0.6<x<1.4,\\
           0, \quad \text{otherwise.}
           \end{cases}
\end{align}
We use a non-uniform moving mesh: $0=x_0(t)<x_1(t)<\dots<x_{N_x}(t)=5.0$, with {$N_x=1000$.} Define $h_{\textrm{coarse}}=\frac{1}{125}$ and $h_{\textrm{refine}}=\frac{1}{500}.$ The grid point at time $t$, {namely, $x_i(t)$ ($1\leq i\leq N_x-1$),} satisfies
\begin{align}
    x_i(t) =\begin{cases}
            x_{i-1}(t)+h_{\textrm{refine}},\quad \text{if\;} 0.5+t\leq x_{i-1}(t)\leq 1.5+t,\\
            x_{i-1}(t)+h_{\textrm{coarse}},\quad \text{otherwise}.
    \end{cases} 
\end{align}
This mesh moves with the underlying wave speed and is always more refined near the discontinuity of the solution.
With this moving mesh, the first order upwind finite volume method with {
$f^n_{j+\frac{1}{2}}=u_j^n$ is applied in space, 
together with the backward Euler method in time} and the time step size $\Dt=h_{\textrm{coarse}}$. Before marching to the next time step, an extra key step of 
%the
{a
moving} mesh method is to interpolate the solution from the mesh at $t^{n-1}$ to the mesh at $t^n$.
{In our simulation,} we follow the Step 2 of the moving mesh algorithm \cite{tang2003adaptive} to do 
this interpolation. When our LP method is applied, we use the same interpolation subroutine to interpolate predicted solution from the old mesh to the new mesh.}
\reviewtwo{Extra cost reduction can be done in this step, and this is not explored here.} 

\reviewtwo{The training set is the snapshot solutions for $t\in[0,0.75]$ \pzcnew{($95$ time snapshots)}, and we will predict the solution from $t=0.75$ to $t=1.5$. With the moving mesh, each solution snapshot {corresponds to} a different mesh. As a result, interpolation of basis functions from different meshes \cite{grassle2018pod} or extra refinement and {re-computation} %s 
with a more refined {common} mesh \cite{yano2019discontinuous,yano2020goal} may be needed for {the} standard POD or RB method. Our neural networks take $x$, $t$ as inputs in a mesh-free manner, { and hence do not need } additional treatments  
{to work with} 
basis functions associated with different meshes. Here, we only present results of the pure learning based method and the LP method.}

\begin{figure}[h!]
\centering
\begin{subfigure}{.45\textwidth}
  \centering
  \includegraphics[width=\textwidth]{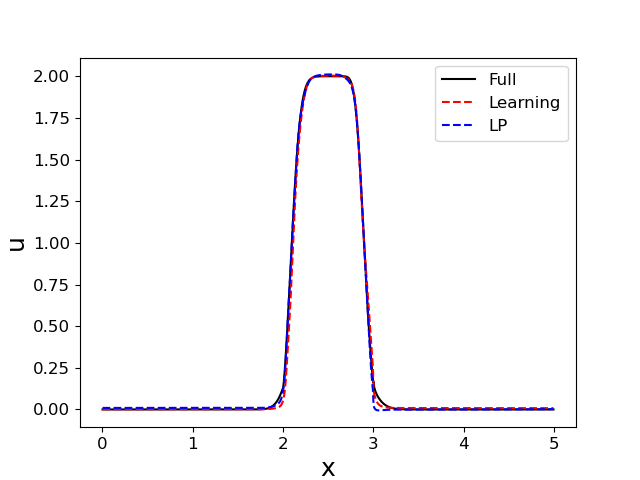}
  \caption{\pzcreview{Predicted solutions with 15 reduced basis functions} }
\end{subfigure}
\begin{subfigure}{.45\textwidth}
  \centering
  \includegraphics[width=\textwidth]{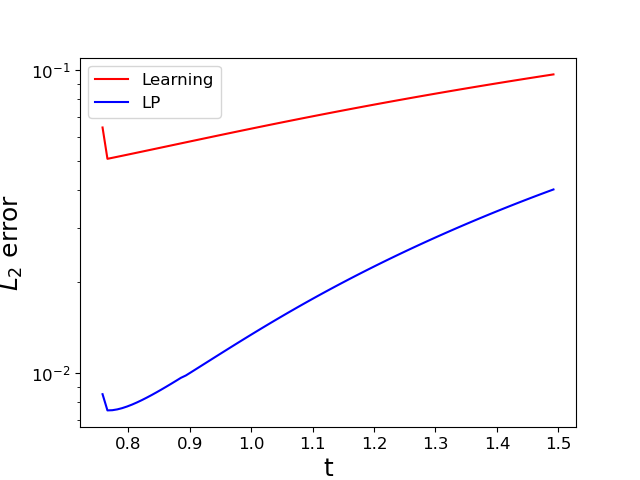}
    \caption{\pzcreview{Error history with 15 reduced basis functions}}
\end{subfigure}
\caption{\pzcreview{Predicted solutions at $t=1.5$ (left) and the error history (right) with 15 reduced basis functions for the moving mesh test in Section \ref{sec:moving_mesh}.}\label{fig:moving_mesh_sol}}
\end{figure}

\begin{figure}[h!]
\centering
\begin{subfigure}[t]{.45\textwidth}
  \centering
  \includegraphics[width=1.13\textwidth,]{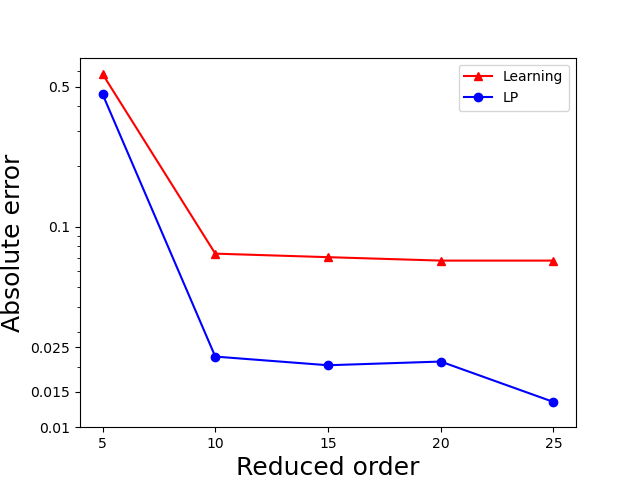}
  \caption{\pzcreview{Average absolute error
  %.
  }\label{fig:moving_mesh_error} }
\end{subfigure}
\hspace{0.01cm}
\begin{subfigure}[t]{.45\textwidth}
  \centering
  \includegraphics[width=\textwidth]{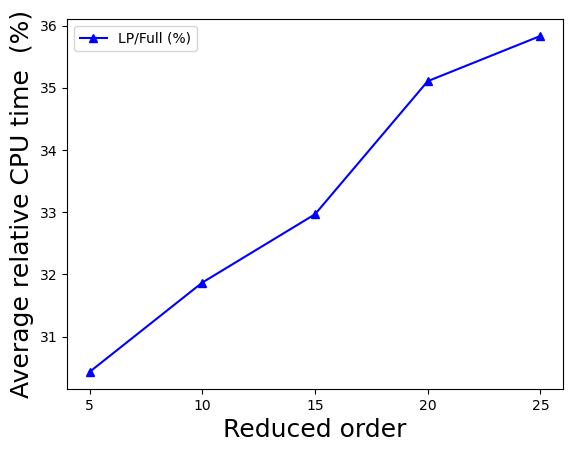}
    \caption{\reviewone{Average relative CPU time to invert the linear system for one step  time marching
    %.
    }\label{fig:moving_mesh_time}}
\end{subfigure}
\begin{subfigure}[t]{.45\textwidth}
  \centering
  \includegraphics[width=1.13\textwidth,]{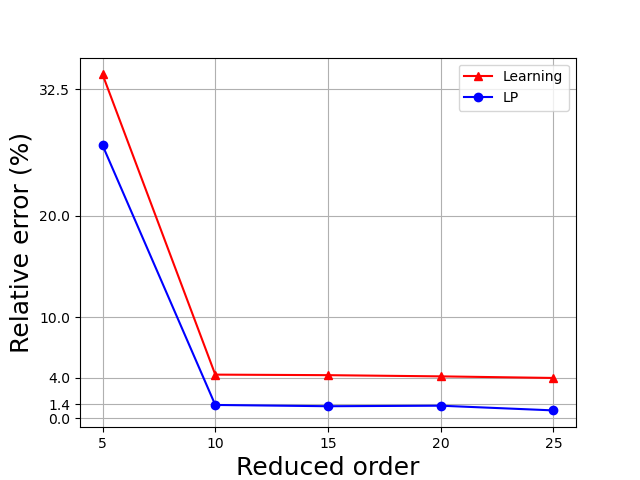}
  \caption{\pzcrev{Average relative error
  %.
  }\label{fig:moving_mesh_rel_error} }
\end{subfigure}
\hspace{0.01cm}
\begin{subfigure}[t]{.45\textwidth}
  \centering
  \includegraphics[width=1.13\textwidth,]{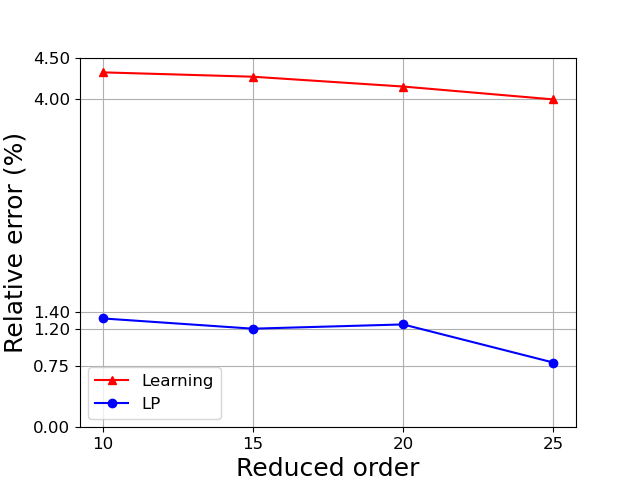}
  \caption{\pzcrev{Average relative error (zoomed-in)
  }\label{fig:moving_mesh_rel_error_zoomed_in} }
\end{subfigure}
\caption{\review{Average absolute/\pzcrev{relative} error  and average relative CPU time to invert the linear system for one step time marching, versus the reduced order $r$  \review{for the moving mesh test in Section \ref{sec:moving_mesh}.}}}
\end{figure}

\reviewtwo{In Figure \ref{fig:moving_mesh_sol}, we present the predicted solution at $t=1.5$ and the error history with {$r=15$. }  Both the learning based and the LP methods match
%matches
the full order solution well, and the error of the LP method is {overall} smaller. \lli{From Figure \ref{fig:moving_mesh_error} (resp. Figures \ref{fig:moving_mesh_rel_error}- 
\ref{fig:moving_mesh_rel_error_zoomed_in}) with $r=5,10,15,20,25$,
 we further observe that the LP method leads to smaller average absolute (resp. relative) errors} than the learning based method.}

\reviewone{Finally we want to compare the computational efficiency between the LP method and the full order method. As it is mentioned, both methods share the same subroutine to interpolate a computed solution from one mesh to another. The total number of the time steps are also the same for both methods. Therefore in Figure \ref{fig:moving_mesh_time}, we only show the relative average CPU time (i.e. relative to that of the full order method) to invert the linear system for one step time marching. With $N_x=1000$, it is observed that with $5$ to $25$ basis functions, it takes the LP method around $36\%$ of the time of the full order method to invert the resulting linear system.}

\pzcreview{}

%%%%%%%%%%%%%%%%%%%%%%%%%%%%%%%%%%%%%%%%%%%%%%%%%%%%%%%%%%%%%%%%%%%%%%%%%%%%%%%%%
% 2D linear advection equation
%%%%%%%%%%%%%%%%%%%%%%%%%%%%%%%%%%%%%%%%%%%%%%%%%%%%%%%%%%%%%%%%%%%%%%%%%%%%%%%%%
\subsection{\review{2D linear advection equation}\label{sec:2d_advection}}
\review{In this example, we consider a 2D linear advection equation:
\begin{align}
u_t + b\cdot \nabla u = 0, \quad (x,y)\in [0,1]\times[0,1]
\end{align}
with $b = (1.0,2.0)^T$ and the zero initial condition. 
The boundary condition is given at the inflow boundary $\Gamma_{\textrm{inflow}}=\{(x,y)\in [0,1]\times[0,1]: xy=0\}$ as below,
\begin{align}
    u(x,y)|_{(x,y)\in\Gamma_{\textrm{inflow}}} = \begin{cases}
              1+\cos(8\pi x),\quad 0.125<x<0.625,\\
              0,\quad\text{otherwise}.
              \end{cases}
\end{align} 
The computational domain is partitioned by a triangular mesh $\mathcal{T}_h=\{\mathcal{T}_i\}_{i=1}^{N_x}$. 
In space, we apply the upwind discontinuous Galerkin (DG) method: seeking 
$u_h\in V_h:=\bigoplus_{\mathcal{T}_i\in\mathcal{T}_h} \mathcal{P}^1(\mathcal{T}_i)$ satisfying 
\begin{align}
  \sum_{\mathcal{T}_i\in\mathcal{T}_h}\left( \int_{\mathcal{T}_i} \partial_t u_h v - (b \cdot \nabla v)u_h  dx+ \int_{\partial \mathcal{T}_i} \widehat{b\cdot \mathbf{n}_i u_h} vds\right) = 0, \quad \forall v \in V_h.
\end{align}
Here, $\mathcal{P}^1(\mathcal{T}_i)$ is the linear polynomial space on $\mathcal{T}_i$, $\mathbf{n}_i$ is the outward \review{unit normal of $\mathcal{T}_i$ along its boundary $\partial\mathcal{T}_i$, and $\widehat{b\cdot \mathbf{n}_i u_h}$ 
%$\hat{u}$ 
is the upwind flux that is based on the upwind mechanism and defined as
$\widehat{b\cdot \mathbf{n}_i u_h}=b\cdot\mathbf{n}_i \textrm{lim}_{\epsilon\rightarrow 0+}u_h(x-(b\cdot \mathbf{n}_i)\mathbf{n}_i\epsilon).$
}
In time, we apply the trapezoidal method \review{that is second order accurate}. The DG code is implemented with the python interface of the c++ finite element library NGSolve \cite{schoberl2014c++}.
}

\begin{figure}[h!]
\centering
%%%%%%%%%%%%%%%%%%%%%%%%%
\begin{subfigure}{.48\textwidth}
  \centering
\includegraphics[width=\textwidth,trim={5cm 5.34cm 3.25cm 3.5cm},clip]{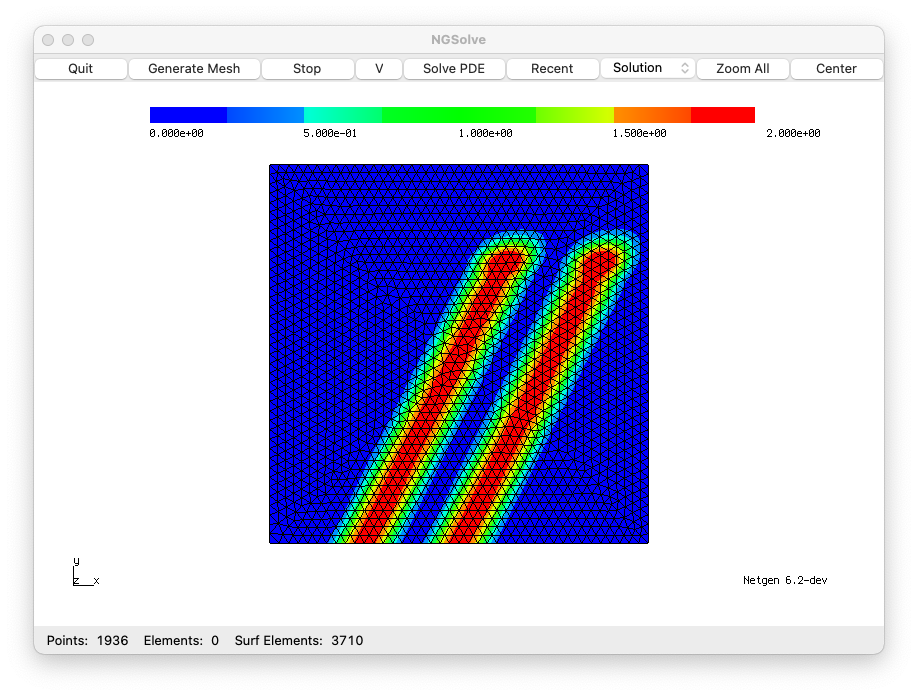}
  \caption{Full order solution
  %.
  \label{fig:2d_advection_full}}
\end{subfigure}
%%%%%%%%%%%%%%%%%%%%%%%%%
\begin{subfigure}{.48\textwidth}
  \centering
\includegraphics[width=\textwidth,trim={5cm 5.34cm 3.25cm 3.5cm},clip]{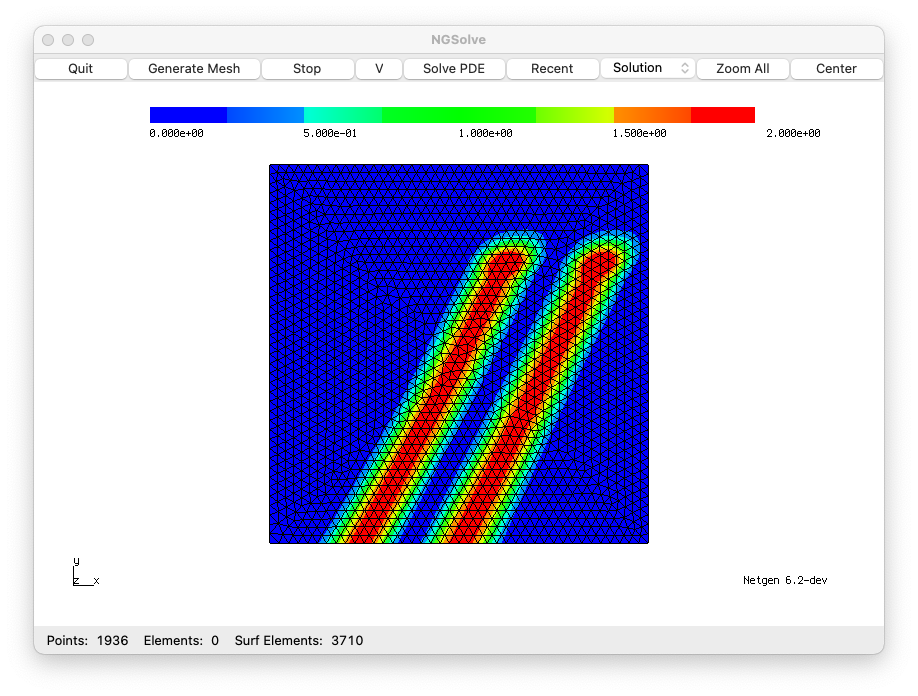}
  \caption{LP solution with $10$ basis functions}
\end{subfigure}
%%%%%%%%%%%%%%%%%%%%%%%%%
\begin{subfigure}{.48\textwidth}
  \centering
\includegraphics[width=\textwidth,trim={5cm 5.34cm 3.25cm 3.5cm},clip]{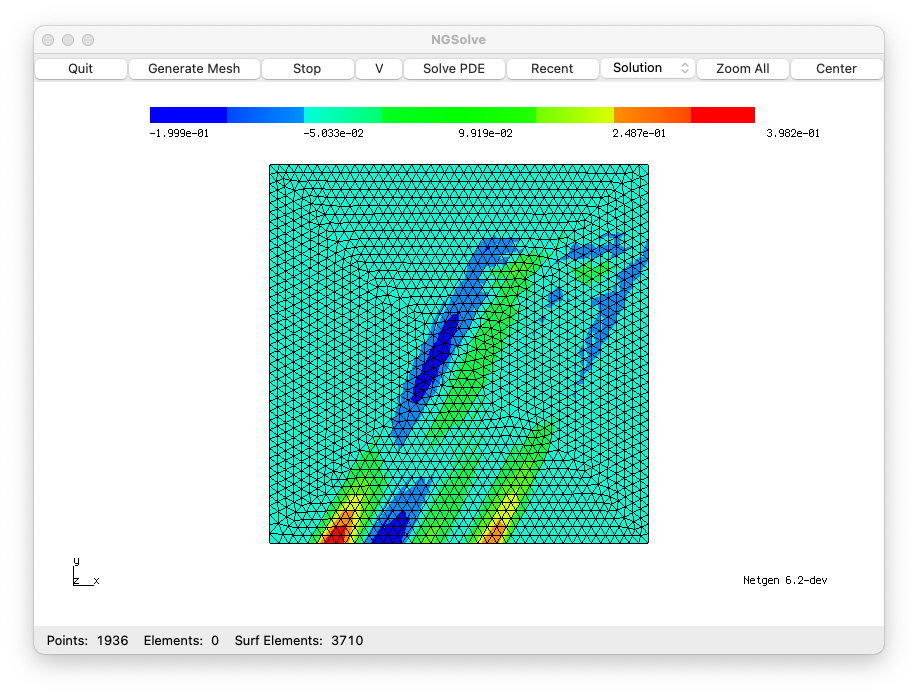}
  \caption{Error of the LP method with $10$ basis functions}
\end{subfigure}
%%%%%%%%%%%%%%%%%%%%%%%%%
\begin{subfigure}{.48\textwidth}
  \centering
\includegraphics[width=\textwidth,trim={5cm 5.34cm 3.25cm 3.5cm},clip]{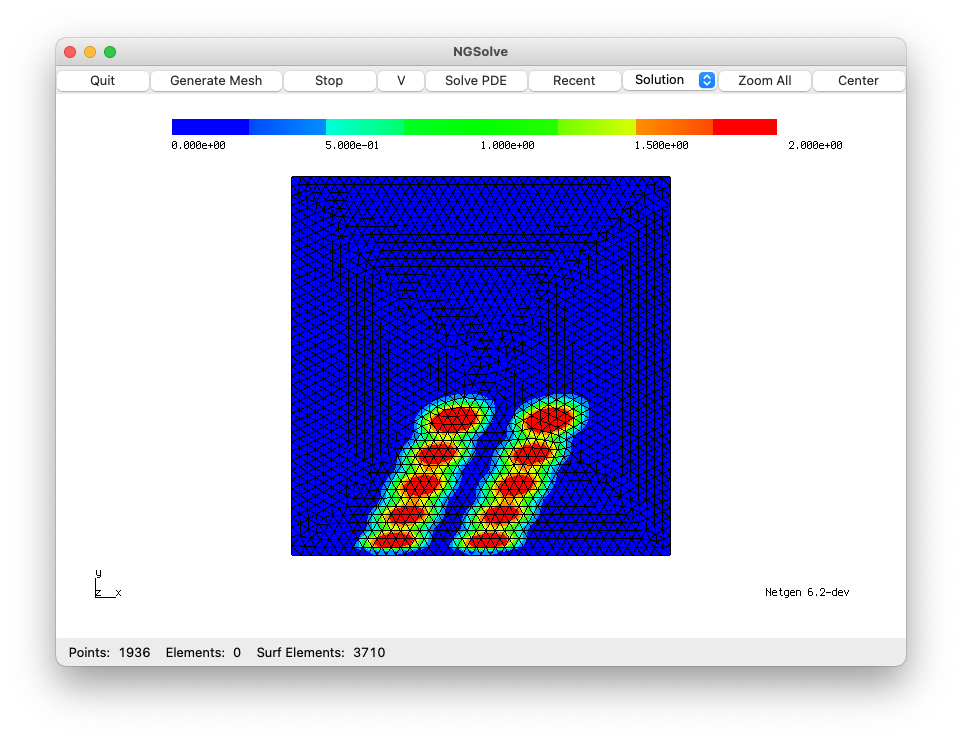}
  \caption{POD solution with $10$ basis functions}
\end{subfigure}
%%%%%%%%%%%%%%%%%%%%%%%%%
\caption{\review{The full order solution, {the LP solution and its error with 10 basis functions, and the POD solution with 10 basis functions, all} 
 at $t=0.4$ for the 2D linear advection equation in Section \ref{sec:2d_advection}.\label{fig:2d_advection_sol}}}
\end{figure}

\review{
The mesh {in our test consists of $N_x=3710$ elements with the maximal edge length 
$h_{\textrm{max}}=0.025$ (see Figure \ref{fig:2d_advection_full}),} and the time step size is set as \review{$\Dt=\frac{0.5 h_{\textrm{max}}}{\| b\|_2}$}. The training data is the nodal values of the full order solution for $t\in[0,0.2]$ {($37$ time snapshots)} and we predict the solution from $t=0.2$ to $t=0.4$. We find it much easier to train the neural network using nodal values instead of modal coefficients as the degrees of freedom. Since NGSolve implements the DG method based on a modal basis instead of a nodal basis, this leads to extra cost to covert between the nodal values and modal coefficients of the solution. }

\begin{figure}[h!]
\centering
%%%%%%%%%%%%%%%%%%%%%%%%%
\begin{subfigure}{.45\textwidth}
  \centering
  \includegraphics[width=\textwidth]{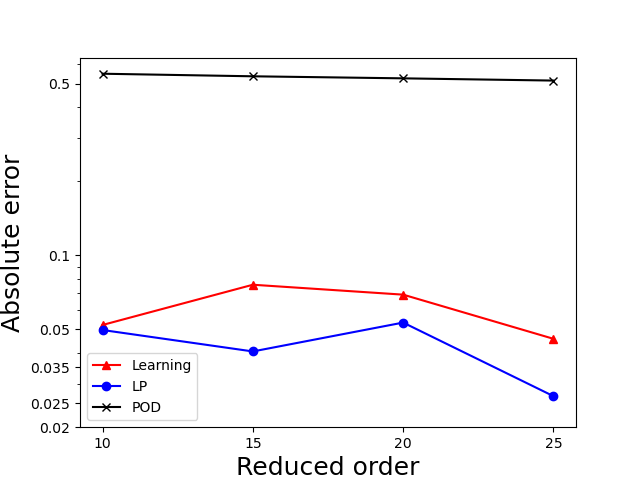}
  \caption{\pzcreview{Absolute error at $t=0.4$
  %.
  }\label{fig:2d_advection_error} }
\end{subfigure}
\begin{subfigure}{.45\textwidth}
  \centering
  \includegraphics[width=\textwidth]{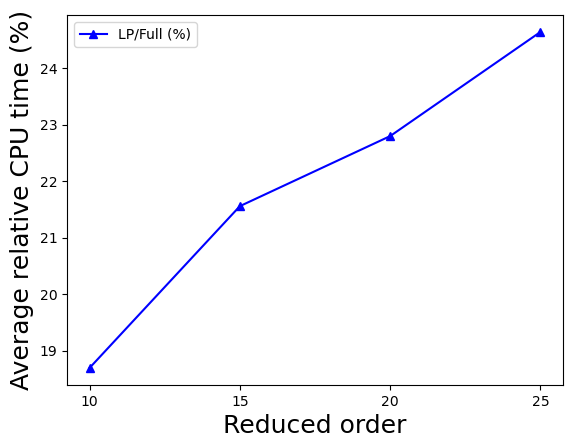}
    \caption{\reviewone{Average relative CPU time
    }\label{fig:2d_advection_relative_time}}
\end{subfigure}
%%%%%%%%%%%%%%%%%%%%%%%%%
\caption{\review{The absolute error at $t=0.4$ and the average relative CPU time for the 2D linear advection equation in Section \ref{sec:2d_advection}.\label{fig:2d_advection_error_time}}}
\end{figure}

\review{In Figure \ref{fig:2d_advection_sol}, we present the full order solution, the LP and POD solutions with $r=10$ at $t=0.4$, as well as the error in the LP solution. It is clear }
\pzcreview{that the LP method matches the full order solution much better than the POD method. In Figure \ref{fig:2d_advection_error_time}, we present the absolute error at $t=0.4$ for different methods and the relative CPU time of the LP method. We observe that the LP method always leads to the smallest
%least 
error, and it only needs less than $25\%$ CPU time of the full order solve even with the extra online overhead to convert nodal values to modal coefficients.} 

\begin{figure}[h!]
\centering
%%%%%%%%%%%%%%%%%%%%%%%%%
\begin{subfigure}{.45\textwidth}
  \centering
  \includegraphics[width=\textwidth]{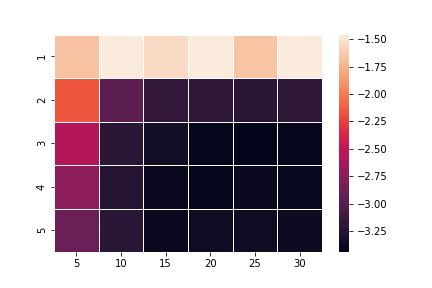}
  \caption{\pzcreview{Final value of the loss function of trained neural network
  %.
  }\label{fig:2d_advection_error} }
\end{subfigure}
\begin{subfigure}{.45\textwidth}
  \centering
  \includegraphics[width=\textwidth]{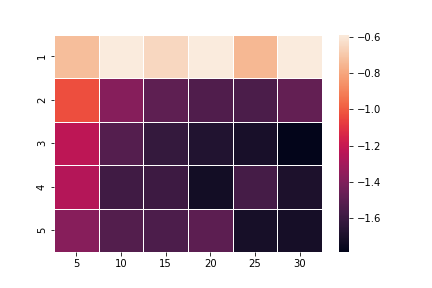}
    \caption{\pzcreview{The $L_2$ prediction error of the trained neural network at $t=0.2$ under $\log$-scale
    %.
    }\label{fig:2d_advection_relative_time}}
\end{subfigure}
%%%%%%%%%%%%%%%%%%%%%%%%%
\caption{\reviewthree{$x$-axis: width of the hidden layers; $y$-axis: depth of the hideen layers; color:   the value of the loss function/error under log scale. \label{fig:nn_heatmap}}}
\end{figure}

\review{We also use this example to illustrate
%show 
that the proposed method is relatively insensitive to the width and the depth of the hidden layer of the neural network. In Figure \ref{fig:nn_heatmap}, we present the final value of the loss function and the prediction error of the trained neural network at $t=0.2$ with different width and depth of the hidden layer.
The two neural networks that represent the basis and the coefficients share the same width and the depth.  These neural networks are also trained with the solution snapshots for $t\in[0,2]$, and the learning rate is $10^{-3}$ and the number of epochs is $30$. For each pair of depth and width, the value of the loss function and the prediction error are the average value of $5$ different random initialization of the neural network.  We observe that when the depth is not less than $2$ and the width is not less than $10$, the final value of the loss function and the prediction error of the trained neural network at $t=0.2$ do not change much with respective to the depth and the width.} 

\begin{figure}[h!]
\centering
%%%%%%%%%%%%%%%%%%%%%%%%%
\begin{subfigure}{.45\textwidth}
  \centering
  \includegraphics[width=\textwidth]{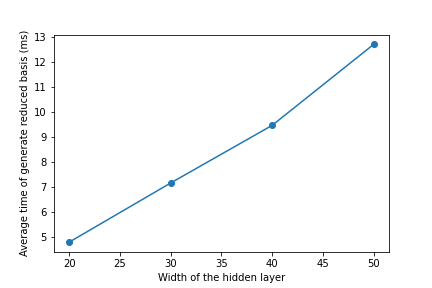}
  \caption{\reviewone{Depth of the hidden layer, $3$
  }\label{fig:time_vs_width} }
\end{subfigure}
\begin{subfigure}{.45\textwidth}
  \centering
  \includegraphics[width=\textwidth]{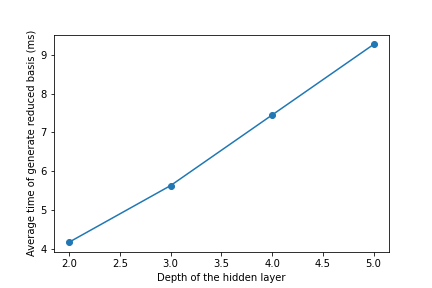}
    \caption{\reviewone{Width of the hidden layer, $25$
    }\label{fig:time_vs_depth}}
\end{subfigure}
%%%%%%%%%%%%%%%%%%%%%%%%%
\caption{\reviewone{$x$-axis: width/depth of the hidden layers; $y$-axis: average CPU time to generate $10$ reduced basis functions online. \label{fig:time_vs_nn}}}
\end{figure}

\reviewone{The last test for this example is to investigate the relation between the cost to generate the reduced basis online and the architecture of the neural networks. In Figure \ref{fig:time_vs_nn}, we compute the average CPU time of $100$ runs to generate reduced basis for $t=0.4$ and $r=10$. In \ref{fig:time_vs_width}, we fix the depth of the hidden layer as $3$ and vary the width of the hidden layer. In \ref{fig:time_vs_depth}, we fix the width of the hidden layer as $25$ and vary the depth. We observe that the computational cost of generating $10$ reduced basis functions is roughly proportional to the width and the depth of the hidden layer. We also would like to point out that, as to be shown in Section \ref{sec:advection_hetero}, the computational time of generating the basis will be less than solving the reduced order projection equation when the mesh is refined enough.}

%%%%%%%%%%%%%%%%%%%%%%%%%%%%%%%%%%%%%%%%%%%%%%%%%%%%%%%%%%%

%%%%%%%%%%%%%%%%%%%%%%%%%%%%%%%%%%%%%%%%%%%%%%%%%%

%%%%%%%%%%%%%%%%%%%%%%%%%%%%%%%%%%%%%%%%%%%%%%%%%%

%%%%%%%%%%%%%%%%%%%%%%%%%%%%%%%%%%%%%%%%%%%%%%%%%%

%%%%%%%%%%%%%%%%%%%%%%%%%%%%%%%%%%%%%%%%%%%%%%%%%%%%%%%%%%%%%%%%%%%%%%%%%%%%%%%%%
% Future prediction Burgers equation
%%%%%%%%%%%%%%%%%%%%%%%%%%%%%%%%%%%%%%%%%%%%%%%%%%%%%%%%%%%%%%%%%%%%%%%%%%%%%%%%%
\subsection{1D linear advection equation with  $(x,\mu)$-dependent
\pzcrev{inhomogeneous} media\label{sec:advection_hetero}}

In this example, we consider the following 1D linear advection equation
\begin{align}\label{eq:p_dependent_advection}
    u_t+(c(x;\mu)u)_x=0
\end{align}
in \lli{an}
\pzcrev{inhomogeneous} media where the propagation speed  $c(x,\mu)$ is $(x,
\mu)$-dependent, given as 
\begin{align}
    c(x;\mu) = 1.25+\mu_1\sin( \mu_2\pi x ),\quad \lli{\mu= (\mu_1,\mu_2)}\in X_{\MP}=[0,0.5]\times[2,4].
\end{align}
The initial condition is
\begin{align}
    u(x,0) = \begin{cases}
              0.5+0.5\cos(5\pi(x-0.25)), \quad 0.05\leq x\leq 0.45,\\
              0,\quad\text{otherwise\mw{,}}
             \end{cases}
\end{align}
together with the zero boundary conditions.
The \mw{spatial} computational domain is \pzcreview{$[0,2.5]$}, and a uniform mesh with \pzcreview{$500$} elements is used. To control the numerical dissipation of the full order model, we use the second-order Crank-Nicolson time discretization, and the numerical flux is chosen as an  upwind-biased flux:
\begin{align}
    f^n_{j+\frac{1}{2}}=\frac{3}{4}c_ju_j^n+\frac{1}{4}c_{j+1}u_{j+1}^n
    \label{eq:advection_hetero_flux}
\end{align}
with $c_j = c(x_j;\mu)>0$. \reviewtwo{The time step size is $\Dt={2\Dx}$.}
%h$.}

Associated with the model parameter $\mu=(\mu_1,\mu_2)$, the training set  is \mw{taken to be} $\{\mu_1^{(i)}=\frac{i}{20}\}_{i=0}^{10}\times\{\mu_2^{(j)}=2+\frac{j}{10}\}_{j=0}^{20}$, and the \mw{testing} set  is \mw{constituted by }  $20$ pairs of $(\mu_1,\mu_2)$ \pzcreview{randomly} sampled from \review{a uniform distribution on the parameter space} $X_\MP$.
For both the training and the testing,  \pzcreview{$t\in[0,1]$} \pzcnew{($231\times101$ time snapshots)} is considered. 

\begin{figure}[h!]
\centering
\begin{subfigure}{.45\textwidth}
  \centering
  \includegraphics[width=\textwidth]{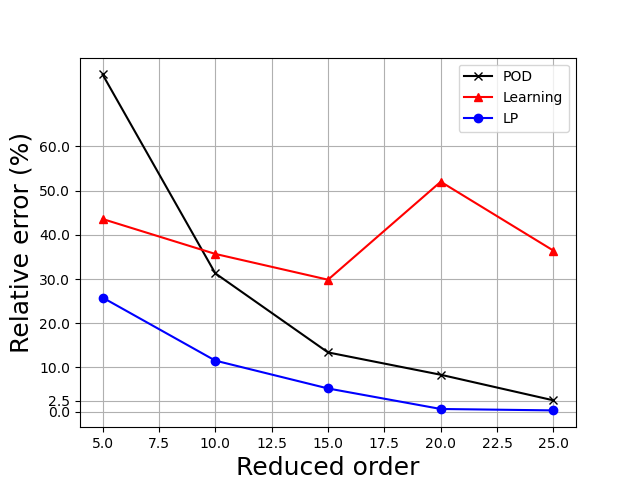}
  \caption{\review{Average relative error}}
  %Comparison between  pure learning, POD and the LP method}
\end{subfigure}
\begin{subfigure}{.45\textwidth}
  \centering
  \includegraphics[width=\textwidth]{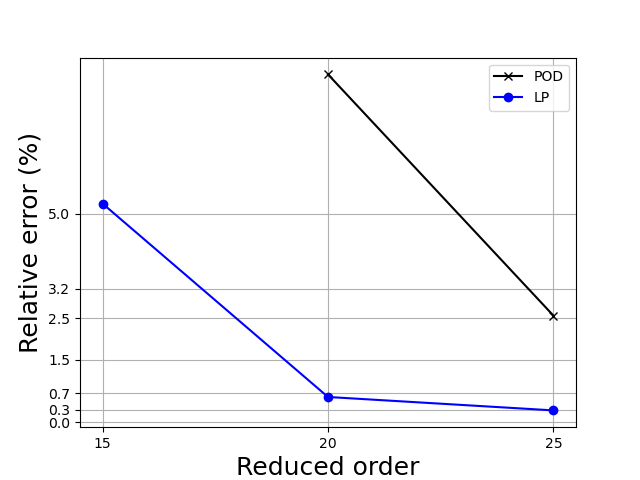}
    \caption{Zoomed-in picture}
\end{subfigure}
\caption{$x$-axis: reduced order $r$; $y$-axis: the average relative error for the test set of the linear advection equation with  $(x,\mu)$-dependent
\pzcrev{inhomogeneous} media in \ref{sec:advection_hetero}.\label{fig:advection_hetero_error}}
\end{figure}

In Figure \ref{fig:advection_hetero_error}, we report \mw{the average relative testing error of the proposed LP method, the learning method, and the standard POD method versus the reduced order $r$}. We observe that the LP method \mw{consistently entails} the least error \mw{for varied $r$}.  When \pzcreview{$r=25$}, the error of the LP method is smaller than $0.3\%$. Possibly due to the overfitting, the error of the pure learning-based method does not decay as the number of basis functions increases. The projection step in the online stage reduces the generalization error dramatically, as the PDE information is utilized through this step.

\begin{figure}[h!]
\centering
\begin{subfigure}{.45\textwidth}
  \centering
  \includegraphics[width=\textwidth]{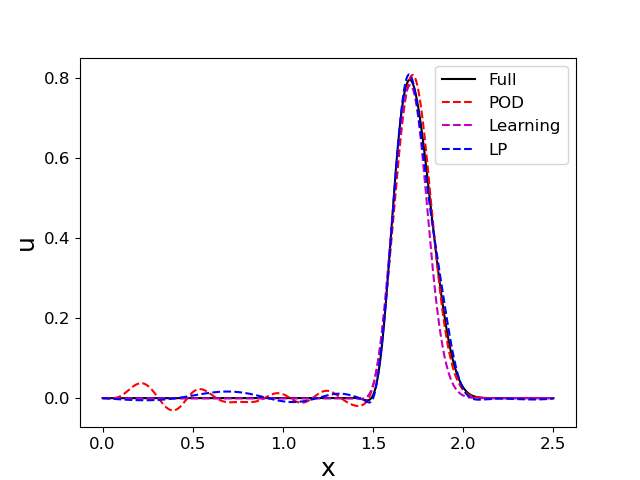}
  \caption{\pzcreview{$15$ basis functions}}
\end{subfigure}
%%%%%%%%%%%%%%%%%%%%%%%%%
\begin{subfigure}{.45\textwidth}
  \centering
  \includegraphics[width=\textwidth]{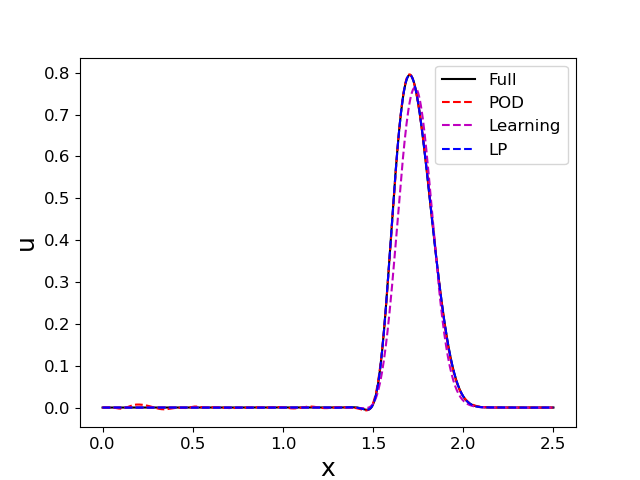}
  \caption{\pzcreview{$25$ basis functions}}
\end{subfigure}
%%%%%%%%%%%%%%%%%%%%%%%%%
\caption{\pzcreview{Predicted solutions of
the linear advection equation with  $(x,\mu)$-dependent
\pzcrev{inhomogeneous} media in Section \ref{sec:advection_hetero} 
at $(\mu_1,\mu_2)=(0.23,3.56)$ and
 $t=1$.}\label{fig:advection_hetero_solution}}
\end{figure}

For the \mw{testing} sample $(\mu_1,\mu_2)=(0.23,3.56)$, we present in  Figure \ref{fig:advection_hetero_solution} the solutions with different reduced order \lli{$r=15,25$} at $t=2$. \pzcrev{With $15$ basis functions, the LP solution matches the full order solution better than the POD method, though}
%small
\lli{some non-physical 
oscillation
%s 
exists. Additional study can be found in  Appendix \ref{sec:appendix} to understand the origin of the oscillation.}
%(see the appendix \ref{sec:appendix} for the reason).}
With $25$ basis functions, the LP solution visually coincides with the full order solution.

\begin{figure}[h!]
\centering
  \includegraphics[height=0.3\textwidth]{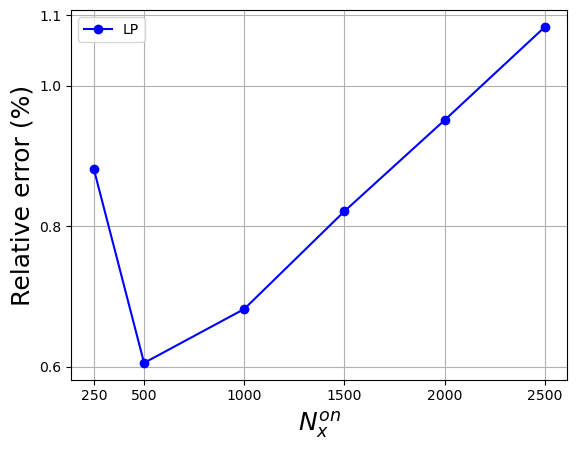}
   \includegraphics[height=0.3\textwidth]{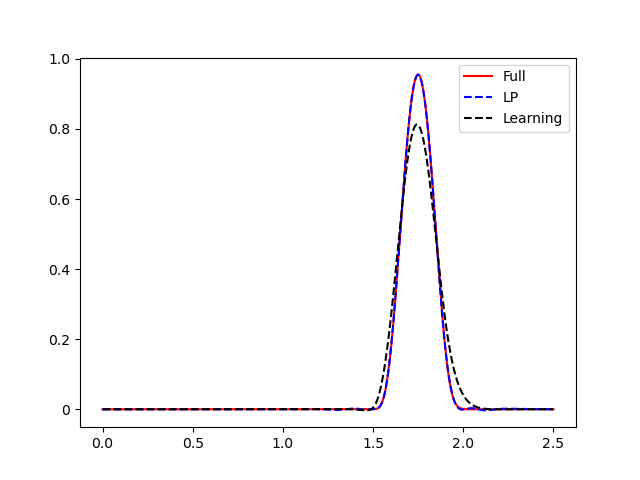}
\caption{\review{Linear advection equation with  $(x,\mu)$-dependent \pzcrev{inhomogeneous} media in Section \ref{sec:advection_hetero}.  Left:} mean relative error \mw{of} the test set for different \mw{online} grid resolutions ($N_x^{\textrm{on}}$). \pzcreview{Right: predicted solution at $\mu=(0.23,3.56)$ and $t=1$ with online mesh $N_x^{\textrm{on}}=2500$.} 
$20$ basis functions are used and are trained \mw{offline} with data subject to constant resolution ($N_x^{\textrm{off}}=500$).
\label{fig:error_grid_advection_hetero}}
\end{figure}

Our next test aims at showing  the flexibility of the LP method when the online spatial meshes are different from the offline ones. Our method avoids extra efforts which are needed by conventional projection based ROMs  to accomplish similar tasks \cite{grassle2018pod,yano2019discontinuous,yano2020goal}. Using the neural network trained with \mw{data generated on a uniform mesh with} \pzcreview{$N_x^{\textrm{off}}=500$} spatial elements offline, we  evaluate basis function \mw{values} on uniform meshes with different $N_x^{\textrm{on}}$ in the online computation. The relative error for different grid resolutions online with $r=20$  is presented in the left picture of Figure \ref{fig:error_grid_advection_hetero}, and it is reasonably small. \pzcreview{In the right picture of Figure \ref{fig:error_grid_advection_hetero}, we present the predicted solution at $\mu=(0.23,3.56)$ and $t=1$ with $N_x^{\textrm{on}}=2500$. The full order solution with $N_x^{\textrm{off}}=500$ is more dissipative than the full order solution with $N_x^{\textrm{on}}=2500$,  but the LP solution still matches the full order solution with $N_x^{\textrm{on}}=2500$ well, while the pure learning-based method provides a more dissipative solution compared with the full order solution. }

\begin{figure}
\begin{subfigure}{.45\textwidth}
  \centering
  \includegraphics[width=\textwidth]{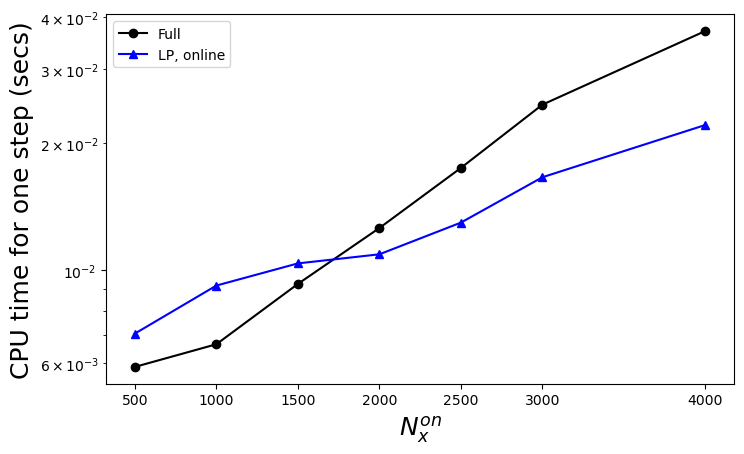}
\end{subfigure}
\begin{subfigure}{.45\textwidth}
  \centering
  \includegraphics[width=\textwidth]{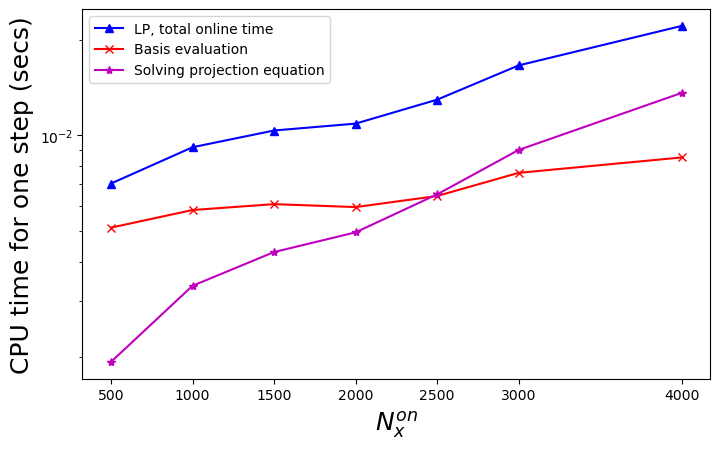}
\end{subfigure}
  \caption{The average CPU time for one time step  as a function of different \mw{online resolution $N_x^{\textrm{on}}$}. Left: full order model \mw{v.s.} the online computation of the LP method.  Right: 
  online basis evaluation, the online \pzc{solve of projection equation} of the LP method and the total time of the online computation of the LP method.
  \label{fig:advection_hetero_cpu_time}}
\end{figure}

The final test we carry out is \mw{intended} to examine the computational efficiency. In Figure \ref{fig:advection_hetero_cpu_time}, for the test sample $(\mu_1,\mu_2)=(0.23,3.56)$, we compare the CPU time of the \mw{one step \pzc{online} computation of the}  LP method  with $20$ basis functions and that of the full order model for different grid resolutions online.  We also compare the times  for basis \mw{evaluation}  in \eqref{eq:online_basisGen} and for solving the projection equation  \eqref{eq:online_projection}  during the online stage of the LP method (see Algorithm \ref{alg:onlinee}). For different  \mw{resolutions online} ($N_x^{\textrm{on}}$), the reduced-order basis is always trained with \pzcreview{$N_x^{\textrm{off}}=500$} elements in the offline stage, and the solutions predicted by our method 
always visually match the full order solution reasonably well. As $N_x^{\textrm{on}}$ grows, the CPU time of the basis evaluation grows much slower compared with that for solving the projection equation. For small $N_x^{\textrm{on}}$, the ROM spends most of the time on the basis generation. As $N_x^{\textrm{on}}$ grows, the computational time of the full order model grows much faster than that of the proposed LP method.

%%%%%%%%%%%%%%%%%%%%%%%%%
%%%%%%%%%%%%%%%%%%%%%%%%%%%%%%%%%%%%%%

%%%%%%%%%%%%%%%%%%%%%%%%%%%%%%%%%%%%%%%%%%%%%%%%%%%%%
%%%%%%%%%%%%%%%%%%%%%%%%%%%%%%%%%%%%%%%%%%%%%%%%%%%%%
%%%%%%%%%%%%%%%%%%%%%%%%%%%%%%%%%%%%%%%%%%%%%%%%%%%%%

%%%%%%%%%%%%%%%%%%%%%%%%%%%%%%%%%%%%%%%%%%%%%%%%%%%%%%%%%%%%%%%%%%%%%%%%%%%%%%%%%
% Burgers equation
%%%%%%%%%%%%%%%%%%%%%%%%%%%%%%%%%%%%%%%%%%%%%%%%%%%%%%%%%%%%%%%%%%%%%%%%%%%%%%%%%
\subsection{1D Burgers Equation}
We consider the Burgers equation:
\begin{align}
    u_t + \left(\frac{1}{2}u^2\right)_x=0,
\end{align}
and \review{perform two sets of experiments. In the first one, we consider a Riemann problem with a parameter in the initial data, and the solutions will be predicted for a test set of parameter and for future times. In the second experiment, we start with a smooth initial profile, and shock discontinuity will form later. We will predict the solutions at future times based on the full order solution on $[0,t_0]$, over which the shock structure may or may not have been ``seen'' by the offline training algorithm.  The first order method in \eqref{eq:rev:med1} will be applied, with the numerical flux chosen as:}
\begin{align}
    f_{j+\half}^n:= \frac{1}{2}\left(f(u_{j+1}^n)+f(u_{j}^{n})\right)-
    \frac{\Dx}{2\Dt}(u_{j+1}^n-u_j^n).
\end{align}

%%%%%%%%%%%%%%%%%%%%%%%%%%%%%%%%%%%%%%%%%%%%%%%%%%%%%%%%%%%%%%%%%%%%%%%%%%%%%%%%%
\subsubsection{Parameter-dependent Riemann problem  \label{sec:burgers_initial_parameter}}
%Now, we 
We consider the computational domain $[0,8]$ with a parameter-dependent initial condition:
\begin{align}
    u(x,0) = \begin{cases}
             1.+\mu,\quad
             %\text{if}\;
             \pzcreview{x\in[0.5,0.75]},\\
             0.,\quad\text{otherwise.}
             \end{cases}
\end{align}
Here, $\mu\in X_\MP=[0,1]$. The solution  has a rarefaction fan and a moving shock. We use a uniform mesh with \pzcreview{$N_x=400$}. The time step size is set \mw{to be} \pzcreview{$\Dt=\Dx$}.

The training set is taken as $\{0,0.1,0.2,\dots,1\}\subset X_\MP$ with \pzcreview{$0\leq t\leq 0.25$} \pzcnew{($11\times 101$ snapshots)}. The test set is taken as $\{0.05\pm0.01,0.15\pm0.01,\dots,0.95\pm0.01\}\subset X_\MP$, and we predict the solution from $t=0$ to \pzcreview{$t=0.5$}. \mw{In} the offline stage of this problem, we train the neural network with \pzcreview{$150$} epochs.

\begin{figure}[h!]
\centering
\begin{subfigure}{.45\textwidth}
  \centering
  \includegraphics[width=\textwidth]{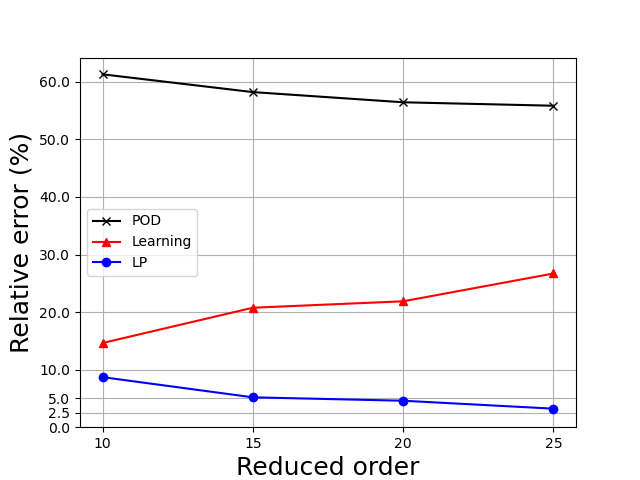}
  \caption{\review{Average relative error}}
  %Comparison between pure learning, POD and the \mw{LP} method}
\end{subfigure}
\begin{subfigure}{.45\textwidth}
  \centering
  \includegraphics[width=\textwidth]{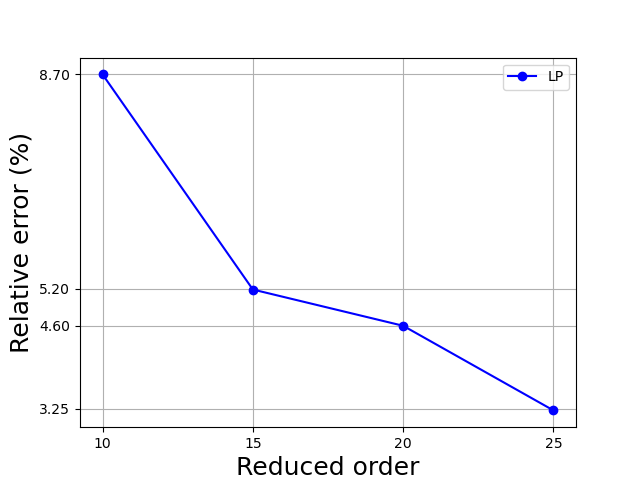}
    \caption{Zoomed-in picture}
\end{subfigure}
\caption{$x$-axis: reduced order $r$; 
$y$-axis: the average relative error for the test set of the Burgers equation (parameter dependent initial profile) in Section \ref{sec:burgers_initial_parameter}.\label{fig:burgers_initial_parameter_error_vs_order}}
\end{figure}

In Figure \ref{fig:burgers_initial_parameter_error_vs_order}, the average relative errors for the test set are presented \review{corresponding to different reduced order $r$.} \mw{Overall, we} observe that for both the POD method and \mw{the} LP method, the error decays as {$r$ } grows. However, this is not the case  for the purely learning-based method. Moreover, the error of the LP method is much smaller than \mw{that of} the other two methods.

\begin{figure}[h!]
\centering
%%%%%%%%%%%%%%%%%%%%%%%%%
\begin{subfigure}{.45\textwidth}
  \centering
  \includegraphics[width=\textwidth]{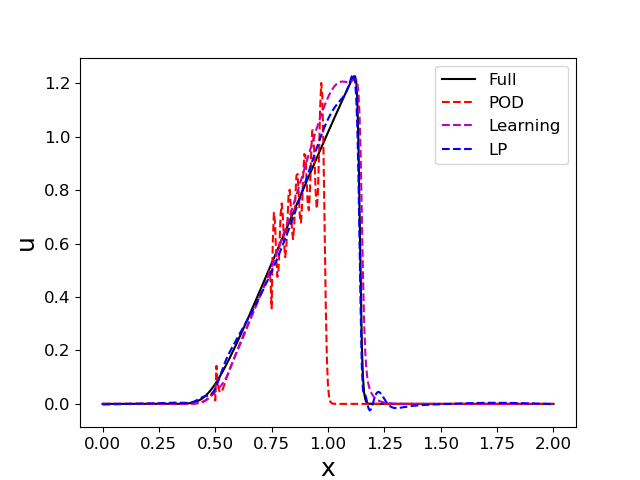}
  \caption{$15$ basis functions}
\end{subfigure}
%%%%%%%%%%%%%%%%%%%%%%%%%
\begin{subfigure}{.45\textwidth}
  \centering
  \includegraphics[width=\textwidth]{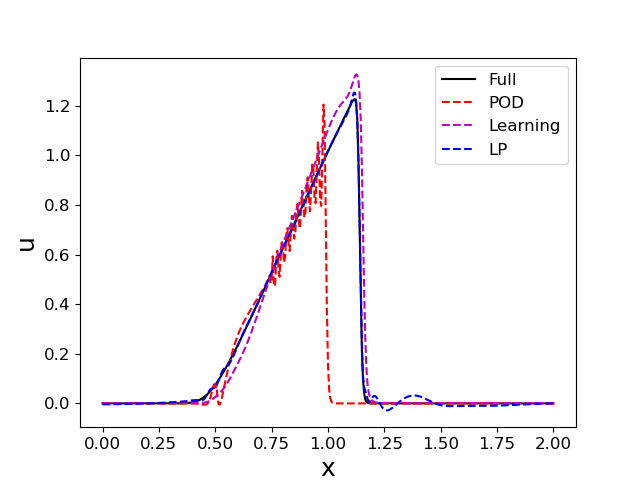}
  \caption{$25$ basis functions}
\end{subfigure}
%%%%%%%%%%%%%%%%%%%%%%%%%
\caption{\pzcreview{Predicted solutions 
for the Burgers equation (parameter dependent initial profile) in Section \ref{sec:burgers_initial_parameter} with $\mu=0.76$ and $t=0.5$.\label{fig:burgers_initial_parameter_solution}}}
\end{figure}

In Figure \ref{fig:burgers_initial_parameter_solution}, we present the predicted solutions for $\mu=0.76$ at \review{$t=0.5$} with \review{$r=15, 25$}.  We observe that the LP method performs better than the other two methods in matching the full order solution, especially in capturing  the shock location. There are small oscillations \review{near the shock location} \review{in the LP solution }  .

\begin{figure}[h!]
\centering
\begin{subfigure}{.45\textwidth}
  \centering
  \includegraphics[width=\textwidth]{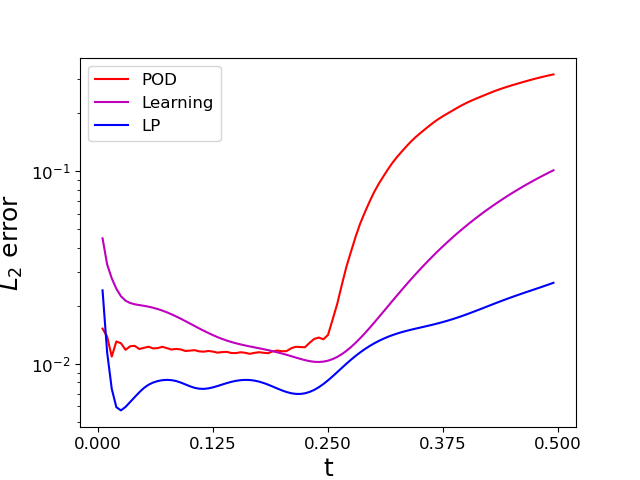}
  \caption{$15$ basis functions}
\end{subfigure}
%%%%%%%%%%%%%%%%%%%%%%%%%
\begin{subfigure}{.45\textwidth}
  \centering
  \includegraphics[width=\textwidth]{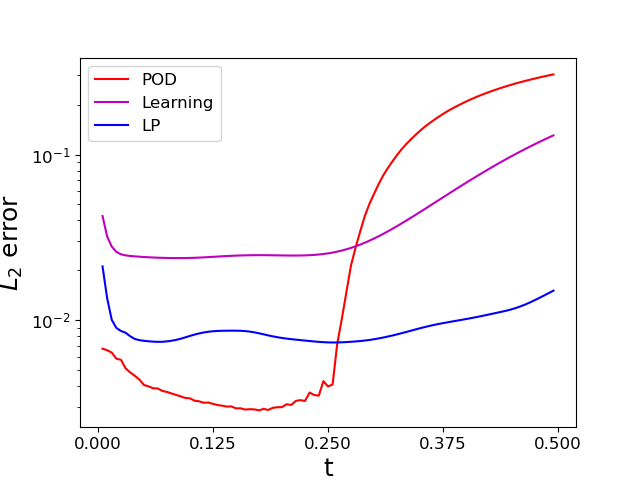}
  \caption{$25$ basis functions}
\end{subfigure}
%%%%%%%%%%%%%%%%%%%%%%%%%
%%%%%%%%%%%%%%%%%%%%%%%%%
\caption{\pzcreview{The time evolution of the prediction error \mw{history} for the Burgers equation (parameter dependent initial profile) in Section \ref{sec:burgers_initial_parameter} with $\mu=0.76$.\label{fig:burgers_initial_parameter_error_history}}}
\end{figure}

In Figure \ref{fig:burgers_initial_parameter_error_history}, we present the time \mw{evolution} of the $L_2$ error in prediction for $\mu=0.76$ under the $log$-scale with \review{$r=15$ and $25$.
Among the LP method, the learning method, and the POD method, the LP method displays the slowest growth in its error over time.  When $r=25$, }
we observe that the error of the POD method is smaller than that of the LP method for \pzcreview{$t\leq0.25$}, then it grows \pzcreview{dramatically} fast and becomes much larger than the LP method. The reason why the POD method becomes \pzcreview{much} worse after \pzcreview{$t=0.25$} is that we only have the training data for \pzcreview{$t\leq 0.25$} in the offline stage. The LP method shows its advantages in providing more accurate and stable long time prediction.

\begin{figure}[h!]
\centering
%%%%%%%%%%%%%%%%%%%%%%%%%
  \centering
  \includegraphics[width=0.45\textwidth]{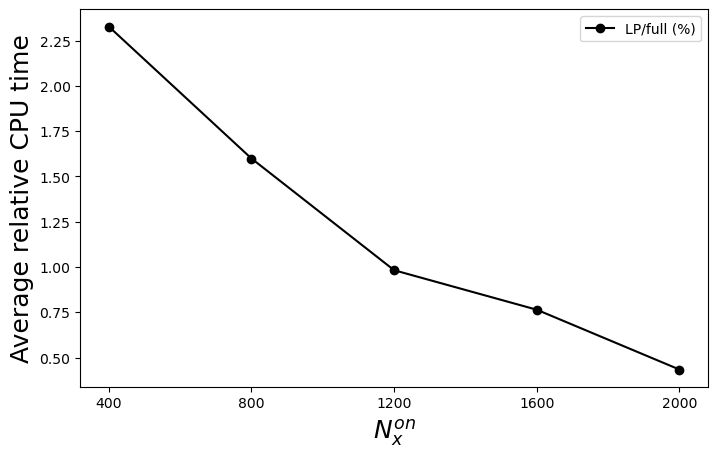}
\caption{\pzcreview{Relative computational time with different resolutions  in the online stage.
%$N_x^{\textrm{on}}$.
\label{fig:burgers_initial_parameter_cpu_time}}}
\end{figure}

\reviewmulti{In Figure \ref{fig:burgers_initial_parameter_cpu_time}, we present the relative CPU time with different number of mesh elements in the online stage $N_x^{\textrm{on}}$ and $r=20$. The neural network is trained with $N_x^{\textrm{off}}=400$. When the mesh is refined enough, the LP method does lead to smaller computational time compared with a full order solver. However, {when the mesh is refined with larger}  $N_x^{\textrm{on}}$, the LP method produces more numerical oscillations. It becomes important on much refined meshes to control non-physical oscillations of reduced order solutions and to improve numerical stability.}

%%%%%%%%%%%%%%%%%%%%%%%%%%%%%%%%%%%%%%%%%%%%%%%%%%%%%%%%%%%%%
%%%%%%%%%%%%%%%%%%%%%%%%%%%%%%%%%%%%%%%%%%%%%%%%%%%%%%%%%
%%%%%%%%%%%%%%%%%%%%%%%%%%%%%%%%%%%%%%%%%%%%%%%%%%%%%%%%%

%%%%%%%%%%%%%%%%%%%%%%%%%%%%%%%%%%%%%%%%%%%%%%%%%%%%%%%%%
%%%%%%%%%%%%%%%%%%%%%%%%%%%%%%%%%%%%%%%%%%%%%%%%%%%%%%%%%%%%%%%%%%%%%%%%%%%%%%%%%
\subsubsection{Future prediction with smooth initial profile and shock formation} \label{sec:burgers_prediction_smooth}
We consider the computational domain $[-1,1]$ and the following initial condition
%the smooth initial data
\begin{align}
    u(x,0) = \frac{1}{4}+\frac{1}{2}\sin(\pi x)
\end{align}
with periodic boundary conditions. Though with a smooth initial 
profile,
a shock will form at $t=t_\star=\frac{2}{{\pi}}$ and then moves as time evolves. We use a uniform mesh with $N_x=400$ and the time step size $\Dt=0.25\Dx$.

We consider two training sets, which are the full order numerical solutions for $t\in[0,t_0]$, with $t_0=1.1$ 
(\mw{when shock has  formed with \pzcnew{$881$ snapshots in the training set}}) 
and $t_0=0.6$
(\mw{when} shock \mw{has} not formed \mw{yet} with \pzcnew{$481$ snapshots in the training set}), respectively.
We predict the solution from $t=t_0$  to $t=1.6$. The neural network is trained with $100$ epochs. 
It is expected that the experiment with the second training data \mw{which} has not ``seen'' any shock \mw{to be} more challenging \mw{in} future prediction. 

\begin{figure}[h!]
\centering
%%%%%%%%%%%%%%%%%%%%%%%%%
\begin{subfigure}{.45\textwidth}
  \centering
  \includegraphics[width=\textwidth]{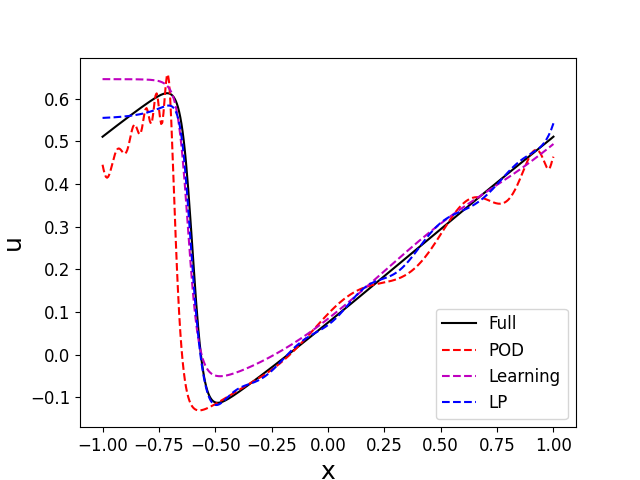}
  \caption{15 basis functions; training data: $t\in[0,1.1]$}  
 % trained with the data for $t\in[0,1.1]$}
\end{subfigure}
%%%%%%%%%%%%%%%%%%%%%%%%%
\begin{subfigure}{.45\textwidth}
  \centering
  \includegraphics[width=\textwidth]{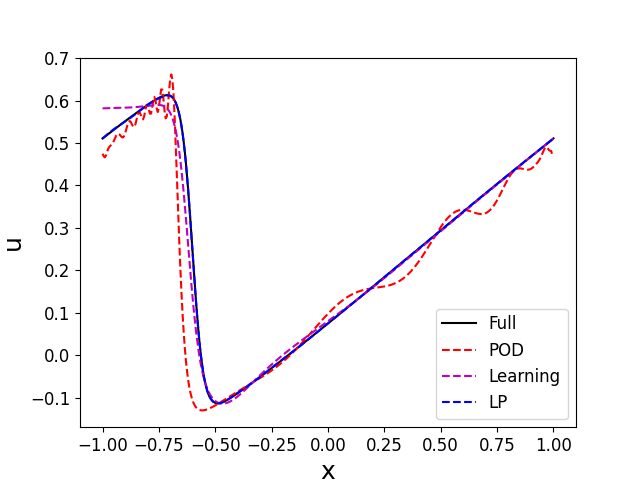}
  \caption{25 basis functions; training data: $t\in[0,1.1]$}
\end{subfigure}
%%%%%%%%%%%%%%%%%%%%%%%%%
\begin{subfigure}{.45\textwidth}
  \centering
  \includegraphics[width=\textwidth]{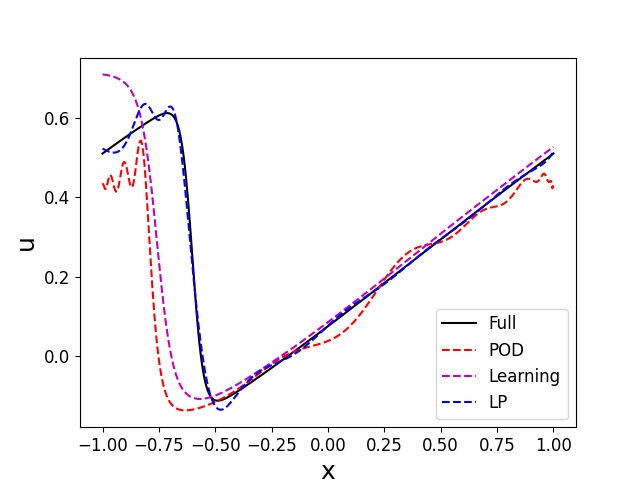}
  \caption{15 basis functions; training data: $t\in[0,0.6]$}
\end{subfigure}
%%%%%%%%%%%%%%%%%%%%%%%%%
\begin{subfigure}{.45\textwidth}
  \centering
  \includegraphics[width=\textwidth]{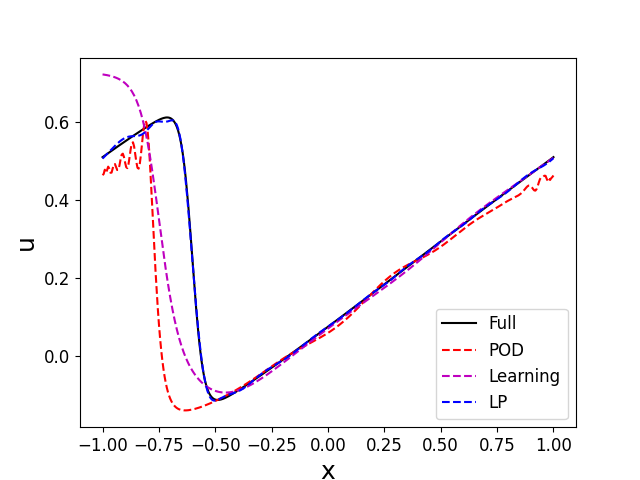}
  \caption{25 basis functions; training data: $t\in[0,0.6]$}
\end{subfigure}
%%%%%%%%%%%%%%%%%%%%%%%%%
\caption{Predicted  solutions of the Burgers equation  with smooth initial profile and shock formation in Section \ref{sec:burgers_prediction_smooth} at $t=1.6$. \label{fig:burgers_smooth_prediction_solution}}
\end{figure}

In Figure \ref{fig:burgers_smooth_prediction_solution}, we present the predicted reduced-order solutions at $t=1.6$. We observe that, with both $t_0=1.1$ and $t_0=0.6$, the solutions by the LP method always match the full order solution the best.   Even though the training data corresponding to $t\in[0,0.6]$ does not include any solutions with shock features, the LP method  captures the correct shock location in the predicted solutions.  Besides, with the training data from $t\in[0,0.6]$, the results of the LP method are relatively  more oscillatory. \mw{The} numerical oscillation is further suppressed with \mw{an} increase of the reduced order $r$.

\begin{figure}[h!]
\centering
\begin{subfigure}{.45\textwidth}
  \centering
  \includegraphics[width=\textwidth]{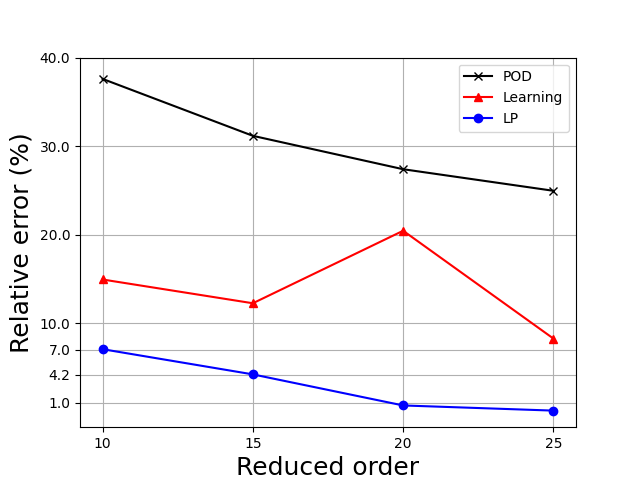}
  \caption{Comparison between pure learning, POD and the LP method, training data:  $t\in[0,1.1]$}
\end{subfigure}
\begin{subfigure}{.45\textwidth}
  \centering
  \includegraphics[width=\textwidth]{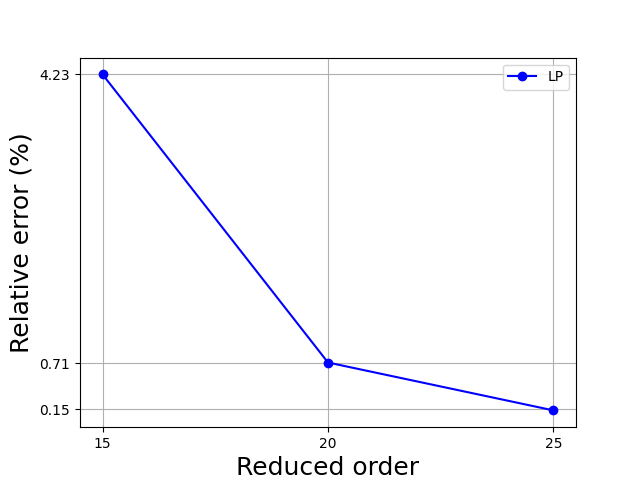}
    \caption{Zoomed-in picture, training data:  $t\in[0,1.1]$}
\end{subfigure}
\begin{subfigure}{.45\textwidth}
  \centering
  \includegraphics[width=\textwidth]{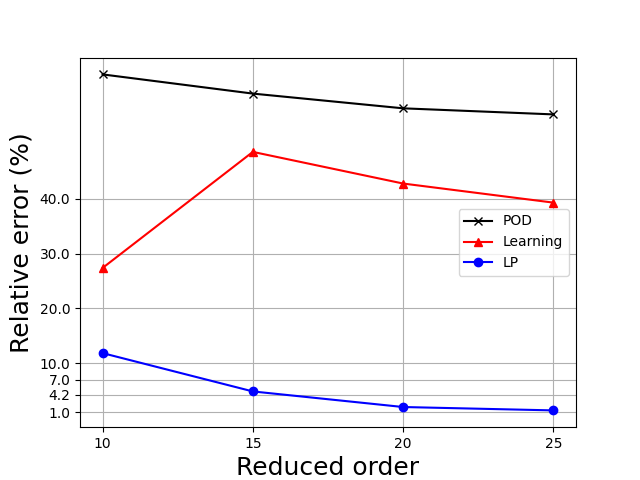}
  \caption{Comparison between pure learning, POD and the LP method, training data: $t\in[0,0.6]$}
\end{subfigure}
\begin{subfigure}{.45\textwidth}
  \centering
  \includegraphics[width=\textwidth]{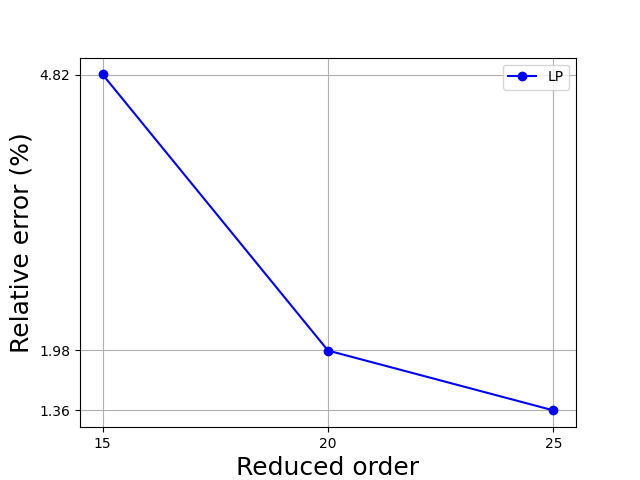}
    \caption{Zoomed-in picture, training data:  $t\in[0,0.6]$}
\end{subfigure}
\caption{$x$-axis: reduced order $r$;  $y$-axis: the largest relative prediction error for the Burgers equation  with smooth initial profile and shock formation  in Section \ref{sec:burgers_prediction_smooth}.\label{fig:burgers_smooth_prediction_error_vs_order}}
\end{figure}

In Figure \ref{fig:burgers_smooth_prediction_error_vs_order}, we present the largest relative prediction error from $t=t_0$ to $t=1.6$ when the ROM is trained with the data \mw{sampled from} $t=[0,t_0]$ \mw{for} $t_0=1.1$ or \mw{$t_0 = 0.6$} versus the reduced order $r$. For both cases, the LP method always has the smallest relative error. When trained with $t_0=1.1$,  the relative error of the LP method is approximately $0.71\%$ with $20$ basis functions and $0.15\%$ with 25 basis functions. For the same number of reduced basis, the relative error of the LP method trained with $t_0=1.1$ is smaller than that with $t_0=0.6$.

%%%%%%%%%%%%%%%%%%%%%%%%%%%%%%%%%%%%%%%%%%%%%%%%%%
%%%%%%%%%%%%%%%%%%%%%%%%%%%%%%%%%%%%%%%%%%%%%%%%%%%%%%

%%%%%%%%%%%%%%%%%%%%%%%%%

\begin{figure}[h!]
\centering
%%%%%%%%%%%%%%%%%%%%%%%%%
\begin{subfigure}{.48\textwidth}
  \includegraphics[width=\textwidth,trim={5cm 5.34cm 3.25cm 3.5cm},clip]{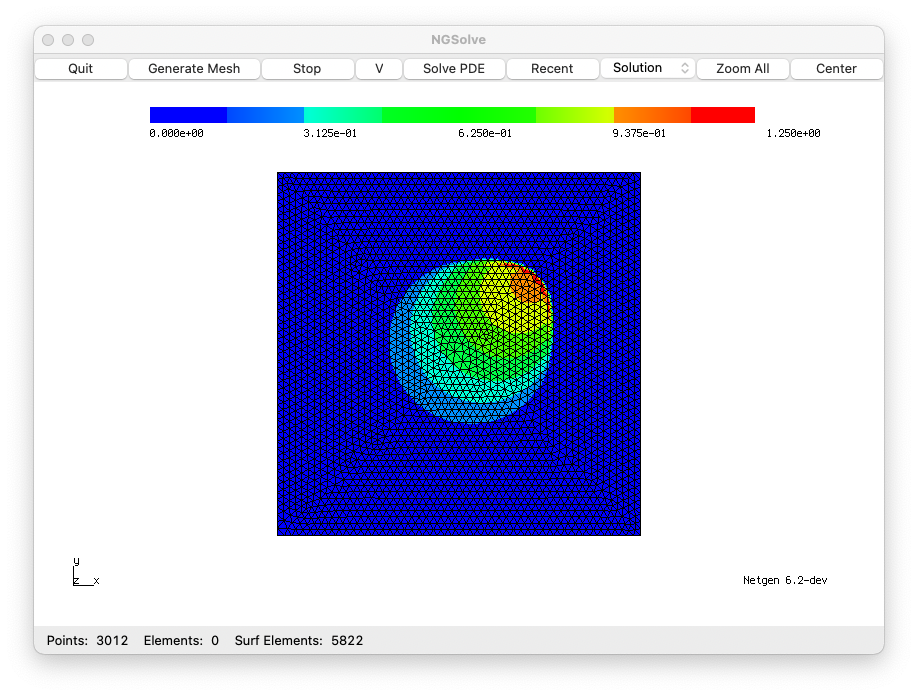}
  \caption{Full order solution without artificial viscosity\label{fig:2d_burgers_no_viscosity}}
\end{subfigure}  
\begin{subfigure}{.48\textwidth}
  \centering
  \includegraphics[width=\textwidth,trim={5cm 5.34cm 3.25cm 3.5cm},clip]{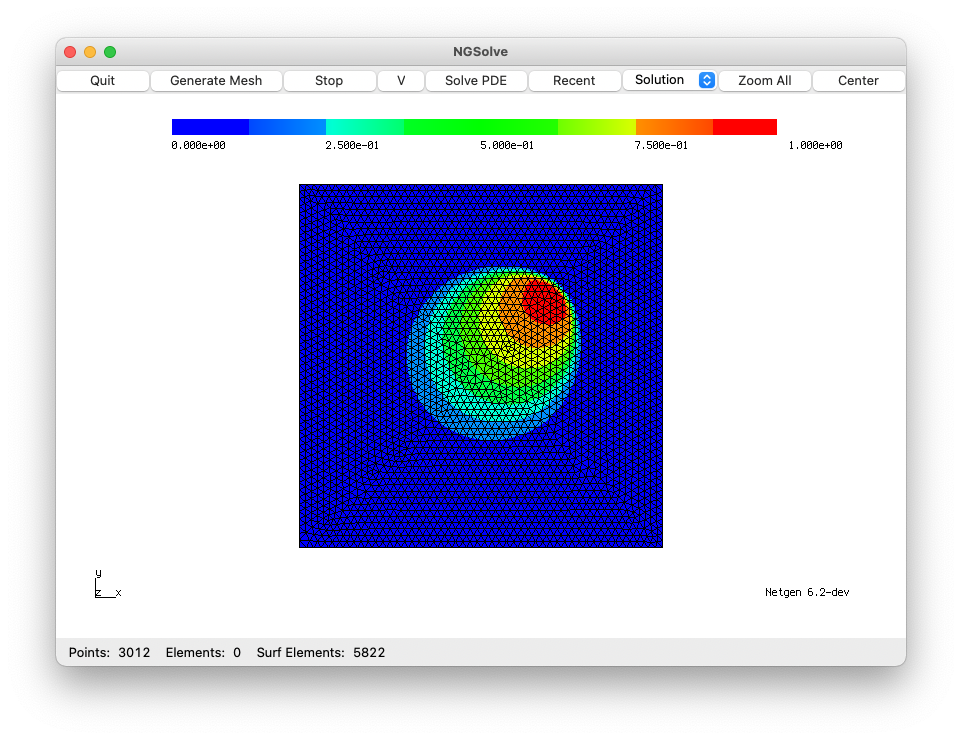}
    \caption{\pzcreview{Full order solution with artificial viscosity
    %.
    }\label{fig:2d_burgers_full}}
\end{subfigure}
\begin{subfigure}{.48\textwidth}
  \centering
  \includegraphics[width=\textwidth,trim={5cm 5.34cm 3.25cm 3.5cm},clip]{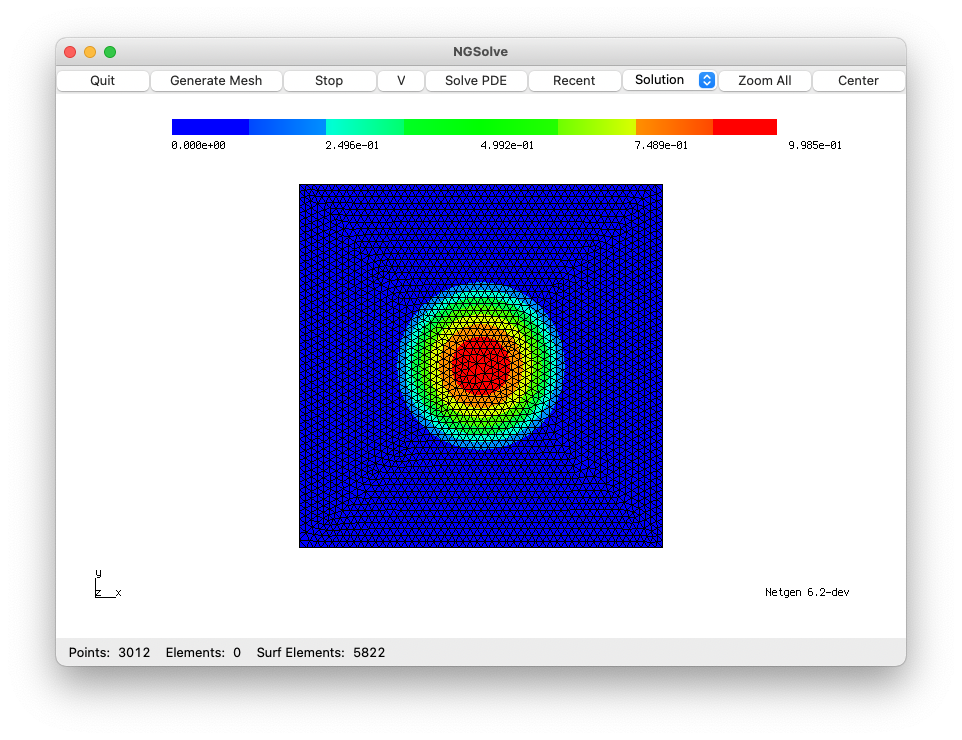}
  \caption{\pzcreview{Initial condition
  %.
  }\label{fig:2d_burgers_initial} }
\end{subfigure}
\caption{\pzcreview{Initial condition, and full order solutions with and without artificial  viscosity at $t=0.42$ for the 2D Burgers equation in Section \ref{sec:2d_burgers}.}\label{fig:2d_burgers_sol-1}}
%%%%%%%%%%%%%%%%%%%%%%%%%
\end{figure}

%%%%%%%%%%%%%%%%%%%%%%%%%%%%%%%%%%%%%%%%%%%%%%%%%%%%%%%%%%%%%%%%%%%%%%%%%%%%%%%%%
% 2D Burgers equation
%%%%%%%%%%%%%%%%%%%%%%%%%%%%%%%%%%%%%%%%%%%%%%%%%%%%%%%%%%%%%%%%%%%%%%%%%%%%%%%%%
\pzcreview{\subsection{2D  Burgers equation\label{sec:2d_burgers}}}
\pzcreview{We consider a 2D  Burgers equation:
\begin{align}
u_t + \nabla\cdot \mathbf{F}(u) = 0,\quad (x,y)\in[-1,1]\times [-1, 1],
\end{align}
with $\mathbf{F}(u)=(\frac{u^2}{2},\frac{u^2}{2})^T$ and the following 
initial condition (also see \ref{fig:2d_burgers_initial}),
\begin{align}
    u(x,0) = \begin{cases}
             \sin\left(\pi(0.5-\sqrt{x^2+y^2})\right),\quad \sqrt{x^2+y^2}\leq0.5,\\
             0,\quad\text{otherwise}.
             \end{cases}
\end{align}
The computational domain is partitioned by a triangular mesh $\mathcal{T}_h=\{\mathcal{T}_i\}_{i=1}^{N_x}$, with $N_x=5822$ and the maximal edge length 
$h_{\textrm{max}}=0.04$.
In space, we apply a DG method: 
seeking $u_h\in V_h := \bigoplus_{\mathcal{T}_i\in\mathcal{T}_h} \mathcal{P}^1(\mathcal{T}_i)$
%,\;k=1$ 
satisfying
\begin{align}
  \sum_{\mathcal{T}_i\in \mathcal{T}_h} \left(\int_{\mathcal{T}_i} (\partial_t u_h v - \mathbf{F}(u_h)\cdot\nabla v)dx + \int_{\partial \mathcal{T}_i} \widehat{\mathbf{F}\cdot\mathbf{n}_i} vds \right)= 0, \quad \forall v \in V_h.
\end{align}
Here we use the local Lax-Fridriches flux
\review{
\begin{align}
\widehat{\mathbf{F}\cdot\mathbf{n}_i}
%   \hat{\mathbf{F}}\cdot\mathbf{n}|_{\partial T_i}
= \frac{(\mathbf{F}(u_h^+)+\mathbf{F}({u_h}^-)}{2} \cdot \mathbf{n}_i + \frac{\lambda_i}{2} (u_h^+ - u_h^-),
\end{align}
 where $\mathbf{n}_i$ is the outward unit normal of $\mathcal{T}_i$ along its boundary $\partial\mathcal{T}_i$, 
 }
$u_h^+$ is the value of $u_h$ on $\partial{\mathcal{T}_i}$ from the interior of $\mathcal{T}_i$, $u_h^-$ is the value of $u_h$ on $\partial{\mathcal{T}_i}$ from the neighboring element of $\mathcal{T}_i$ sharing the same edge, and
$\lambda_i=
|[1,1]^T\cdot \mathbf{n}_i| \max (|u_h^+|, |u_h^-|)$.
In time, we apply the forward Euler method with the time step size $\Dt=0.001$. To stabilize 
%both of the full order solver and the reduced order solver, 
both the full and reduced order solvers,  an artificial viscosity term $\frac{1}{8} h\Delta u$ is added to the equation in an element-wise fashion (with $h|_{\mathcal{T}_i}$ as the element diameter), and is discretized by the symmetric interior penalty method \cite{wheeler1978elliptic}.
The training data is the solution snapshots from $t=0.2$ to $t=0.35$ \pzcnew{($151$ time snapshots)}, and we predict the solution at $t=0.42$. The coefficient neural network has $2$ hidden layers with $25$ neurons per layer, and the basis neural network has $3$ layers with $25$ neurons per layer. In Figure \ref{fig:2d_burgers_sol-1}, we present the full order solution with and without artificial viscosity. Without artificial viscosity, one can see the non-physical oscillations  near the wave front in the full order solution, and the LP method produces even more oscillatory solutions. 
}

\begin{figure}[h!]
\centering
%%%%%%%%%%%%%%%%%%%%%%%%%
%%%%%%%%%%%%%%%%%%%%%%%%%
\begin{subfigure}{.48\textwidth}
  \centering
  \includegraphics[width=\textwidth,trim={6cm 5.34cm 3.8cm 3.5cm},clip]{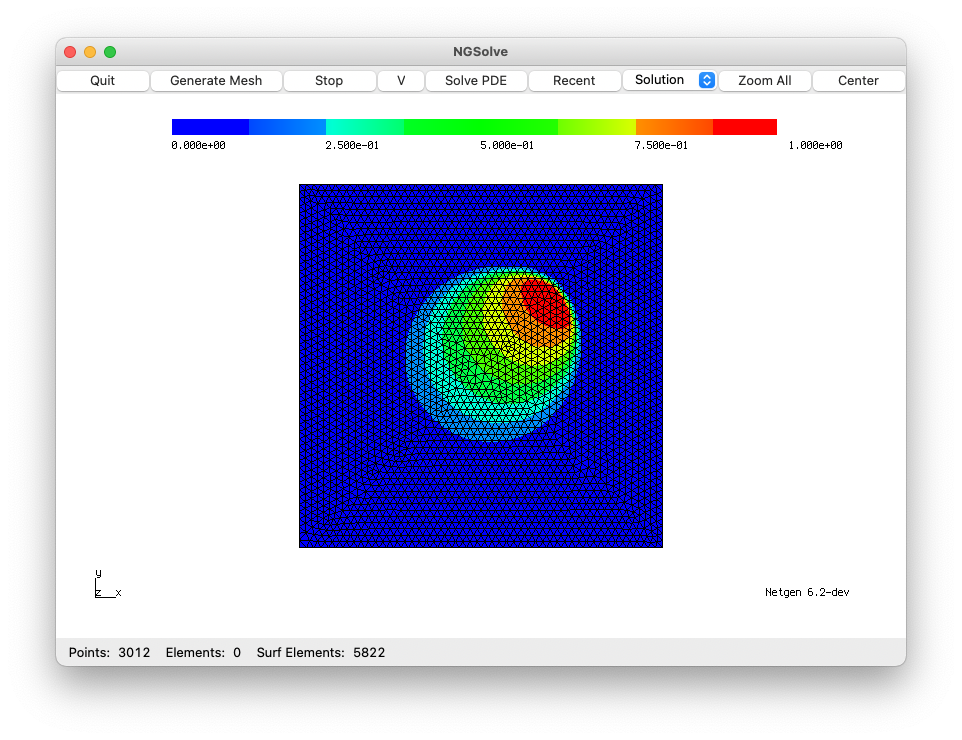}
    \caption{\pzcreview{LP solution with 10  basis functions
    %.
    }\label{fig:2d_burgers_full}}
\end{subfigure}
%%%%%%%%%%%%%%%%%%%%%%%%%
\begin{subfigure}{.48\textwidth}
  \centering
  \includegraphics[width=\textwidth,trim={6cm 5.34cm 3.8cm 3.5cm},clip]{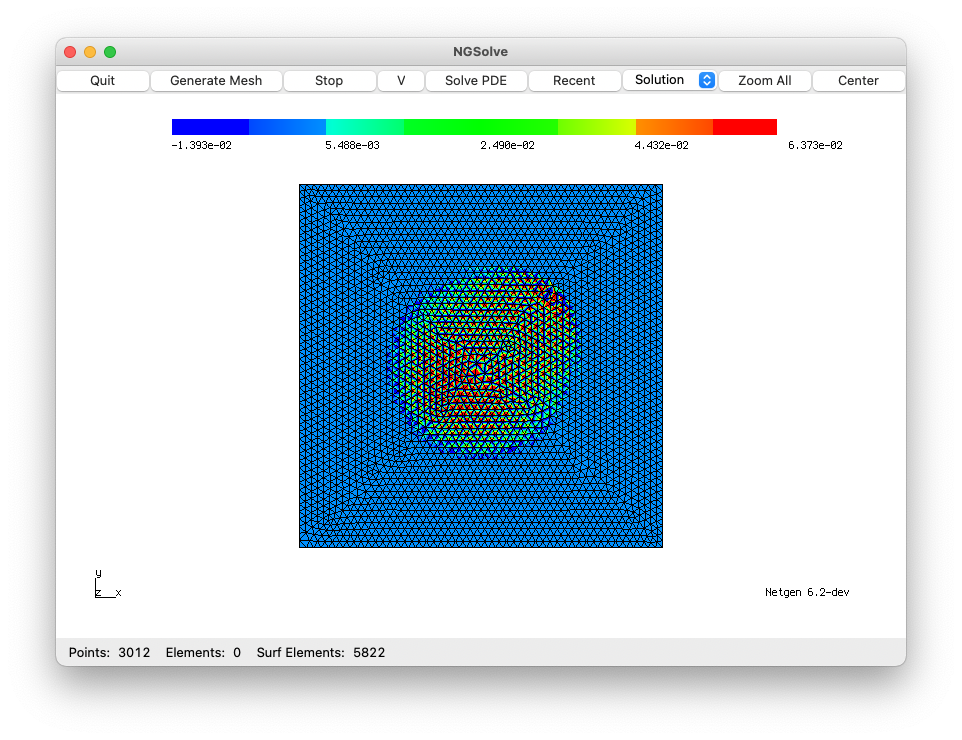}
    \caption{\pzcreview{POD solution with $10$  basis functions
    %.
    }}
\end{subfigure}
%%%%%%%%%%%%%%%%%%%%%%%%%
\begin{subfigure}{.48\textwidth}
  \centering
  \includegraphics[width=\textwidth,trim={6cm 5.34cm 3.8cm 3.5cm},clip]{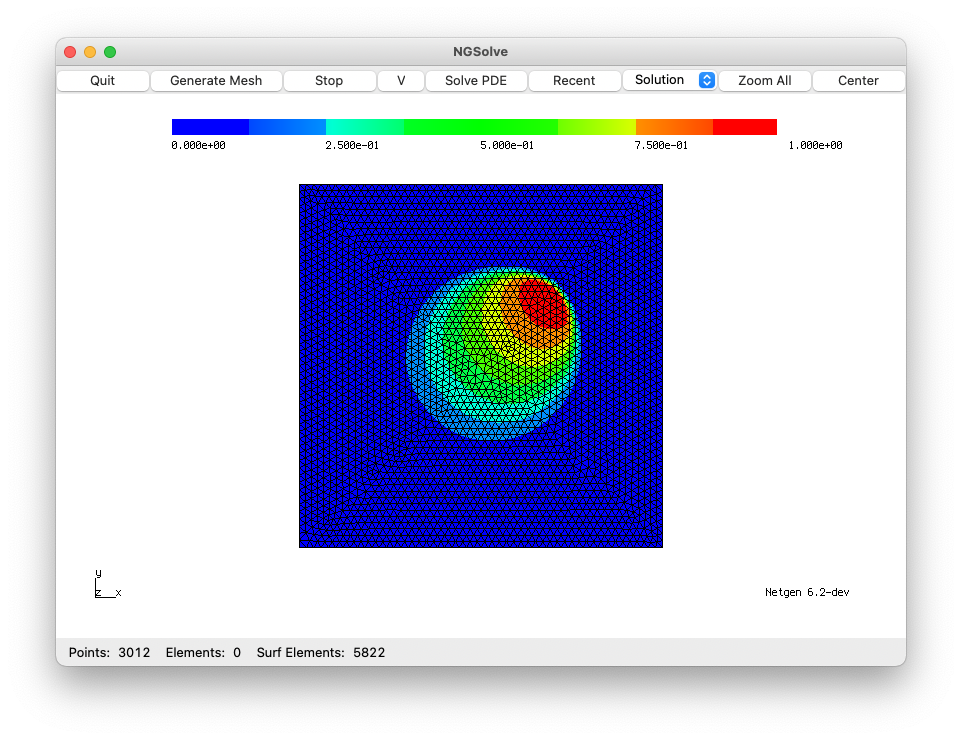}
    \caption{\pzcreview{
    %Pure 
    Learning solution with 10 basis functions
    %.
    }}
\end{subfigure}
%%%%%%%%%%%%%%%%%%%%%%%%%
\begin{subfigure}{.48\textwidth}
  \centering
  \includegraphics[width=\textwidth,trim={6cm 5.34cm 3.8cm 3.5cm},clip]{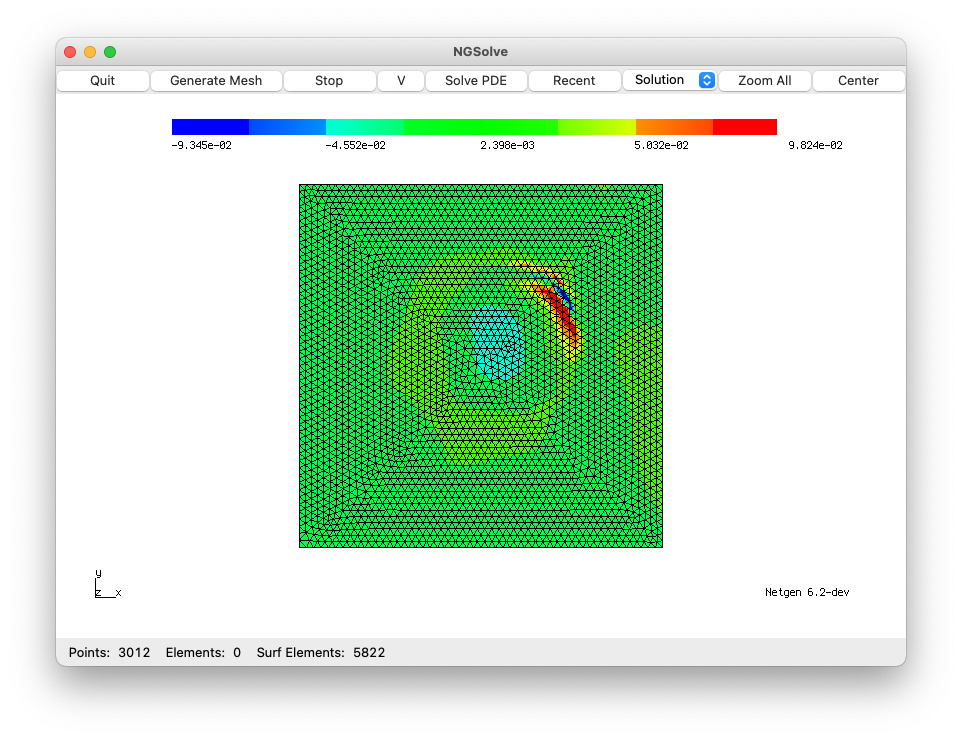}
    \caption{\pzcreview{Error of the LP method with $10$  basis functions
    %.
    }}
\end{subfigure}
%%%%%%%%%%%%%%%%%%%%%%%%%
%%%%%%%%%%%%%%%%%%%%%%%%%
\caption{\pzcreview{Full order and reduced order solutions and error at $t=0.42$ for the 2D Burgers equation  in Section \ref{sec:2d_burgers}.}\label{fig:2d_burgers_sol}}
%%%%%%%%%%%%%%%%%%%%%%%%%
\end{figure}

\begin{figure}[h!]
\centering
%%%%%%%%%%%%%%%%%%%%%%%%%
\begin{subfigure}{.45\textwidth}
  \centering
  \includegraphics[width=\textwidth]{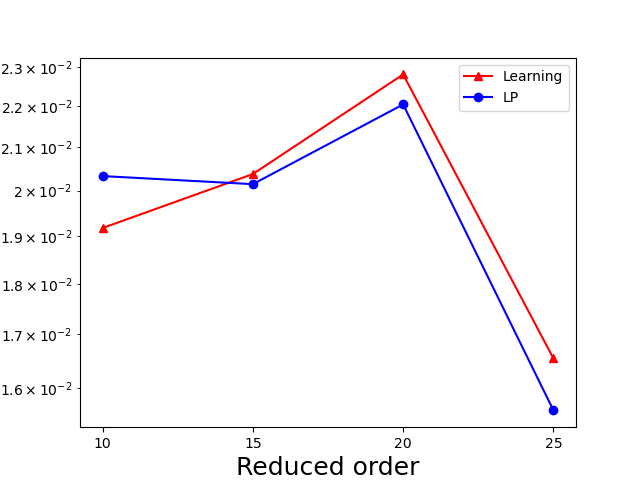}
  \caption{\pzcreview{Absolute error at $t=0.42$
  %.
  }\label{fig:2d_advection_error} }
\end{subfigure}
%%%%%%%%%%%%%%%%%%%%%%%%%
\caption{\pzcreview{The absolute error at $t=0.42$ for the 2D Bureges equation in Section \ref{sec:2d_burgers}.\label{fig:2d_burgers_error_time}}}
\end{figure}

\reviewtwo{As mentioned before, there will be no saving of computational cost with the forward Euler time integrator. The purpose of this test is to check whether the neural network is capable of capturing the underlying low-rank structure in 2D nonlinear hyperbolic problems.}
\review{
The reduced order solutions with 10 basis functions at $t=0.42$ are presented in Figure \ref{fig:2d_burgers_sol},  all with the artificial viscosity strategy applied. The LP method and the pure learning-based method match
%matches 
the full order solution (in Figure \ref{fig:2d_burgers_full}) reasonably well, while the POD method with the same number of basis functions fails to capture the solution structure when going beyond the training set. In  Figure \ref{fig:2d_burgers_error_time}, the absolute errors at $t=0.42$ versus the reduced order $r$  are presented for the LP and the pure learning methods. The errors of the two method are close to each other, but the online stage of the pure learning method is faster, as it only needs function evaluation.}

%%%%%%%%%%%%%%%%%%%%%%%%%

%%%%%%%%%%%%%%%%%%%%%%%%%

%%%%%%%%%%%%%%%%%%%%%%%%%%%%%%%%%%%%%%%%%%%%%%%%%%

%%%%%%%%%%%%%%%%%%%%%%%%%%%%%%%%%%%%%%%%%%%%%%%%%%%%%%%%%%%%%%%%%%%%%%%%
\subsection{
Compressible Euler system
 \label{sec:euler}}
 Finally we apply the proposed LP method to the  {1D} compressible Euler system of gas dynamics:
\begin{align}
\begin{cases}
 &\rho_t + (\rho u)_x=0,\\
 &(\rho u)_t + (\rho u^2+p)_x=0,\\
 &E_t + \left(u(E+p)\right)_x=0.
 \end{cases}
\end{align}
Here $\rho, u, p$ are the density, velocity, pressure, respectively, with the total energy 
$E=\frac{p}{\gamma-1}+\frac{1}{2}\rho u^2$ and $\gamma=3$ for the ideal gas.
We consider the Sod's shock tube problem, with the initial conditions:
\begin{align}
\rho(x,0)=\begin{cases}
          1.0,\quad x<0,\\
          0.125,\quad x\geq0,
          \end{cases}
\quad
u(x,0)=0.0,
\quad
p(x,0)=\begin{cases}
          1.0,\quad x<0,\\
          0.1,\quad x\geq0.
          \end{cases}
\end{align}
For this example, the full order method is based on the fifth-order WENO finite difference method with the Lax-Friedrichs flux \cite{jiang1996efficient} and the explicit third order strong stability preserving Runge-Kutta method \cite{shu1988efficient}. \reviewtwo{With an  explicit time integrator applied, the  LP method will not reduce the computational time. We want to use this example to demonstrate that the proposed method can capture the hidden low-rank structure of nonlinear hyperbolic systems. }

We consider the computational domain $[-1.5,1.5]$, and it is partitioned by a uniform mesh with $N_x=1500$ elements. \review{The CFL number is taken as $0.6$.} The training data are the solutions for $t\in[0,0.2]$, and we use ROMs to predict solutions for $t\in[0.2,0.25]$. In the offline stage of the LP method, we train three neural networks $\rho_{\NN}$, $(\rho u)_{\NN}$ and $E_{\NN}$, with each using the same architecture as we specify previously. 
The basis blocks of these neural networks further provide $t$-dependent reduced basis for $\rho$, $\rho u$ and $E$ respectively. 
$100$ epochs \mw{of training} are \mw{conducted} in the offline \mw{stage}.  

\begin{figure}[h!]
\centering
%%%%%%%%%%%%%%%%%%%%%%%%%
\begin{subfigure}{.32\textwidth}
  \centering
  \includegraphics[width=\textwidth]{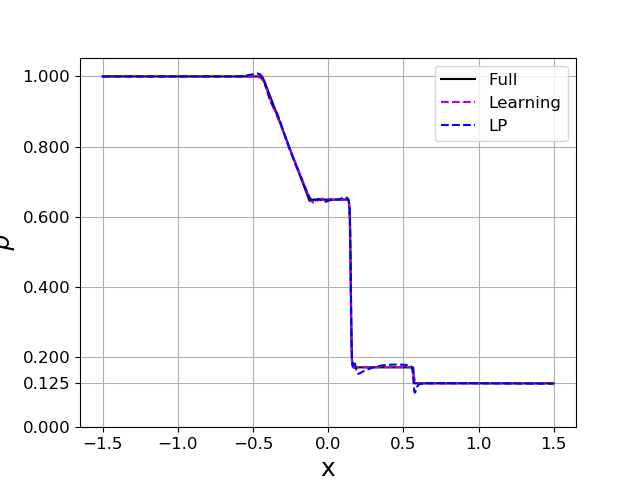}
  \caption{20 basis functions, density}
\end{subfigure}
%%%%%%%%%%%%%%%%%%%%%%%%%
\begin{subfigure}{.32\textwidth}
  \centering
  \includegraphics[width=\textwidth]{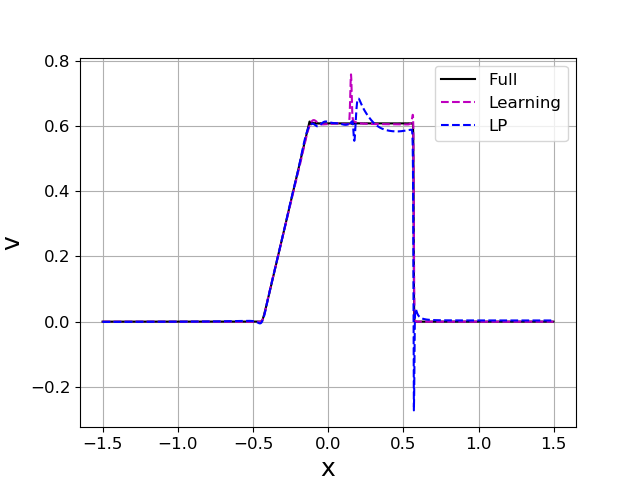}
  \caption{20 basis functions, velocity}
\end{subfigure}
%%%%%%%%%%%%%%%%%%%%%%%%%
\begin{subfigure}{.32\textwidth}
  \centering
  \includegraphics[width=\textwidth]{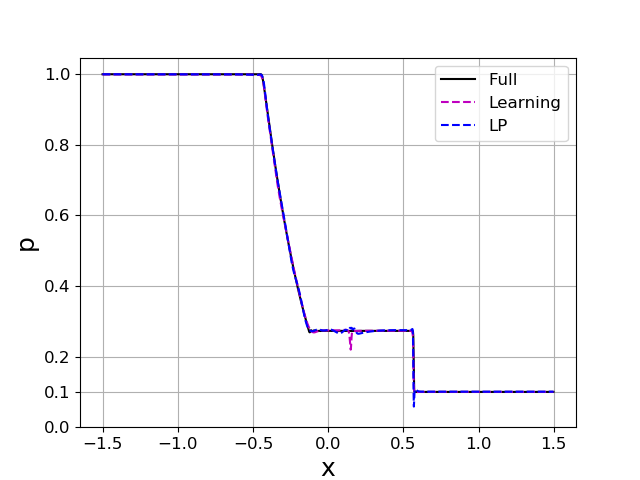}
  \caption{20 basis functions, pressure}
\end{subfigure}
%%%%%%%%%%%%%%%%%%%%%%%%%
\begin{subfigure}{.32\textwidth}
  \centering
  \includegraphics[width=\textwidth]{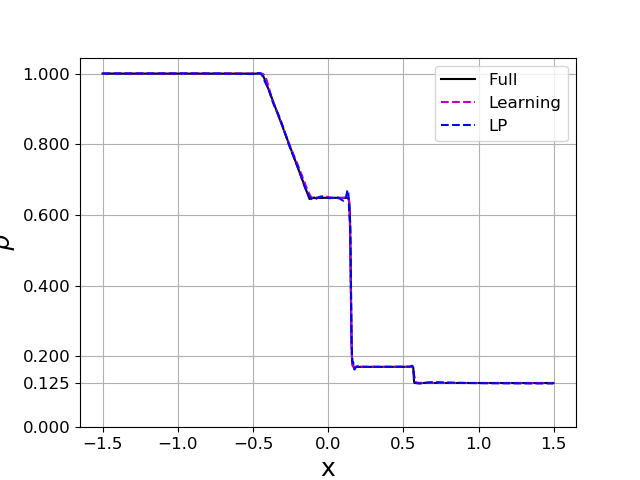}
  \caption{30 basis functions, density}
\end{subfigure}
%%%%%%%%%%%%%%%%%%%%%%%%%
\begin{subfigure}{.32\textwidth}
  \centering
  \includegraphics[width=\textwidth]{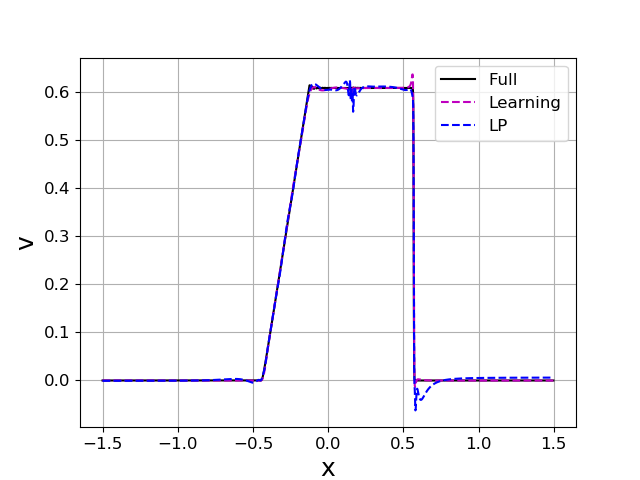}
  \caption{30 basis functions, velocity}
\end{subfigure}
%%%%%%%%%%%%%%%%%%%%%%%%%
\begin{subfigure}{.32\textwidth}
  \centering
  \includegraphics[width=\textwidth]{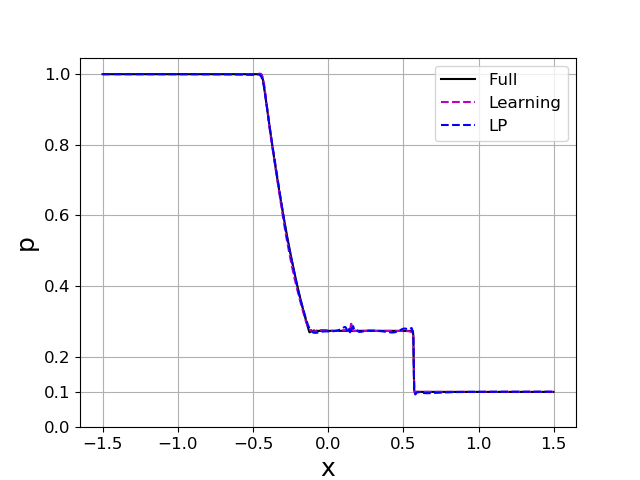}
  \caption{30 basis               functions, pressure}
\end{subfigure}
%%%%%%%%%%%%%%%%%%%%%%%%%
\caption{Predicted solutions  for the Euler equation in Section \ref{sec:euler}
at $t=0.25$.\label{fig:euler_solution}}
\end{figure}

In Figure \ref{fig:euler_solution}, we present the density, velocity, and pressure of the  reduced order solutions  at $t=0.25$ by the LP and the pure learning methods with $r=20, 30$. The solution includes a rarefaction wave,  a contact discontinuity, and a shock.  Overall, the predicted solutions by both methods capture {these} wave structures well.  Numerical oscillations are observed near the rightward moving contact discontinuity and the shock, and they are further suppressed  when the reduced order $r$ increases.

\begin{figure}[h!]
\centering
%%%%%%%%%%%%%%%%%%%%%%%%%
\begin{subfigure}{.32\textwidth}
  \centering
  \includegraphics[width=\textwidth]{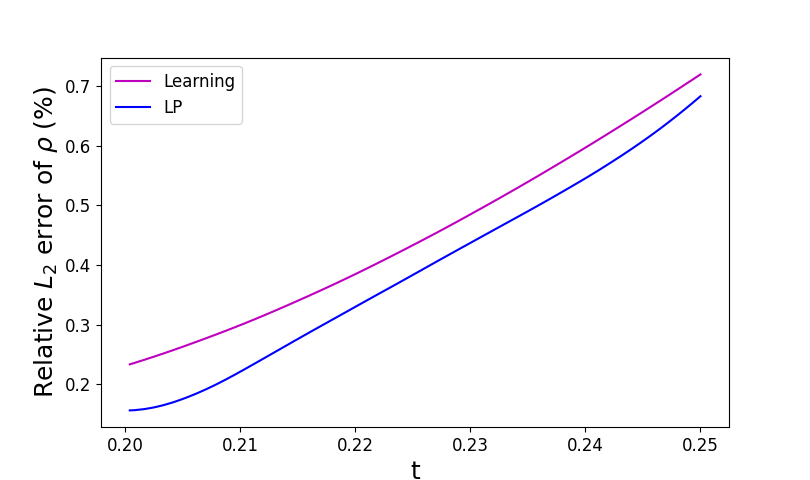}
  \caption{20 basis functions, density}
\end{subfigure}
%%%%%%%%%%%%%%%%%%%%%%%%%
\begin{subfigure}{.32\textwidth}
  \centering
  \includegraphics[width=\textwidth]{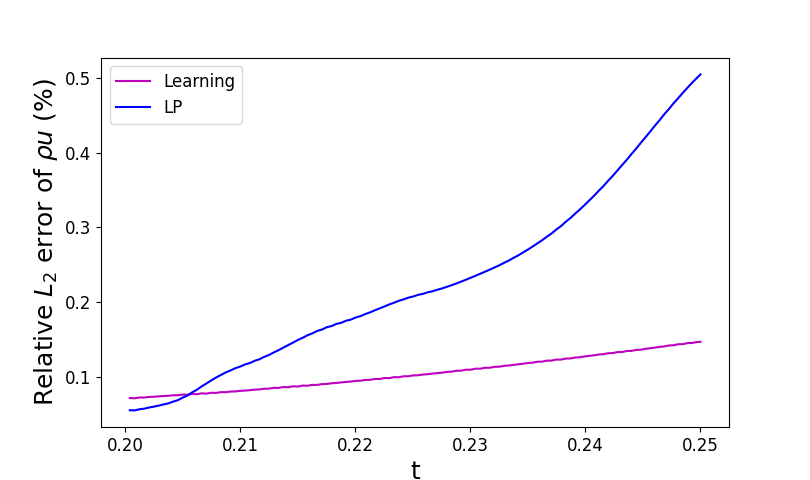}
  \caption{20 basis functions, velocity}
\end{subfigure}
%%%%%%%%%%%%%%%%%%%%%%%%%
\begin{subfigure}{.32\textwidth}
  \centering
  \includegraphics[width=\textwidth]{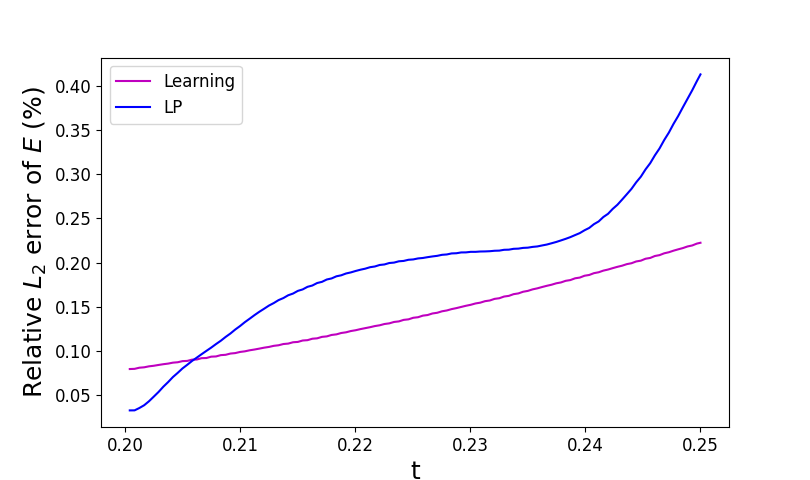}
  \caption{20 basis functions, pressure}
\end{subfigure}
%%%%%%%%%%%%%%%%%%%%%%%%%
\begin{subfigure}{.32\textwidth}
  \centering
  \includegraphics[width=\textwidth]{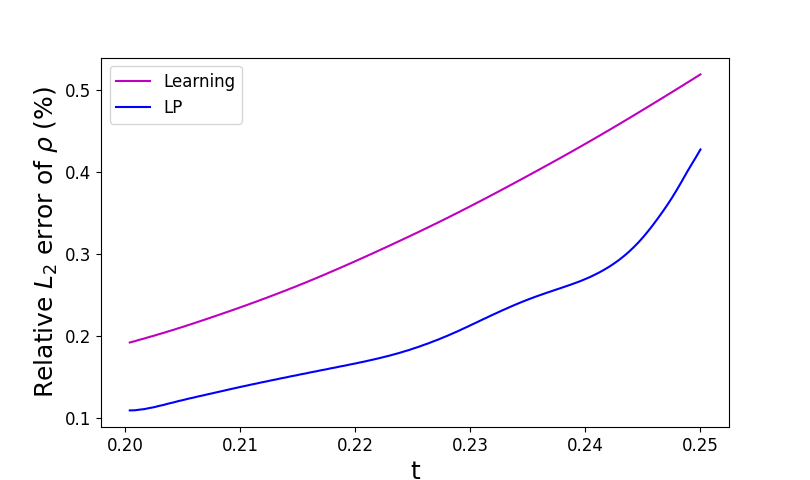}
  \caption{30 basis, density}
\end{subfigure}
%%%%%%%%%%%%%%%%%%%%%%%%%
\begin{subfigure}{.32\textwidth}
  \centering
  \includegraphics[width=\textwidth]{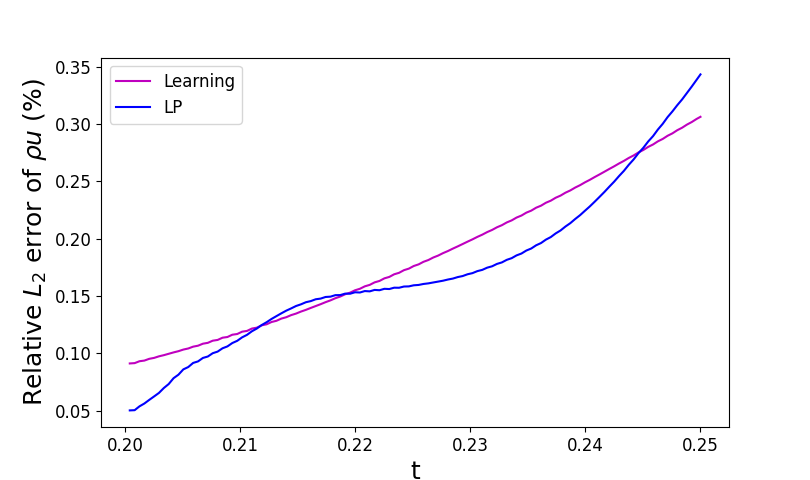}
  \caption{30 basis functions, velocity}
\end{subfigure}
%%%%%%%%%%%%%%%%%%%%%%%%%
\begin{subfigure}{.32\textwidth}
  \centering
  \includegraphics[width=\textwidth]{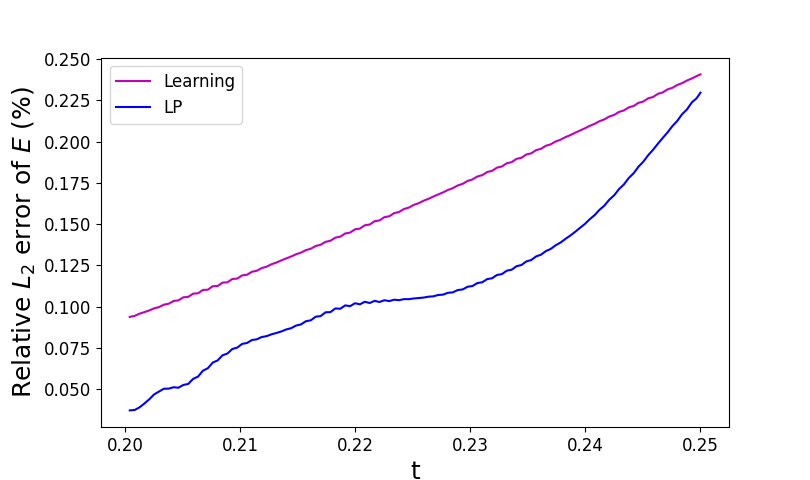}
  \caption{30 basis functions, pressure}
\end{subfigure}
%%%%%%%%%%%%%%%%%%%%%%%%%
\caption{The evolution of the prediction relative error of the Euler equation in Section \ref{sec:euler}
from $t=0.2$ to $t=0.25$.\label{fig:euler_error}}
\end{figure}

In Figure \ref{fig:euler_error}, we plot the time evolution of the 
relative error  by the LP and the pure learning methods in prediction from  $t=0.2$ to $t=0.25$, and the errors are  reasonably small. Compared with examples in previous sections, the error of the LP method grows faster \mw{in} this example. We believe this is due to the insufficient control of numerical oscillations in \mw{the} online projection step \eqref{eq:online_projection}.
We also want to mention that with $10$ basis functions, the LP method suffers from negative pressure after $t=0.2315$ due to the numerical oscillation. 

Note that the results \mw{of} the POD method are not reported \mw{in this example}. \mw{In fact, when} $r=20, 30$, the online computation with the POD basis fails due to the presence of non-physical negative pressure. Indeed, before this \mw{phenomenon} occurs, the POD results are more oscillatory than \mw{that of} the other two methods. It is an important task to explore how to effectively control numerical oscillation and to preserve positivity of certain physical quantities for ROMs.

%%%%%%%%%%%%%%%%%%%%%%%%%%%%%%%%%%%%%%%%%%%%%%%%%%%%%%%%%
%%%%%%%%%%%%%%%%%%%%%%%%%%%%%%%%%%%%%%%%%%%%%%%%%%%%%%%%%

%%%%%%%%%%%%%%%%%%%%%%%%%%%%%%%%%%%%%%%%%%%%%%%%%%%%%%%%%
%%%%%%%%%%%%%%%%%%%%%%%%%%%%%%%%%%%%%%%%%%%%%%%%%%%%%%%%%

\section{Conclusions\label{sec:conclusions}}
Following the offline-online decomposition, we propose a learning-based projection method for transport problems. 

As a combination of deep learning and the  projection-based method, the proposed method has the following advantages:
\begin{enumerate}
\item Compared with the standard POD methods, it achieves smaller long-time prediction error with the same number of basis functions.
\pzcreview{
\item Compared with a pure learning-based method using the same neural network, the online projection step may reduce the generalization error dramatically.
\item
Unlike previous learning-based projection methods \cite{lee2019deep,lee2020model,kim2020fast}, there is no need to take derivatives of neural networks in the online projection step. 
\item With $x$ as an input of the neural network in a mesh free manner, it is flexible and simple for the proposed method
%s
to work with different meshes and even different schemes in the online and offline stages.
}
\end{enumerate}
\pzcreview{However, for problems with strong shocks, the projection based online stage may produce more numerical oscillation compared with a pure learning-based online stage.}

\pzcreview{There are still many aspects  \mw{in this topic that could be further explored.} 
\begin{enumerate}
\item To design
%Design 
an efficient and robust hyper-reduction strategy to save the computational cost with explicit time integrators.
\item The effective and efficient control of numerical oscillations is
%are 
crucial, especially when combined with hyper-reduction. 
\item To apply and test the proposed method with
%to 
more complicated 2D and 3D system hyperbolic conservation laws. 
\item To further combine the proposed method with the adaptive moving mesh  and the $h,p$-adaptive numerical methods which move
%moves 
or regenerate computational meshes for different time and/or  parameter samples. 
\item  One drawback of the proposed method is that the learned basis is  not orthonormal. To improve the robustness associated with this, 
one may introduce a regularization term similar to \cite{wang2020reduced} in the loss function.
\end{enumerate}
}

\begin{appendices}
\section{\pzcrev{More on the oscillation in Figure \ref{fig:advection_hetero_solution} with $r=15$}}
\label{sec:appendix}

\begin{figure}
       \centering
       \begin{subfigure}{.38\textwidth}
        \includegraphics[width=\textwidth]{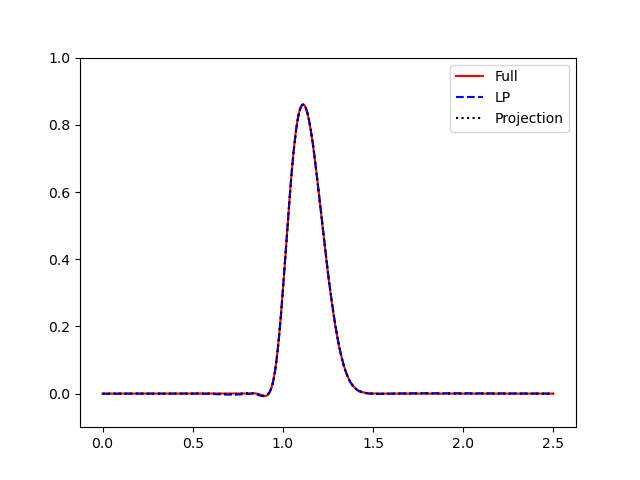}
        \caption{$t=0.6$        \label{fig:app.a}}
       \end{subfigure}
       \begin{subfigure}{.38\textwidth}
        \includegraphics[width=\textwidth]{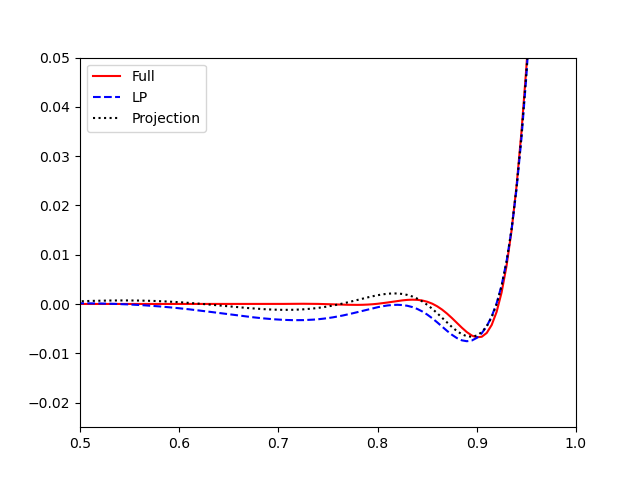}
        \caption{$t=0.6$, zoomed-in}
       \end{subfigure}
        \begin{subfigure}{.38\textwidth}
        \includegraphics[width=\textwidth]{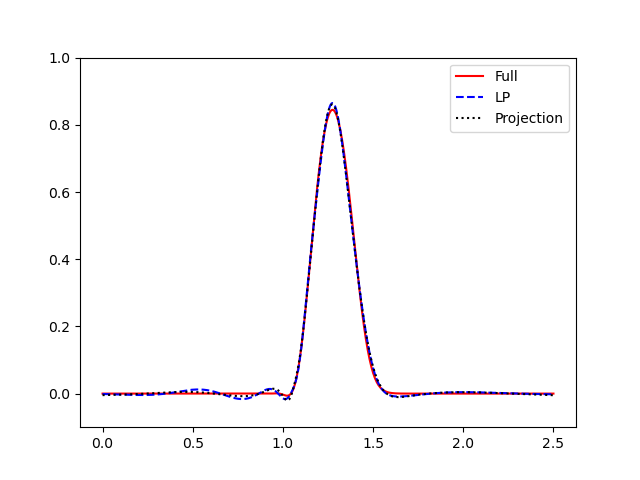}
        \caption{$t=0.7$}
       \end{subfigure}
       \begin{subfigure}{.38\textwidth}
        \includegraphics[width=\textwidth]{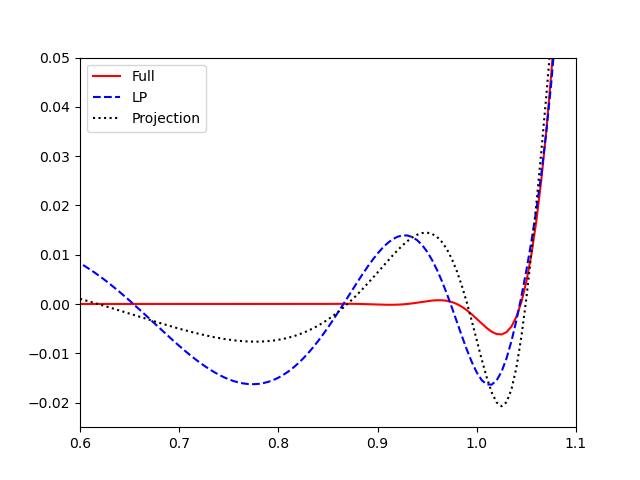}
        \caption{$t=0.7$, zoomed-in}
       \end{subfigure}
        \begin{subfigure}{.38\textwidth}
        \includegraphics[width=\textwidth]{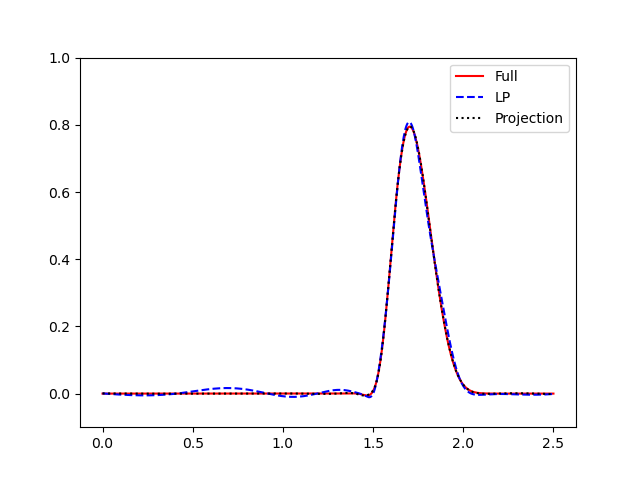}
        \caption{$t=1.0$}
       \end{subfigure}
       \begin{subfigure}{.38\textwidth}
        \includegraphics[width=\textwidth]{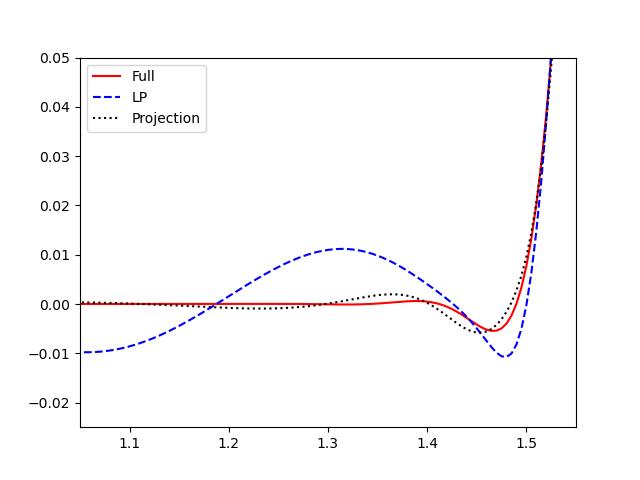}
        \caption{$t=1.0$, zoomed-in
                \label{fig:app.f}}
       \end{subfigure}
      \begin{subfigure}{.38\textwidth}
        \includegraphics[width=\textwidth]{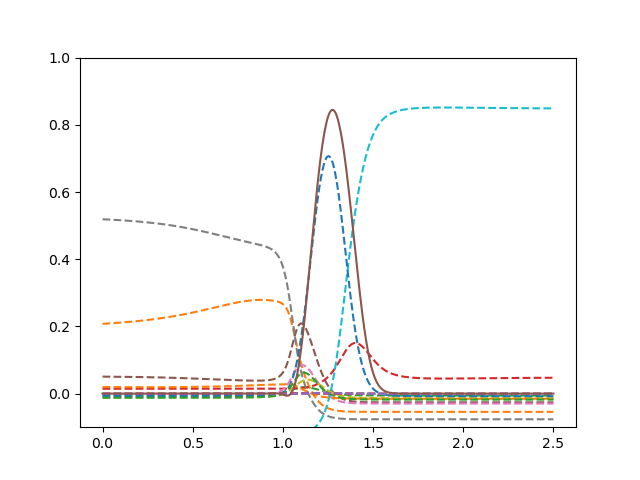}
        \caption{$15$ basis functions (dashed lines) and the full order solution (solid line) at $t=0.7$  \label{fig:sol-basis-overlay0.6}}
       \end{subfigure}
       \begin{subfigure}{.38\textwidth}
        \includegraphics[width=\textwidth]{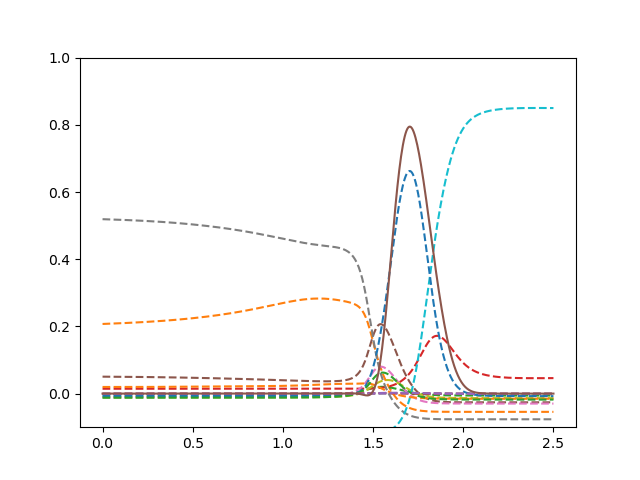}
        \caption{$15$ basis functions (dashed lines) and the full order solution (solid line) at $t=1.0$ \label{fig:sol-basis-overlay0.7}}
       \end{subfigure}
       \caption{\lli{The inhomogeneous media example in Section \ref{sec:advection_hetero}. (a)-(f):  solutions at different times, where ``Full'' is the full order solution, ``LP'' is by the learning-based projection method, and ``Projection'' is  the $L^2$ projection of the full order solution onto the learned reduced subspace. (g)-(h): $15$ learned basis functions (dashed lines) and the full order solution (solid line) at $t=0.7$ and $t=1.0$. \label{fig:sol_osc}}}
 \end{figure}

\pzcrev{
In Figure \ref{fig:advection_hetero_solution}, non-physical oscillation can be observed
in the reduced order solution generated by the LP method with $r=15$. To gain some understanding to the origin of such oscillation,  in  Figures         \ref{fig:app.a}-\ref{fig:app.f}, we present the evolution of the full order solution, reduced order solution by the proposed LP method, and the $L^2$ projection of the full order solution onto the learned adaptive reduced space. 
At $t=0.6$, all three curves match well, with slight difference observed in the zoomed-in plot.   
 At $t=0.7$, the oscillation in both the LP solution and the projection becomes quite visible. At $t=1.0$, the oscillation in the LP solution remains and spreads out, while the projection of the full order solution displays little difference from the full order solution.
In Figures \ref{fig:sol-basis-overlay0.6}-\ref{fig:sol-basis-overlay0.7}, we further overlay the 15 learned basis functions and the full order solution at $t=0.7$ and  $t=1.0$, respectively. One can see that the learned basis functions faithfully capture the local nature of the full order solution, without visible oscillation in themselves. The basis functions at these two different times also capture the propagation of the solution as well as its slight reduction in amplitude as time increases.  From the observations above, we attribute the numerical oscillation in the LP solution with $r=15$ to the accumulation of the errors over time,  the  extrapolation of the basis functions from the training set to a test sample, and also to the online projection step. 
 The underlying media is inhomogeneous, and this also contributes to the difference observed at $t=0.7$ and $t=1.0$ in the shape and amplitude of the basis functions and hence of the oscillation in the LP solution. Adaptively adding artificial viscosity in the online phase of the LP method to control these oscillations is worthy future research.}

\end{appendices}

\bigskip\noindent
{\bf Acknowledgement.}
The authors want to thank Dr. Juntao Huang from Michigan State University for generously sharing his WENO codes solving the compressible Euler system.

%%%%%%%%%%%%%%%%%%%%%%%%%%%%%%%%%%%%%%%%%%%%%%%%%%%%%%%%
% Reference
%%%%%%%%%%%%%%%%%%%%%%%%%%%%%%%%%%%%%%%%%%%%%%%%%%%%%%%%
\newpage
\renewcommand\refname{References}
\bibliographystyle{plain}
\bibliography{bib_peng}

\begin{thebibliography}{10}

\bibitem{abadi2016tensorflow}
Mart{\'\i}n Abadi, Paul Barham, Jianmin Chen, Zhifeng Chen, Andy Davis, Jeffrey
  Dean, Matthieu Devin, Sanjay Ghemawat, Geoffrey Irving, Michael Isard, et~al.
\newblock {Tensorflow: A system for large-scale machine learning}.
\newblock In {\em 12th {USENIX} symposium on operating systems design and
  implementation ({OSDI} 16)}, pages 265--283, 2016.

\bibitem{battisti2022wasserstein}
Beatrice Battisti, Tobias Blickhan, Guillaume Enchery, Virginie Ehrlacher,
  Damiano Lombardi, and Olga Mula.
\newblock {Wasserstein model reduction approach for parametrized flow problems
  in porous media}.
\newblock 2022.

\bibitem{benner2017model}
Peter Benner, Mario Ohlberger, Albert Cohen, and Karen Willcox.
\newblock {\em {Model reduction and approximation: theory and algorithms}}.
\newblock SIAM, 2017.

\bibitem{berkooz1993proper}
Gal Berkooz, Philip Holmes, and John~L Lumley.
\newblock {The proper orthogonal decomposition in the analysis of turbulent
  flows}.
\newblock {\em Annual review of fluid mechanics}, 25(1):539--575, 1993.

\bibitem{beyn2004freezing}
W-J. Beyn and V.~Th{\"u}mmler.
\newblock {Freezing solutions of equivariant evolution equations}.
\newblock {\em SIAM Journal on Applied Dynamical Systems}, 3(2):85--116, 2004.

\bibitem{buffa2012priori}
Annalisa Buffa, Yvon Maday, Anthony~T Patera, Christophe Prud’homme, and
  Gabriel Turinici.
\newblock {A priori convergence of the greedy algorithm for the parametrized
  reduced basis method}.
\newblock {\em ESAIM: Mathematical Modelling and Numerical
  Analysis-Mod{\'e}lisation Math{\'e}matique et Analyse Num{\'e}rique},
  46(3):595--603, 2012.

\bibitem{cagniart2018few}
N.~Cagniart.
\newblock {\em {A few nonlinear approaches in model order reduction}}.
\newblock PhD thesis, Sorbonne Universit{\'e}, 2018.

\bibitem{cagniart2017model}
N.~Cagniart, R.~Crisovan, Y.~Maday, and R.~Abgrall.
\newblock {Model order reduction for hyperbolic problems: a new framework}.
\newblock 2017.

\bibitem{cagniart2019model}
Nicolas Cagniart, Yvon Maday, and Benjamin Stamm.
\newblock {Model order reduction for problems with large convection effects}.
\newblock In {\em Contributions to partial differential equations and
  applications}, pages 131--150. Springer, 2019.

\bibitem{carlberg2015adaptive}
K.~Carlberg.
\newblock {Adaptive $h$-refinement for reduced-order models}.
\newblock {\em International Journal for Numerical Methods in Engineering},
  102(5):1192--1210, 2015.

\bibitem{chen2021eim}
Yanlai Chen, Sigal Gottlieb, Lijie Ji, and Yvon Maday.
\newblock {An EIM-degradation free reduced basis method via over collocation
  and residual hyper reduction-based error estimation}.
\newblock {\em arXiv preprint arXiv:2101.05902}, 2021.

\bibitem{choi2020sns}
Youngsoo Choi, Deshawn Coombs, and Robert Anderson.
\newblock {SNS: a solution-based nonlinear subspace method for time-dependent
  model order reduction}.
\newblock {\em SIAM Journal on Scientific Computing}, 42(2):A1116--A1146, 2020.

\bibitem{cohen2021optimal}
Albert Cohen, Ronald Devore, Guergana Petrova, and Przemyslaw Wojtaszczyk.
\newblock {Optimal stable nonlinear approximation}.
\newblock {\em Foundations of Computational Mathematics}, pages 1--42, 2021.

\bibitem{dahmen2021nonlinear}
Wolfgang Dahmen, Min Wang, and Zhu Wang.
\newblock {Nonlinear Reduced DNN Models for State Estimation}.
\newblock {\em arXiv preprint arXiv:2110.08951}, 2021.

\bibitem{ehrlacher2019nonlinear}
V.~Ehrlacher, D.~Lombardi, O.~Mula, and F-X. Vialard.
\newblock {Nonlinear model reduction on metric spaces. Application to
  one-dimensional conservative PDEs in Wasserstein spaces}.
\newblock {\em ESAIM: Mathematical Modelling and Numerical Analysis}, 2019.

\bibitem{ferrero2021registration}
Andrea Ferrero, Tommaso Taddei, and Lei Zhang.
\newblock Registration-based model reduction of parameterized two-dimensional
  conservation laws.
\newblock {\em arXiv preprint arXiv:2105.02024}, 2021.

\bibitem{fresca2021comprehensive}
Stefania Fresca, Andrea Manzoni, et~al.
\newblock {A comprehensive deep learning-based approach to reduced order
  modeling of nonlinear time-dependent parametrized PDEs}.
\newblock {\em Journal of Scientific Computing}, 87(2):1--36, 2021.

\bibitem{gerbeau2014approximated}
J-F. Gerbeau and D.~Lombardi.
\newblock {Approximated Lax pairs for the reduced order integration of
  nonlinear evolution equations}.
\newblock {\em Journal of Computational Physics}, 265:246--269, 2014.

\bibitem{gerbeau2015reduced}
J-F. Gerbeau, D.~Lombardi, and E.~Schenone.
\newblock {Reduced order model in cardiac electrophysiology with approximated
  Lax pairs}.
\newblock {\em Advances in Computational Mathematics}, 41(5):1103--1130, 2015.

\bibitem{grassle2018pod}
Carmen Gr{\"a}{\ss}le and Michael Hinze.
\newblock {POD reduced-order modeling for evolution equations utilizing
  arbitrary finite element discretizations}.
\newblock {\em Advances in Computational Mathematics}, 44(6):1941--1978, 2018.

\bibitem{greif2019decay}
C.~Greif and K.~Urban.
\newblock {Decay of the Kolmogorov N-width for wave problems}.
\newblock {\em Applied Mathematics Letters}, 96:216--222, 2019.

\bibitem{gubisch2017proper}
Martin Gubisch and Stefan Volkwein.
\newblock {Proper orthogonal decomposition for linear-quadratic optimal
  control}.
\newblock {\em Model reduction and approximation: theory and algorithms}, 5:66,
  2017.

\bibitem{han2019uniformly}
Jiequn Han, Chao Ma, Zheng Ma, and E~Weinan.
\newblock Uniformly accurate machine learning-based hydrodynamic models for
  kinetic equations.
\newblock {\em Proceedings of the National Academy of Sciences},
  116(44):21983--21991, 2019.

\bibitem{hesthaven2016certified}
Jan~S Hesthaven, Gianluigi Rozza, Benjamin Stamm, et~al.
\newblock {\em {Certified reduced basis methods for parametrized partial
  differential equations}}, volume 590.
\newblock Springer, 2016.

\bibitem{jiang1996efficient}
Guang-Shan Jiang and Chi-Wang Shu.
\newblock Efficient implementation of weighted eno schemes.
\newblock {\em Journal of computational physics}, 126(1):202--228, 1996.

\bibitem{kim2020fast}
Y.~Kim, Y.~Choi, D.~Widemann, and T.~Zohdi.
\newblock {A fast and accurate physics-informed neural network reduced order
  model with shallow masked autoencoder}.
\newblock {\em arXiv:2009.11990}, 2020.

\bibitem{kirby1992reconstructing}
M.~Kirby and D.~Armbruster.
\newblock {Reconstructing phase space from PDE simulations}.
\newblock {\em Zeitschrift f{\"u}r angewandte Mathematik und Physik ZAMP},
  43(6):999--1022, 1992.

\bibitem{lee2019deep}
K.~Lee and K.~Carlberg.
\newblock {Deep Conservation: A latent dynamics model for exact satisfaction of
  physical conservation laws}.
\newblock {\em arXiv:1909.09754}, 2019.

\bibitem{lee2020model}
K.~Lee and K.~Carlberg.
\newblock {Model reduction of dynamical systems on nonlinear manifolds using
  deep convolutional autoencoders}.
\newblock {\em Journal of Computational Physics}, 404:108973, 2020.

\bibitem{lee2021partition}
Kookjin Lee, Nathaniel~A Trask, Ravi~G Patel, Mamikon~A Gulian, and Eric~C Cyr.
\newblock {Partition of unity networks: deep hp-approximation}.
\newblock {\em arXiv:2101.11256}, 2021.

\bibitem{lu2020lagrangian}
Hannah Lu and Daniel~M Tartakovsky.
\newblock {Lagrangian dynamic mode decomposition for construction of
  reduced-order models of advection-dominated phenomena}.
\newblock {\em Journal of Computational Physics}, 407:109229, 2020.

\bibitem{lu2021dynamic}
Hannah Lu and Daniel~M Tartakovsky.
\newblock {Dynamic Mode Decomposition for Construction of Reduced-Order Models
  of Hyperbolic Problems with Shocks}.
\newblock {\em Journal of Machine Learning for Modeling and Computing}, 2(1),
  2021.

\bibitem{lu2019deeponet}
Lu~Lu, Pengzhan Jin, and George~Em Karniadakis.
\newblock {Deeponet: Learning nonlinear operators for identifying differential
  equations based on the universal approximation theorem of operators}.
\newblock {\em arXiv:1910.03193}, 2019.

\bibitem{maday2013locally}
Yvon Maday and Benjamin Stamm.
\newblock Locally adaptive greedy approximations for anisotropic parameter
  reduced basis spaces.
\newblock {\em SIAM Journal on Scientific Computing}, 35(6):A2417--A2441, 2013.

\bibitem{mccann1997convexity}
R.~McCann.
\newblock {A convexity principle for interacting gases}.
\newblock {\em Advances in mathematics}, 128(1):153--179, 1997.

\bibitem{melenk2000n}
Jens~M Melenk.
\newblock {On $n$-widths for elliptic problems}.
\newblock {\em Journal of mathematical analysis and applications},
  247(1):272--289, 2000.

\bibitem{mojgani2017lagrangian}
Rambod Mojgani and Maciej Balajewicz.
\newblock {Lagrangian basis method for dimensionality reduction of convection
  dominated nonlinear flows}.
\newblock {\em arXiv:1701.04343}, 2017.

\bibitem{mojgani2020physics}
Rambod Mojgani and Maciej Balajewicz.
\newblock {Physics-aware registration based auto-encoder for convection
  dominated PDEs}.
\newblock {\em arXiv preprint arXiv:2006.15655}, 2020.

\bibitem{nair2019transported}
N.~Nair and M.~Balajewicz.
\newblock Transported snapshot model order reduction approach for parametric,
  steady-state fluid flows containing parameter-dependent shocks.
\newblock {\em International Journal for Numerical Methods in Engineering},
  117(12):1234--1262, 2019.

\bibitem{nonino2019overcoming}
M.~Nonino, F.~Ballarin, G.~Rozza, and Y.~Maday.
\newblock {Overcoming slowly decaying Kolmogorov $n$-width by transport maps:
  application to model order reduction of fluid dynamics and fluid--structure
  interaction problems}.
\newblock {\em arXiv:1911.06598}, 2019.

\bibitem{ohlberger2013nonlinear}
M.~Ohlberger and S.~Rave.
\newblock {Nonlinear reduced basis approximation of parameterized evolution
  equations via the method of freezing}.
\newblock {\em Comptes Rendus Mathematique}, 351(23-24):901--906, 2013.

\bibitem{ohlberger2016reduced}
Mario Ohlberger and Stephan Rave.
\newblock {Reduced basis methods: success, limitations and future challenges}.
\newblock In {\em Proceedings of ALGORITMY}, pages 1--12, 2016.

\bibitem{peherstorfer2020model}
B.~Peherstorfer.
\newblock {Model reduction for transport-dominated problems via online adaptive
  bases and adaptive sampling}.
\newblock {\em SIAM Journal on Scientific Computing}, 42(5):A2803--A2836, 2020.

\bibitem{pinkus2012n}
Allan Pinkus.
\newblock {\em {N-widths in Approximation Theory}}, volume~7.
\newblock Springer Science \& Business Media, 2012.

\bibitem{reiss2018shifted}
Julius Reiss, Philipp Schulze, J{\"o}rn Sesterhenn, and Volker Mehrmann.
\newblock {The shifted proper orthogonal decomposition: A mode decomposition
  for multiple transport phenomena}.
\newblock {\em SIAM Journal on Scientific Computing}, 40(3):A1322--A1344, 2018.

\bibitem{rim2018transport}
D.~Rim, S.~Moe, and R.~LeVeque.
\newblock {Transport reversal for model reduction of hyperbolic partial
  differential equations}.
\newblock {\em SIAM/ASA Journal on Uncertainty Quantification}, 6(1):118--150,
  2018.

\bibitem{rim2019manifold}
D.~Rim, B.~Peherstorfer, and K.~Mandli.
\newblock {Manifold approximations via transported subspaces: model reduction
  for transport-dominated problems}.
\newblock {\em arXiv:1912.13024}, 2019.

\bibitem{rim2020depth}
D.~Rim, L.~Venturi, J.~Bruna, and B.~Peherstorfer.
\newblock {Depth separation for reduced deep networks in nonlinear model
  reduction: distilling shock waves in nonlinear hyperbolic problems}.
\newblock {\em arXiv:2007.13977}, 2020.

\bibitem{romor2022non}
Francesco Romor, Giovanni Stabile, and Gianluigi Rozza.
\newblock {Non-linear manifold ROM with Convolutional Autoencoders and Reduced
  Over-Collocation method}.
\newblock {\em arXiv preprint arXiv:2203.00360}, 2022.

\bibitem{rowley2003reduction}
C.~Rowley, I.~Kevrekidis, J.~Marsden, and K.~Lust.
\newblock {Reduction and reconstruction for self-similar dynamical systems}.
\newblock {\em Nonlinearity}, 16(4):1257, 2003.

\bibitem{rowley2000reconstruction}
C.~Rowley and J.~Marsden.
\newblock {Reconstruction equations and the Karhunen--Lo{\`e}ve expansion for
  systems with symmetry}.
\newblock {\em Physica D: Nonlinear Phenomena}, 142(1-2):1--19, 2000.

\bibitem{schoberl2014c++}
Joachim Sch{\"o}berl.
\newblock C++ 11 implementation of finite elements in ngsolve.
\newblock {\em Institute for Analysis and Scientific Computing, Vienna
  University of Technology}, 30, 2014.

\bibitem{shu1988efficient}
Chi-Wang Shu and Stanley Osher.
\newblock Efficient implementation of essentially non-oscillatory
  shock-capturing schemes.
\newblock {\em Journal of computational physics}, 77(2):439--471, 1988.

\bibitem{taddei2020registration}
T.~Taddei.
\newblock A registration method for model order reduction: data compression and
  geometry reduction.
\newblock {\em SIAM Journal on Scientific Computing}, 42(2):A997--A1027, 2020.

\bibitem{taddei2015reduced}
T.~Taddei, S.~Perotto, and A.~Quarteroni.
\newblock Reduced basis techniques for nonlinear conservation laws.
\newblock {\em ESAIM: Mathematical Modelling and Numerical Analysis},
  49(3):787--814, 2015.

\bibitem{tang2003adaptive}
Huazhong Tang and Tao Tang.
\newblock {Adaptive mesh methods for one-and two-dimensional hyperbolic
  conservation laws}.
\newblock {\em SIAM Journal on Numerical Analysis}, 41(2):487--515, 2003.

\bibitem{torlo2020model}
Davide Torlo.
\newblock {Model reduction for advection dominated hyperbolic problems in an
  ALE framework: Offline and online phases}.
\newblock {\em arXiv preprint arXiv:2003.13735}, 2020.

\bibitem{villani2008optimal}
C.~Villani.
\newblock {\em {Optimal transport: old and new}}, volume 338.
\newblock Springer Science \& Business Media, 2008.

\bibitem{wang2020reduced}
Min Wang, Siu~Wun Cheung, Wing~Tat Leung, Eric~T Chung, Yalchin Efendiev, and
  Mary Wheeler.
\newblock {Reduced-order deep learning for flow dynamics. The interplay between
  deep learning and model reduction}.
\newblock {\em Journal of Computational Physics}, 401:108939, 2020.

\bibitem{welper2017h}
G.~Welper.
\newblock {$ h $ and $ hp $-adaptive interpolation by transformed snapshots for
  parametric and stochastic hyperbolic PDEs}.
\newblock {\em arXiv:1710.11481}, 2017.

\bibitem{welper2017interpolation}
G.~Welper.
\newblock {Interpolation of functions with parameter dependent jumps by
  transformed snapshots}.
\newblock {\em SIAM Journal on Scientific Computing}, 39(4):A1225--A1250, 2017.

\bibitem{welper2020transformed}
G~Welper.
\newblock {Transformed snapshot interpolation with high resolution transforms}.
\newblock {\em SIAM Journal on Scientific Computing}, 42(4):A2037--A2061, 2020.

\bibitem{wheeler1978elliptic}
Mary~Fanett Wheeler.
\newblock An elliptic collocation-finite element method with interior
  penalties.
\newblock {\em SIAM Journal on Numerical Analysis}, 15(1):152--161, 1978.

\bibitem{yano2019discontinuous}
Masayuki Yano.
\newblock {Discontinuous Galerkin reduced basis empirical quadrature procedure
  for model reduction of parametrized nonlinear conservation laws}.
\newblock {\em Advances in Computational Mathematics}, 45(5):2287--2320, 2019.

\bibitem{yano2020goal}
Masayuki Yano.
\newblock {Goal-oriented model reduction of parametrized nonlinear partial
  differential equations: Application to aerodynamics}.
\newblock {\em International Journal for Numerical Methods in Engineering},
  121(23):5200--5226, 2020.

\end{thebibliography}

\end{document}